\documentclass[a4paper, 11pt, twoside, final]{scrartcl}

\usepackage{etex}
\usepackage{calc}
\usepackage{ifthen}
\usepackage{xargs}
\usepackage{remreset}

\usepackage[T1]{fontenc}
\usepackage[utf8]{inputenc}
\usepackage[ngerman,british]{babel}
\usepackage{eurosym}
\DeclareInputText{128}{\EUR{}}

\usepackage[intlimits]{amsmath}
\usepackage{amssymb}
\usepackage{mathtools}
\usepackage{bm}
\usepackage[framed, hyperref, thref, amsmath]{ntheorem}
\usepackage{units}
\usepackage{mleftright}

\usepackage[svgnames]{xcolor}
\usepackage[ntheorem, framemethod=ps]{mdframed}

\allowdisplaybreaks

\usepackage[p]{libertine}
\usepackage[libertine,slantedGreek]{newtxmath} 
\usepackage[bb=px,bbscaled=0.94,cal=dutchcal,calscaled=0.95]{mathalfa}
\useosf

\usepackage{anyfontsize}

\usepackage[singlespacing]{setspace}

\usepackage[headsepline,footsepline,plainfootsepline]{scrpage2}

\usepackage[shortlabels]{enumitem}
\setlist{listparindent=\parindent,parsep=\parskip}

\usepackage{graphicx}

\usepackage{placeins}

\usepackage{caption}
\DeclareCaptionFont{gray}{\color{black}}
\usepackage{subfig}
\usepackage{pstricks-add}
\usepackage{pst-plot}
\usepackage{pst-func}
\psset{unit=1cm}

\usepackage{array}
\usepackage{colortbl}
\usepackage{longtable}
\usepackage{ltablex}
\usepackage{varioref}

\usepackage{hyperref}
\usepackage[quadpoints=false]{hypdvips}

\hypersetup{hidelinks}
\hypersetup{breaklinks=true}
\hypersetup{bookmarksopen,bookmarksopenlevel=1}

\addto\extrasbritish{%
}

\usepackage{breakurl}
\usepackage{bookmark}

\KOMAoption{BCOR}{5mm}
\KOMAoption{DIV}{classic}
\KOMAoption{headinclude}{false}
\KOMAoption{footinclude}{false}
\KOMAoption{pagesize}{auto}

\recalctypearea

\KOMAoption{headings}{small}
\KOMAoption{numbers}{endperiod}

\newcommand{\AUTHORsave}{}
\newcommand{\AUTHOR}[1]{\renewcommand{\AUTHORsave}{#1}}

\newcommand{\SHORTAUTHORsave}{}
\newcommand{\SHORTAUTHOR}[1]{%
  \ifthenelse{\equal{#1}{}}%
  {\renewcommand{\SHORTAUTHORsave}{\AUTHORsave}}%
  {\renewcommand{\SHORTAUTHORsave}{#1}}%
}

\newcommand{\TITLEsave}{}
\newcommand{\TITLE}[1]{\renewcommand{\TITLEsave}{#1}}

\newcommand{\SHORTTITLEsave}{}
\newcommand{\SHORTTITLE}[1]{%
  \ifthenelse{\equal{#1}{}}%
  {\renewcommand{\SHORTTITLEsave}{\TITLEsave}}%
  {\renewcommand{\SHORTTITLEsave}{#1}}%
}

\newcommand{\INSTITUTEsave}{}
\newcommand{\INSTITUTE}[1]{\renewcommand{\INSTITUTEsave}{#1}}

\newcommand{\CORRESPONDENCEsave}{}
\newcommand{\CORRESPONDENCE}[1]{\renewcommand{\CORRESPONDENCEsave}{#1}}

\newcommand{\ABSTRACTsave}{}
\newcommand{\ABSTRACT}[1]{\renewcommand{\ABSTRACTsave}{#1}}

\newcommand{\KEYWORDSsave}{}
\newcommand{\KEYWORDS}[1]{\renewcommand{\KEYWORDSsave}{#1}}

\newcommand{\AMSCLASSsave}{}
\newcommand{\AMSCLASS}[1]{\renewcommand{\AMSCLASSsave}{#1}}

\newcommand{\Inst}[1]{\textsuperscript{#1}}
\newcommand{\email}[1]{\href{mailto:#1}{\nolinkurl{#1}}}

\newcommand{\TitleHeader}{%
  \noindent%
  \rule{\linewidth}{1pt}\\[10pt]%
  {\sffamily\bfseries\Large\TITLEsave}\\[10pt]%
  {\sffamily\bfseries\small\AUTHORsave}\\%
  {\sffamily\scriptsize\INSTITUTEsave}\\[15pt]%
  \parbox{\linewidth}{\sffamily\scriptsize {\bfseries
      Correspondence}\\\CORRESPONDENCEsave}\\[10pt]%
  \rule{\linewidth}{1pt}\\[15pt]%
  \hspace*{40pt}%
  \parbox{\linewidth-40pt}{%
    \small%
    {\sffamily\bfseries Abstract.} \quad \ABSTRACTsave%
  }\\[5pt]%
  \hspace*{40pt}%
  \parbox{\linewidth-40pt}{%
    \small%
    {\sffamily\bfseries Keywords.} \quad \KEYWORDSsave%
  }\\[5pt]%
  \hspace*{40pt}%
  \parbox{\linewidth-40pt}{%
    \small%
    {\sffamily\bfseries AMS subject classification.} \quad \AMSCLASSsave%
  }\\[5pt]%
  \rule{\linewidth}{1pt}%
}

\clearscrheadfoot

\setheadsepline{1pt}
\setfootsepline{1pt}

\ohead{\pagemark}
\cehead{\footnotesize\itshape\scshape \SHORTAUTHORsave}
\cohead{\footnotesize\itshape\scshape \SHORTTITLEsave}

\newcommand{\qed}{}

\newlength{\frametopsep}
\colorlet{ThmFrame}{CornflowerBlue!50}
\colorlet{ThmBody}{CornflowerBlue!15}

\mdfdefinestyle{thmframe}{innerleftmargin=5pt, innerrightmargin=5pt,
  innermargin=-5.5pt, outermargin=-5.5pt, skipabove=\frametopsep,
  skipbelow=\frametopsep, 
  roundcorner=5pt, linewidth=0.5pt, linecolor=ThmFrame,
  backgroundcolor=ThmBody, splittopskip=12.5pt, splitbottomskip=5pt}

\makeatletter
\newtheoremstyle{plainheaderbreak}
  {\item[]{\theorem@headerfont ##1\ ##2\theorem@separator}\hskip\labelsep}%
  {\item[]{\theorem@headerfont ##1\ ##2\
    (##3)}\theorem@separator\hskip\labelsep}
\makeatother

\theoremheaderfont{\sffamily\bfseries}
\theorembodyfont{}
\theoremstyle{plainheaderbreak}
\theoremseparator{.}
\theoremprework{}
\theorempostwork{}
\theoremframepreskip{\frametopsep}
\theoremframepostskip{\frametopsep}
\theoreminframepreskip{0pt}
\theoreminframepostskip{0pt}
\newmdtheoremenv[style=thmframe]{Definition}{Definition}[section]

\newmdtheoremenv[style=thmframe]{Problem}[Definition]{Problem}

\newmdtheoremenv[style=thmframe]{Assumption}[Definition]{Assumption}

\theoremheaderfont{\sffamily\bfseries\upshape}
\theorembodyfont{\itshape}

\newmdtheoremenv[style=thmframe]{Theorem}[Definition]{Theorem}

\newmdtheoremenv[style=thmframe]{Proposition}[Definition]{Proposition}

\newmdtheoremenv[style=thmframe]{Lemma}[Definition]{Lemma}

\newmdtheoremenv[style=thmframe]{Corollary}[Definition]{Corollary}

\theoremheaderfont{\sffamily\bfseries\upshape}
\theorembodyfont{\sffamily}
\theoremstyle{break}

\newmdtheoremenv[style=thmframe]{Algorithm}[Definition]{Algorithm}

\newlength{\unframetopsep}
\theoremheaderfont{\sffamily\bfseries}
\theorembodyfont{}
\theoremstyle{plain}
\theoremseparator{.}
\theorempreskip{\unframetopsep}
\theorempostskip{\unframetopsep}
\theoremprework{%
  \renewcommand{\qed}{%
    \hspace*{\fill}\nolinebreak[3]%
    \nopagebreak[3]\hspace*{\fill}{\raisebox{0.25pt}{\scalebox{0.9}{$\medcirc$}}}}}
\theorempostwork{}

\newtheorem{Example}[Definition]{Example}

\theoremprework{%
  \renewcommand{\qed}{%
    \hspace*{\fill}\nolinebreak[3]%
    \nopagebreak[3]\hspace*{\fill}{\raisebox{0.25pt}{\scalebox{0.9}{$\medcirc$}}}}}
\theorempostwork{}

\newtheorem{Counterexample}[Definition]{Counterexample}

\theoremprework{%
  \renewcommand{\qed}{%
    \hspace*{\fill}\nolinebreak[3]%
    \nopagebreak[3]\hspace*{\fill}{\raisebox{0.25pt}{\scalebox{0.9}{$\medcirc$}}}}}
\theorempostwork{}

\newtheorem{Remark}[Definition]{Remark}

\makeatletter
\newtheoremstyle{nonumberplainof}%
{\item[\theorem@headerfont\hskip\labelsep ##1\theorem@separator]}%
{\item[\theorem@headerfont\hskip \labelsep ##1\ of\ ##3\theorem@separator]}
\makeatother
\theoremstyle{nonumberplainof}

\theoremprework{%
    \renewcommand{\qed}{%
        \hspace*{\fill}\nolinebreak[3]%
        \nopagebreak[3]\hspace*{\fill}{$\square$}}}
\theorempostwork{}

\newtheorem{Proof}{Proof}

\newcommand{\pers}[1]{%
    {\scshape#1}%
}
\newcommand{\PBanach}{\pers{Banach}}

\newcommand{\PBregman}{\pers{Bregman}}
\newcommand{\PCaratheodory}{\pers{Carath\'eodory}}
\newcommand{\PCartes}{\pers{Cartes}}

\newcommand{\PCauchy}{\pers{Cauchy}}

\newcommand{\PEuklid}{\pers{Euklid}}

\newcommand{\PFermat}{\pers{Fermat}}
\newcommand{\PFredholm}{\pers{Fredholm}}

\newcommand{\PFourier}{\pers{Fourier}}
\newcommand{\PFrobenius}{\pers{Frobenius}}
\newcommand{\PGateaux}{\pers{Gâ\-teaux}}
\newcommand{\PGauss}{\pers{Gauß}}

\newcommand{\PHankel}{\pers{Hankel}}
\newcommand{\PHausdorff}{\pers{Hausdorff}}

\newcommand{\PHilbert}{\pers{Hilbert}}

\newcommand{\PJensen}{\pers{Jensen}}

\newcommand{\PNewton}{\pers{Newton}}

\newcommand{\PPaley}{\pers{Paley}}

\newcommand{\PSchauder}{\pers{Schauder}}
\newcommand{\PSchmidt}{\pers{Schmidt}}
\newcommand{\PSchwarz}{\pers{Schwarz}}

\newcommand{\PSobolew}{\pers{Sobolew}}

\newcommand{\PTaylor}{\pers{Tay\-lor}}
\newcommand{\PTikhonov}{\pers{Tik\-ho\-nov}}
\newcommand{\PTitchmarsh}{\pers{Titchmarsh}}
\newcommand{\PToeplitz}{\pers{Toeplitz}}

\newcommand{\PWiener}{\pers{Wiener}}

\newcommand{\ie}{%
    i.e.%
}
\newcommand{\etal}{%
    et~al.%
}
\newcommand{\cf}{%
    cf.%
}

\newcommand{\R}{%
    \ensuremath{\mathbb{R}}%
}

\newcommand{\N}{%
    \ensuremath{\mathbb{N}}%
}

\newcommand{\C}{%
    \ensuremath{\mathbb{C}}%
}
\newcommand{\iProd}[2]{%
    \left\langle#1,#2\right\rangle%
}

\newcommand{\iProdn}[2]{%
    \langle#1,#2\rangle%
}
\newcommand{\iProdb}[2]{%
    \bigl\langle#1,#2\bigr\rangle%
}

\newcommand{\absn}[1]{%
    |\hspace{1pt}#1\hspace{1pt}|%
}

\newcommand{\absbb}[1]{%
    \biggl|\hspace{1pt}#1\hspace{1pt}\biggr|%
}

\newcommand{\pNorm}[1]{%
    \left|\left|\hspace{1pt}#1\hspace{1pt}\right|\right|%
}

\newcommand{\pNormn}[1]{%
    ||\hspace{1pt}#1\hspace{1pt}||%
}

\newcommand{\Fourier}{%
    \mathop{\kern0pt\mathcal{F}}\nolimits%
}

\newcommand{\Fresnel}{%
    \mathop{\kern0pt\mathcal{E}}\nolimits%
}

\newcommand{\Laplace}{%
    \mathop{\kern0pt\mathcal{L}}\nolimits%
}
\newcommand{\LCT}{%
    \mathop{\kern0pt\mathcal{C}}\nolimits%
}

\newcommand{\z}{%
    \mathop{\kern0pt\mathcal{Z}}\nolimits%
}
\newcommand{\Landau}{%
    \mathop{\kern0pt\mathcal{O}}\nolimits%
}
\newcommand{\Ind}{%
    \vvmathbb{1}%
}
\newcommand{\e}{%
    \ensuremath{\mathrm{e}}%
}
\newcommand{\I}{%
    \ensuremath{\mathrm{i}}%
}

\newcommand{\argmin}{%
    \ensuremath{\operatorname*{argmin}}%
}
\newcommand{\conv}{%
    \ensuremath{\operatorname{conv}}%
}

\newcommand{\dom}{%
  \ensuremath{\operatorname{dom}}%
}

\newcommand{\Id}{%
    \ensuremath{\operatorname{Id}}%
}

  \newcommand{\minimize}{%
    \ensuremath{\operatorname*{minimize}}%
}

\newcommand{\ran}{%
    \ensuremath{\operatorname{ran}}%
}

\newcommand{\sinc}{%
    \ensuremath{\operatorname{sinc}}%
}

\newcommand{\Span}{%
    \ensuremath{\operatorname{span}}%
}

\newcommand{\supp}{%
    \ensuremath{\operatorname{supp}}%
}

\newcommand{\weakrightharpoonup}{%
  \xrightharpoonup{\raisebox{-4pt}[0pt][0pt]{$*$}}%
}

\newcommand{\diff}{%
    \mathop{}\!\mathrm{d}%
}

\newcommand{\T}{%
    \ensuremath{\mathrm{T}}%
}

\newcommand{\addmathskip}[1][5pt]{%
    \vspace*{#1}%
}
\newcommand{\submathskip}[1][-5pt]{%
    \vspace*{#1}%
}

\newlength{\subalignskip}
\newlength{\subalignaboveskip}
\newlength{\subalignbelowskip}
\newlength{\fsmallskip}
\newlength{\fskip}
\newlength{\subintertextskip}
\newlength{\addintertextskip}
\newcommand{\bq}[1]{%
    `#1'%
}

\AUTHOR{Robert \pers{Beinert}\Inst{1} and Kristian
  \pers{Bredies}\Inst{1}}

\SHORTAUTHOR{R. Beinert and K. Bredies}

\TITLE{Non-convex regularization of bilinear and quadratic inverse
  problems by tensorial lifting}

\SHORTTITLE{Regularization of bilinear and quadratic inverse problems}

\INSTITUTE{\Inst{1}\parbox[t]{0.95\linewidth}{Institut für Mathematik
    und Wissenschaftliches Rechnen\\
    \pers{Karl}-\pers{Franzens}-Universität Graz\\
    Heinrichstraße 36\\
    8010 Graz, Austria
  }
}

\CORRESPONDENCE{R. \pers{Beinert}: \email{robert.beinert@uni-graz.at}\\
  K. \pers{Bredies}: \email{kristian.bredies@uni-graz.at}}

\ABSTRACT{Considering the question: how non-linear may a non-linear
  operator be in order to extend the linear regularization theory, we
  introduce the class of dilinear mappings, which covers linear,
  bilinear, and quadratic operators between \PBanach\ spaces.  The
  corresponding dilinear inverse problems cover blind deconvolution,
  deautoconvolution, parallel imaging in MRI, and the phase retrieval
  problem.  Based on the universal property of the tensor product, the
  central idea is here to lift the non-linear mappings to linear
  representatives on a suitable topological tensor space.  At the same
  time, we extend the class of usually convex regularization functionals
  to the class of diconvex functionals, which are likewise defined by
  a tensorial lifting.  Generalizing the concepts of subgradients and
  \PBregman\ distances from convex analysis to the new framework, we
  analyse the novel class of dilinear inverse problems and establish
  convergence rates under similar conditions than in the linear
  setting.  Considering the deautoconvolution problem as specific
  application, we derive satisfiable source conditions and validate
  the theoretical convergence rates numerically.}


\AMSCLASS{45Q05, 46A32, 47J06, 65J15}

\begin{document}
\thispagestyle{empty}
\TitleHeader

\section{Introduction}
\label{sec:introduction}

Nowadays the theory of inverse problems has become one of the central
mathematical approaches to solve recovery problems in medicine,
engineering, and life sciences.  Some of the main applications are
computed tomography (CT), magnetic resonance imaging (MRI), and
deconvolution problems in microscopy, see for instance \cite{BB98,
  MS12, Ram05, SW13, Uhl03, Uhl13} being recent monographs as well as
many other publications.

The beginnings of the modern regularization theory for ill-posed
problems are tracing back to the pioneering works of A.~N.~\PTikhonov\
\cite{Tik63b, Tik63a}.  Between then and now, the theory has been
heavily extended and covers linear and non-linear formulations in the
\PHilbert\ space as well as in the \PBanach\ space setting.  In order
to name at least a few of the numerous monographs, we refer to
\cite{BG94, EHN96, Hof86, Lou89, LP13, Mor84, SKHK12, TA77,TLY98}.
Due to the enormous relevance of the research topic, the published
literature embraces many further monographs as well as a vast number
of research articles.

Especially for linear problem formulations, the analytical framework
is highly developed and allows the general treatment of inverse
problems with continuous operators, see for instance \cite{EHN96,
  MS12, SKHK12} and references therein.  In addition to the
sophisticated analysis, efficient numerical implementations of the
solution schemes are available for practitioners \cite{EHN96, Sch11}.
The interaction between analysis and numerics are one reason for the
great, interdisciplinary success of the linear regularization theory
for ill-posed inverse problems.

If the linear operator in the problem formulation is replaced by a
non-linear operator, the situation changes dramatically.  Depending on
the operator, there are several regularization approaches with
different benefits and drawbacks \cite{EHN96, Gra11, HKPS07, SGG+09}.  One
standard approach to regularize non-linear operators is to introduce
suitable non-linearity conditions and to restrict the set of
considered operators.  These constraints are mostly based on
properties of the remainder of the first-order \PTaylor\ expansion
\cite{BS04, EKN89, HKPS07, RS06}.  In an abstract way, this approach
allows the generalization of well-understood linear results by
controlling the deviation from the linear setting.  Unfortunately, the
validation of the required assumptions for a specific non-linear
operator is a non-trivial task.

In order to extend the linear theory to the non-linear domain further,
our idea is to introduce a class of operators that covers many
interesting applications for practitioners and, at the same time,
allows a general treatment of the corresponding inverse problems.
More precisely, we introduce the class of dilinear operators that
embraces linear, bilinear, and quadratic mappings between \PBanach\
spaces.  Consequently, our novel class of dilinear inverse problems
covers formulations arising in imaging and physics \cite{SGG+09} like
blind deconvolution \cite{BS01, JR06}, deautoconvolution \cite{GH94,
  FH96, GHB+14, ABHS16}, parallel imaging in MRI \cite{BBM+}, or phase
retrieval \cite{DF87, Mil90, SSD+06}.

The central idea behind the class of dilinear operators is the
universal property of the topological tensor product, which enables us
to lift a continuous but non-linear mapping to a linear operator.
Owing to the lifting, we get immediate access to the linear
regularization theory.  On the downside, a simple lifting of the
non-linear inverse problem causes an additional non-convex rank-one
constraint, which is similarly challenging to handle than the original
non-linear problem.  For this reason, most results of the linear
regularization theory are not applicable for the lifted problem and
cannot be transferred to the original (unlifted) inverse problem.  In
order to overcame this issues, we use the tensorial lifting indirectly
and generalize the required concepts from convex analysis to the new
framework.

The recent literature already contains some further ideas to handle
inverse problems arising from bilinear or quadratic operators.  For
instance, each quadratic mapping between separable \PHilbert\ spaces
may be factorized into a linear operator and a strong quadratic
isometry so that the corresponding inverse problem can be decomposed
into a possibly ill-posed linear and a well-posed quadratic part, see
\cite{Fle14}.  In order to determine a solution, one can now apply a
two-step method.  Firstly, the ill-posed linear part is solved by a
linear regularization method.  Secondly, the well-posed quadratic
problem is solved by projection onto the manifold of symmetric
rank-one tensors.  The main drawback of this approach is that the
solution of the linear part has not to lie in the range of the
well-posed quadratic operator such that the second step may corrupt
the obtained solution.  This issue does not occur if the forward
operator of the linear part is injective, which is generally not
true for quadratic inverse problems.

Besides the forward operator, one can also generalize the used
regularization from usually convex to non-convex functionals.  An
abstract analysis of non-convex regularization methods for bounded
linear operators between \PHilbert\ spaces has been introduced in
\cite{Gra10}, where the definition of the subdifferential and the
\PBregman\ distance have been extended with respect to an arbitrary
set of functions.  On the basis of a variational source condition, one
can further obtain convergence rates for these non-convex
regularization methods.  Similarly to \cite{Gra10}, we employ
non-convex regularizations -- however -- with the tensorial
lifting in mind.

As mentioned above, the deautoconvolution problem is one specific
instance of a dilinear inverse problem \cite{GHB+14, ABHS16}.
Although the unregularized problem can have two different solutions at
the most, the deautoconvolution problem is everywhere locally ill
posed \cite{GH94, FH96, Ger11}.  Nevertheless, with an appropriate
regularization, very accurate numerical solutions can be obtained
\cite{CL05, ABHS16}.  Establishing theoretical convergence rates for
the applied regularization is, unfortunately, very challenging since
most conditions for the non-linear theory are not fulfilled
\cite{BH15a, ABHS16}.  For more specific classes of true solutions,
for instance, the class of all trigonometric polynomials or some
subset of \PSobolew\ spaces, however, the regularized solutions
converge to the true solution with a provable rate, see \cite{BFH16}
or \cite{Jan00, DL08} respectively.  Applying our novel regularization
theory, we establish convergence rates under a source-wise
representation of the subdifferential of the regularization
functional.  In other words, we generalize the classical range source
condition in a specific manner fitting the necessities of dilinear
inverse problems.

In this paper, we show that the essential results of the classical
regularization theory with bounded linear operators and convex
regularization functionals may be extended to bilinear and quadratic
forward operators.  At the same time, we allow the regularization
functional to be non-convex in a manner being comparable with the
non-linearity of the considered operator.

Since our analysis is mainly based on the properties of the
topological tensor product \cite{DF93, Rya02}, we firstly give a brief
survey of tensor spaces and the tensorial lifting in
\autoref{sec:tensor-products}.  Our main foci are here the different
interpretations of a specific tensor.  Further, we introduce the set
of dilinear operators and show that each dilinear operator may be
uniquely lifted to a linear operator.  Analogously, in
\autoref{sec:gener-subgr}, the class of diconvex functionals is
defined through a convex lifting.  In order to study diconvex
regularization methods, we generalize the usually convex
subdifferential and \PBregman\ distance with respect to dilinear
mappings.  Further, we derive sum and chain rules for the new dilinear
subdifferential calculus.

Our main results about dilinear -- bilinear and quadratic -- inverse
problems are presented in \autoref{sec:dil-inverse-probl}.  Under
suitable assumptions being comparable to the assumptions in the
classical linear theory, dilinear inverse problems are well posed,
stable, and consistent.  Using a source-wise representation as in the
classical range condition, we obtain convergence rates for the data
fidelity term and the \PBregman\ distance between the true solution of
the undisturbed and the regularized problem.

In \autoref{sec:deaut-probl}, we apply our non-convex regularization
theory for dilinear inverse problems to the deautoconvolution problem and
study the required assumptions in more detail.  Further, we
reformulate the source-wise representation and obtain an equivalent
source condition, which allows us to construct suitable true
solutions.  With the numerical experiments in \autoref{sec:numer-simul}, we
verify the derived rates for the deautoconvolution numerically and
give some examples of possible signals fulfilling the source
condition.

\section{Tensor products and dilinear mappings}
\label{sec:tensor-products}

The calculus of tensors has been invented more than a century ago.
Since then tensors have been extensively studied from an algebraic and
topological point of view, see for instance
\cite{DF93,DFS08,Lic66,Rom08,Rya02} and references therein.  One of
the most remarkable results in tensor analysis, which is the starting
point for our study of non-linear inverse problems, is the universal
property of the tensor product that allows us to lift a bilinear
mapping to a linear one.  In order to get some intuition about tensor
calculus on \PBanach\ spaces and about the different interpretations
of a specific tensor, we briefly survey the required central ideas
about tensor products and adapt the lifting approach to the class of
dilinear forward operators.

The tensor product of two real vector spaces $V_1$ and $V_2$ denoted
by $V_1 \otimes V_2$ can be constructed in various ways.  Following
the presentations of \pers{Ryan} \cite{Rya02} and \pers{Diestel}
\etal\ \cite{DFS08}, we initially define a tensor as a linear operator
acting on bilinear forms.  For this purpose, we recall that a mapping
from $V_1 \times V_2$ into the vector space $W$ is \emph{bilinear} if
it is linear with respect to each variable.  The corresponding vector
space of all bilinear mappings is denoted by $B(V_1 \times V_2, W)$.
In the special case $W = \R$, we simply write $B(V_1 \times V_2)$ for
the vector space of all \emph{bilinear forms}.  Given two elements
$u \in V_1$ and $v \in V_2$, we now define the \emph{tensor}
$u \otimes v$ as the linear functional that evaluates a bilinear form
at the point $(u,v)$.  The tensor $u\otimes v$ thus acts on a bilinear
form $A \in B(V_1 \times V_2)$ by
\begin{equation*}
  (u \otimes v)(A) = A(u,v).
\end{equation*}

On the basis of this construction, the \emph{tensor product} $V_1
\otimes V_2$ of the real vector spaces $V_1$ and $V_2$ now consists of
all finite linear combinations
\begin{equation*}
  w = \sum_{n=1}^N \lambda_n \, u_n \otimes v_n
\end{equation*}
with $u_n \in V_1$, $v_n \in V_2$, $\lambda_n \in \R$, and $N \in \N$.
Thus, the tensor product $V_1 \otimes V_2$ is the subspace of the
algebraic dual $B(V_1 \times V_2)^\#$ spanned by the rank-one tensors
$u \otimes v$ with $u \in V_1$ and $v \in V_2$.  Using the definition
of a tensor as linear functional on $B(V_1 \times V_2)$, one can
easily verify that the mapping $(u,v) \mapsto u \otimes v$ itself is
bilinear.  In other words, we have
\begin{enumerate}[(i)]
\item $(u_1 + u_2) \otimes v = u_1 \otimes v + u_2 \otimes v$,
\item $u \otimes (v_1 + v_2) = u \otimes v_1 + u \otimes v_2$,
\item $\lambda (u \otimes v) = (\lambda u) \otimes v = u \otimes
  (\lambda v)$, and
\item $0 \otimes v = u \otimes 0 = 0$.
\end{enumerate}

One of the most central concepts behind the calculus of tensors and
most crucial for our analysis of non-linear inverse problems is the
linearization of bilinear mappings.  Given a bilinear mapping
$A \colon V_1 \times V_2 \rightarrow W$, we define the linear mapping
$\breve A \colon V_1 \otimes V_2 \rightarrow W$ by
$\breve A (\sum_{n=1}^N u_n \otimes v_n) = \sum_{n=1}^N A(u_n, v_n)$.
One can show that the mapping $\breve A$ is well defined and is the
only linear mapping in $L(V_1 \otimes V_2, W)$ such that
$A(u,v) = \breve A(u \otimes v)$ for all $u \in V_1$ and $v \in V_2$,
see for instance \cite{Rya02}.

\begin{Proposition}[Lifting of bilinear mappings]
  \label{prop:lift-bil-map}
  Let $A \colon V_1 \times V_2 \rightarrow W$ be a bilinear mapping,
  where $V_1$, $V_2$, and $W$ denote real vector spaces.  Then there exists a
  unique linear mapping $\breve A \colon V_1 \otimes V_2 \rightarrow
  W$ so that $A(u,v) = \breve A(u \otimes v)$ for all $(u,v)$ in $V_1
  \times V_2$.
\end{Proposition}

Besides the definition of a tensor as linear functional on
$B(V_1 \times V_2)$, we can interpret a tensor itself as a bilinear
form.  For this purpose, we associate to each $u \in V_1$ and
$v \in V_2$ the bilinear form
$B_{u \otimes v} \colon V_1^\# \times V_2^\# \rightarrow \R$ with
$B_{u \otimes v}(\phi,\psi) = \phi(u) \, \psi(v)$.  Moreover, the
mapping $(u,v) \mapsto B_{u \otimes v}$ is also bilinear, which implies that
there is a linear mapping from $V_1 \otimes V_2$ into
$B(V_1^\# \times V_2^\#)$.  Since this mapping is injective, see
\cite{Rya02}, each tensor $w$ in $V_1 \otimes V_2$ corresponds
uniquely to a bilinear mapping, and hence
$V_1 \otimes V_2 \subset B(V_1^\# \times V_2^\#)$.

Further, the tensor product $V_1 \otimes V_2$ can be seen as vector
space of linear mappings \cite{DFS08,Rya02}.  For this, we observe
that each bilinear mapping $A \in B(V_1 \times V_2)$ generates the
linear mappings $L_A \colon V_1 \rightarrow V_2^\#$ with
$u \mapsto A(u,\cdot)$ and $R_A \colon V_2 \rightarrow V_1^\#$ with
$v \mapsto A(\cdot, v)$ by fixing one of the components.  In this
context, we may write
\begin{equation}
  \label{eq:tensor-lin-op}
  A(u,v) = \iProd{L_A(u)}{v} = \iProd{R_A(v)}{u}.
\end{equation}
On the basis of this idea, each tensor
$w = \sum_{n=1}^N \lambda_n \, u_n \otimes v_n$ in $V_1 \otimes V_2$
generates the linear mappings
\begin{equation}
  \label{eq:left-right-map}
  L_w (\phi) \coloneqq \sum_{n=1}^N \lambda_n \, \phi(u_n) \, v_n
  \qquad\text{and}\qquad
  R_w (\psi) \coloneqq \sum_{n=1}^N \lambda_n \, \psi(v_n) \, u_n.
\end{equation}
In other words, we have $V_1 \otimes V_2 \subset L(V_1^\# , V_2)$ and
$V_1 \otimes V_2 \subset L(V_2^\#, V_1)$.  If one of the real vector spaces $V_1$
or $V_2$ is already a dual space, we have a natural embedding into the
smaller linear function spaces $V_1^\# \otimes V_2 \subset L(V_1,
V_2)$ and $V_1 \otimes V_2^\# \subset L(V_2,V_1)$, see \cite{Rya02}.

There are several approaches to define an appropriate norm on the
tensor product $V_1 \otimes V_2$.  Here we employ the projective norm,
which allow us to lift each bounded bilinear operator to a bounded
linear operator on the tensor product.  To define the projective norm,
we assume that the real vector spaces $V_1$ and $V_2$ above are
already equipped with an appropriate norm.  More precisely, we replace
the arbitrary real vector spaces $V_1$ and $V_2$ by two \PBanach\
spaces $X_1$ and $X_2$.  The projective norm is now defined in the
following manner, see for instance \cite{DF93,Rya02}.

\begin{Definition}[Projective norm]
  \label{def:proj-norm}
  Let $X_1$ and $X_2$ be real \PBanach\ spaces.  The \emph{projective norm}
  $\pNorm{\cdot}_\uppi$ on the tensor product $X_1 \otimes X_2$ is
  defined by
  \begin{equation*}
    \pNorm{w}_\uppi \coloneqq
    \inf \biggl\{ \sum_{n=1}^N \pNorm{u_n} \pNorm{v_n} : w =
    \sum_{n=1}^N u_n \otimes v_n \biggr\},
  \end{equation*}
  where the infimum is taken over all finite representations of
  $w \in X_1 \otimes X_2$.
\end{Definition}

The projective norm on $X_1 \otimes X_2$ belongs to the reasonable
crossnorms, which means that
$\pNorm{u \otimes v}_\uppi = \pNorm{u} \pNorm{v}$ for all $u \in X_1$
and $v \in X_2$, see \cite{DFS08,Rya02}.  Based on the projective
norm, we obtain the projective tensor product $X_1 \otimes_\uppi X_2$.

\begin{Definition}[Projective tensor product]
  \label{def:proj-tensor-prod}
  Let $X_1$ and $X_2$ be real \PBanach\ spaces.  The \emph{projective
    tensor product} $X_1  \otimes_\uppi X_2$ of $X_1$ and $X_2$ is the
  completion of the tensor product $X_1 \otimes X_2$ with respect to
  the projective norm $\pNorm{\cdot}_\uppi$.
\end{Definition}

Figuratively, we complete the tensor product $X_1 \otimes X_2$
consisting of all finite-rank tensors by the infinite-rank tensors $w$ with
\begin{equation*}
  \pNorm{w}_\uppi =
  \inf \biggl\{ \sum_{n=1}^\infty \pNorm{u_n} \pNorm{v_n} : w =
  \sum_{n=1}^\infty u_n \otimes v_n \biggr\} < +\infty.
\end{equation*}
Similarly as above, if one of the \PBanach\ spaces $X_1$ or $X_2$ is a
dual space, the projective tensor product can be embedded into the
space of bounded linear operators.  More precisely, we have
$X_1^* \otimes_\uppi X_2 \subset \mathcal L(X_1, X_2)$ and
$X_1 \otimes_\uppi X_2^* \subset \mathcal L(X_2,X_1)$, see for
instance \cite{Won79}. 
 In the \PHilbert\ space setting, the
projective tensor product $H_1 \otimes_\uppi H_2$ of the \PHilbert\
spaces $H_1$ and $H_2$ corresponds to the trace class operators, and
the projective norm is given by
$\pNorm{w}_\uppi = \sum_{n=1}^\infty \sigma_n(w)$ for all
$w \in H_1 \otimes_\uppi H_2$, where $\sigma_n(w)$ denotes the $n$th
singular value of $w$, see \cite{Wer02}.

The main benefit of equipping the tensor product with the projective
norm is that \thref{prop:lift-bil-map} remains valid for bounded
bilinear and linear operators.  For this, we recall that a bilinear
operator $A \colon X_1 \times X_2 \rightarrow Y$ is \emph{bounded} if
there exists a real constant $C$ such that
$\pNorm{A(u,v)} \le C \pNorm{u} \pNorm{v}$ for all $(u,v)$ in
$X_1 \times X_2$.  Following the notation in \cite{Rya02}, we denote
the \PBanach\ space of all bounded bilinear operators from
$X_1 \times X_2$ into $Y$ equipped with the norm
$\pNorm{A} \coloneqq
\sup\{\nicefrac{\pNorm{A(u,v)}}{\pNorm{u}\pNorm{v}} : u \in X_1
\setminus \{0\},
v \in X_2 \setminus\{0\} \}$
by $\mathcal B(X_1 \times X_2, Y)$.  In the special case $Y = \R$, we
again write $\mathcal B(X_1 \times X_2)$.  More precisely, the lifting
of bounded bilinear operators can be stated in the following form, see
for instance \cite{Rya02}.

\begin{Proposition}[Lifting of bounded bilinear mappings]
  \label{prop:lift-boun-bil-map}
  Let $A \colon X_1 \times X_2 \rightarrow Y$ be a bounded bilinear operator,
  where $X_1$, $X_2$, and $Y$ denote real \PBanach\ spaces.  Then
  there exists a unique bounded linear operator $\breve A \colon X_1
  \otimes_\uppi X_2 \rightarrow Y$ so that $A(u,v)=\breve A(u \otimes v)$
  for all $(u,v)$ in $X_1 \times X_2$.
\end{Proposition}

The lifting of a bilinear operator in \thref{prop:lift-boun-bil-map}
is also called the \emph{universal property} of the projective tensor
product.  Vice versa, each bounded linear mapping
$\breve A \colon X_1 \otimes_\uppi X_2 \rightarrow Y$ uniquely defines
a bounded bilinear mapping $A$ by $A(u,v) = \breve A(u \otimes v)$,
which gives the canonical identification
\begin{equation*}
  \mathcal B(X_1 \times X_2, Y) 
  = \mathcal L(X_1 \otimes_\uppi X_2, Y).
\end{equation*}
Consequently, the topological dual of the projective tensor product
$X_1 \otimes_\uppi X_2$ is the space $\mathcal B (X_1 \times X_2)$ of bounded
bilinear forms, where a specific bounded bilinear mapping
$A \colon X_1 \times X_2 \rightarrow \R$ acts on an arbitrary tensor
$w = \sum_{n=1}^\infty \lambda_n \, u_n \otimes v_n$ by
\begin{equation*}
  \iProd{A}{w} = \sum_{n=1}^\infty \lambda_n \, A(u_n,v_n),
  \submathskip
\end{equation*}
see for instance \cite{Rya02}.

In order to define the novel class of dilinear operators, we restrict
the projective tensor product to the subspace of symmetric tensors.
Assuming that $X$ is a real \PBanach\ space, we call a tensor
$w \in X \otimes_\uppi X$ \emph{symmetric} if and only if $w = w^\T$,
where the \emph{transpose} of
$w = \sum_{n=1}^\infty \lambda_n \, u_n \otimes v_n$ is given by
$w^\T \coloneqq \sum_{n=1}^\infty \lambda_n \, v_n \otimes u_n$.  In
this context, a tensor $w$ is symmetric if and only if the related
bilinear form $B_w$ defined by
$B_w(\phi,\psi) = \sum_{n=1}^\infty \lambda_n \, \phi(u_n) \,
\psi(v_n)$
is symmetric, since
$B_w(\phi,\psi) = B_{w^\T}(\phi,\psi) = B_w (\psi, \phi)$ for all
$\phi$ and $\psi$ in $X^*$.  In the following, the closed subspace of the
symmetric tensors spanned by $u \otimes u$ with $u \in X$ is denoted
by $X \otimes_{\uppi,\mathrm{sym}} X$.  Considering that the
\emph{annihilator}
\begin{equation*}
  (X \otimes_{\uppi, \mathrm{sym}} X)^\perp
  \coloneqq 
  \{ A \in \mathcal B (X \times X) : A(w) = 0 \ \text{for all}\  w
  \in X \otimes_{\uppi, \mathrm{sym}} X \} 
\end{equation*}
consists of all antisymmetric bilinear forms, which means
$A(u,v) = - A(v,u)$, the topological dual of
$X \otimes_{\uppi, \mathrm{sym}} X$ is isometric isomorphic to the
space $\mathcal B_{\mathrm{sym}}(X \times X)$ of symmetric, bounded
bilinear forms, see for instance \cite{Meg98}.

On the basis of these preliminary considerations, we are now ready to
define the class of dilinear operators, which embraces linear,
quadratic, and bilinear mappings as discussed below.

\begin{Definition}[Dilinear mappings]
  \label{def:dilinear-map}
  Let $X$ and $Y$ be real \PBanach\ spaces.  A mapping
  $K \colon X \rightarrow Y$ is \emph{dilinear} if there exists a
  linear mapping
  $\breve K \colon X \times (X \otimes_{\uppi, \mathrm{sym}} X)
  \rightarrow Y$
  such that $K(u) = \breve K (u,u\otimes u)$ for all $u$ in $X$.
\end{Definition}

Hence, the dilinear mappings are the restrictions of the linear
operators from $X \times (X \otimes_{\uppi,\mathrm{sym}} X)$ into $Y$
to the \emph{diagonal} $\{ (u, u \otimes u) : u \in X \}$.  Since the
\emph{representative} $\breve K$ acts on a \PCartes{}ian product, we
can always find two linear mappings $\breve A \colon X \rightarrow Y$
and $\breve B \colon X \otimes_{\uppi, \mathrm{sym}} X \rightarrow Y$
so that $\breve K (u,w) = \breve A(u) + \breve B(w)$.  A dilinear
mapping $K$ is \emph{bounded} if the representative linear mapping
$\breve K$ is bounded.  It is worth to mention that the representative
$\breve K$ of a bounded dilinear operator $K$ is uniquely defined.

\begin{Lemma}[Lifting of bounded dilinear mappings]
  \label{lem:lift-boun-dil-map}
  Let $K \colon X \rightarrow Y$ be a bounded dilinear mapping, where
  $X$ and $Y$ denotes real \PBanach\ spaces.  Then the (bounded)
  representative $\breve K$ is unique.
\end{Lemma}

\begin{Proof}
  Let us suppose that there exist two bounded representatives
  $\breve K_1$ and $\breve K_2$ for the bounded dilinear mapping $K$.
  Since both representatives $\breve K_\ell$ can be written as
  $\breve K_\ell(u) = \breve A_\ell (u) + \breve B_\ell (u \otimes
  u)$,
  where $A_\ell \colon X \rightarrow Y$ and
  $B_\ell \colon X \otimes_{\uppi,\mathrm{sym}} X \rightarrow Y$ are
  bounded linear mappings, we have
  \begin{equation*}
    \breve A_1 (u) - \breve A_2 (u)
    = \breve B_2 ( u \otimes u) - \breve B_1 (u \otimes u)
    \addmathskip
  \end{equation*}
  for all $u$ in $X$.  By replacing $u$ by $t u$ with $t \ge 0$ and
  $\pNorm{u}=1$, this identity is equivalent to
  \begin{equation*}
    t \, (\breve A_1 - \breve A_2)(u) 
    = t^2 \, (\breve B_2 - \breve B_1) (u \otimes u),
  \end{equation*}
  which already implies $\breve A_1 = \breve A_2$.  Due to the
  continuity and linearity of $\breve B_\ell$, the mappings
  $\breve B_\ell$ coincide on the symmetric subspace
  $X \otimes_{\uppi, \mathrm{sym}} X$, which yields the assertion.
\qed
\end{Proof}

\begin{Remark}[Unique dilinear lifting]
  \label{rem:uni-dil-lift}
  Because of the uniqueness in \thref{lem:lift-boun-dil-map}, the
  breve $\breve \cdot$ henceforth denotes the unique lifting of a
  bounded dilinear mapping.
  \qed
\end{Remark}

\begin{Example}[Linear mappings]
  \label{ex:lin-amp}
  One of the easiest examples of a dilinear operator are the linear
  operators $A \colon X \rightarrow Y$ with the representative
  $\breve A (u, u \otimes u) = A(u)$.  Consequently, the dilinear
  operators can be seen as a generalization of linear mappings.  \qed
\end{Example}

\begin{Example}[Quadratic mappings]
  \label{ex:quad-map}
  A further example of dilinear mappings are the quadratic mappings.
  For this, we recall that a mapping $Q \colon X \rightarrow Y$ is
  \emph{quadratic} if there exists a bounded symmetric bilinear
  mapping $A \colon X \times X \rightarrow Y$ such that
  $Q(u) = A(u,u)$ for all $u$ in $X$.  Since each bilinear operator is
  uniquely liftable to the tensor product $X \otimes_\uppi X$ by
  \thref{prop:lift-boun-bil-map}, the representative of $Q$ is just
  given by $\breve Q (u , u \otimes u) = \breve A( u \otimes u)$,
  where $\breve A$ is the restriction of the lifting of $A$ to the
  subspace $X \otimes_{\uppi, \mathrm{sym}} X$.  \qed
\end{Example}

\begin{Example}[Bilinear mappings]
  \label{ex:bil-map}
  Finally, the dilinear mappings also cover the class of bounded
  bilinear operators.  For this, we replace the \PBanach\ space $X$ by
  the \pers{Cartes}\-ian product $X_1 \times X_2$, where $X_1$ and $X_2$
  are arbitrary real \PBanach\ spaces.  Given a bounded bilinear operator
  $A \colon X_1 \times X_2 \rightarrow Y$, we define the symmetric
  bilinear mapping
  $B \colon (X_1 \times X_2) \times (X_1 \times X_2) \rightarrow Y$ by
  $B((u_1,v_1),(u_2,v_2)) = \nicefrac{1}{2} \, A(u_1, v_2) +
  \nicefrac{1}{2} \, A(u_2, v_1)$.
  Using the lifting $\breve B$ of $B$, we obtain the representative
  $\breve A ((u,v), (u,v)\otimes(u,v)) = \breve B ((u,v) \otimes
  (u,v))$ for all $(u,v)$ in $X_1 \times X_2$.  \qed
\end{Example}

\section{Generalized subgradient}
\label{sec:gener-subgr}

One of the drawbacks of the dilinear operators in
\thref{def:dilinear-map} is that dilinear mappings
$K \colon X \rightarrow \R$ do not have to be convex.  Hence, the
application of the usual subgradient to dilinear operators is limited.
To surmount this issue, we generalize the concept of convexity and of
the usual subgradient, see for instance \cite{BC11, ET79, Roc70}, to
our setting.  In the following, we denote the real numbers extended by
$+\infty$ and $-\infty$ by $\overline \R$.  For a mapping $F$ between
a real \PBanach\ space and the extended real numbers $\overline \R$,
the \emph{effective domain} is the section
\begin{equation*}
  \dom(F) \coloneqq
  \{u : F(u) < + \infty\}.
\end{equation*}
The mapping $F \colon X \rightarrow \overline \R$ is \emph{proper} if
it is never $-\infty$ and not everywhere $+\infty$.

\begin{Definition}[Diconvex mappings]
  \label{def:dicon-map}
  Let $X$ be real \PBanach\ spaces.  A mapping
  $F \colon X \rightarrow \overline \R$ is \emph{diconvex} if there
  exists a proper, convex mapping
  $\breve F \colon X \times (X \otimes_{\uppi,\mathrm{sym}} X)
  \rightarrow \overline \R$
  such that $F(u) = \breve F(u,u\otimes u)$ for all $u$ in $X$.
\end{Definition}

Since each proper, convex mapping
$F \colon X \rightarrow \overline \R$ may be represented by the convex
mapping $\breve F(u,u \otimes u) = F(u)$ on
$X \times ( X \otimes_{\uppi,\mathrm{sym}} X)$, we can view the
diconvex mappings as a generalization of the set of proper, convex
mappings. The central notion behind this definition is that each
dilinear functional is by definition also diconvex.  However,
differently from the dilinear operators, the \emph{representative}
$\breve F$ of a diconvex mapping $F$ does not have to be unique.  As
we will see, one sufficient condition for diconvexity is the existence
of a continuous, \emph{diaffine} minorant $A$, where diaffine means that there
exists a continuous, affine mapping
$\breve A \colon X \times (X \otimes_{\uppi, \mathrm{sym}} X)
\rightarrow Y$
such that $A(u) = \breve A(u, u \otimes u)$ for all $u$ in $X$.  In
this context, a continuous mapping $\breve A$ is \emph{affine} if and
only if $\breve A(u,w) \coloneqq T(u,w) + t_0$ for a bounded linear
operator
$T \colon X \times (X \otimes_{\uppi, \mathrm{sym}} X) \rightarrow \R$
and a constant $t_0 \in R$.

In order to prove this assertion, we will exploit that each vector
$(u, u \otimes u)$ is an \emph{extreme point} of the convex hull of
the diagonal $\{(u, u \otimes u) : u \in X\}$, which means that
$(u,u \otimes u)$ cannot be written as a non-trivial convex
combination of other points, see \cite{Roc70}.

\begin{Lemma}[Extreme points of the convex hull of the diagonal]
  \label{lem:extr-points-diag}
  Let $X$ be a real \PBanach\ space.  Each point $(u,u \otimes u)$
  with $u \in X$ is an extreme point of
  $\conv\{(u,u \otimes u) : u \in X \}$.
\end{Lemma}

\begin{Proof}
  For an element $(u, u \otimes u)$ with $u \in X$, we consider an
  arbitrary convex combination
  \begin{equation*}
    (u, u \otimes u) = \sum_{n=1}^N \alpha_n \, (u_n, u_n \otimes u_n)
    \addmathskip
  \end{equation*}
  with $u_n \in X$, $\alpha_n \in [0,1]$, and
  $\sum_{n=1}^N \alpha_n = 1$.  Applying the linear functionals
  $(\phi, 0 \otimes 0)$ and $(0, \phi \otimes \phi)$ with
  $\phi \in X^*$ of the dual space
  $X^* \times (X \otimes_{\uppi, \mathrm{sym}} X)^*$, we get the
  identity
  \begin{equation*}
    \biggl( \sum_{n=1}^N \alpha_n \iProd{\phi}{u_n} \biggr)^2 
    = \iProd{\phi}{u}^2 
    = \iProd{\phi \otimes \phi}{u \otimes u}
    = \sum_{n=1}^N  \alpha_n \iProd{\phi}{u_n}^2.
    \addmathskip
  \end{equation*}
  Due to the strict convexity of the square, this equation can only
  hold if $\iProd{\phi}{u_n} = \iProd{\phi}{u}$ for every $n$ between
  $1$ and $N$, and for every $\phi \in X^*$.  Consequently, all $u_n$
  coincide with $u$, which shows that the considered convex
  combination is trivial, and that $(u,u \otimes u)$ is an extreme
  point.  \qed
\end{Proof}

With the knowledge that the diagonal $\{(u,u \otimes u) : u \in X\}$
contains only extreme points of its convex hull, we may now give an
sufficient condition for a mapping being diconvex.

\begin{Proposition}
  \label{prop:diaff-minor}
  Let $X$ be a real \PBanach\ space.  If the mapping
  $F : X \rightarrow \overline R$ has a continuous, diaffine minorant,
  then $F$ is diconvex.
\end{Proposition}

\begin{Proof}                   
  If $F$ has a continuous, diaffine minorant $G$ with representative
  $\breve G$, we can construct a representative $\breve F$ by
  \begin{equation*}
    \breve F(u,w) \coloneqq
    \begin{cases}
      F(u) & w = u \otimes u, \\
      \breve G(u,w) & w \ne u \otimes u,
      ~\text{but}~ (u,w) \in \conv \{ (u,u \otimes u) : u \in X\},\\
      + \infty & \text{else.}
    \end{cases}
  \end{equation*}
  Since we firstly restrict the convex mapping $\breve G$ to the
  convex set $\{(u, u\otimes u) : u \in X\}$ and secondly increase the
  function values of the extreme points on the diagonal, the
  constructed mapping $\breve F$ is convex.  Obviously, the functional
  $\breve F$ is also proper and thus a valid representative.  \qed
\end{Proof}

\begin{Remark}
  \label{rem:diaff-minor:1}
  If the \PBanach\ space $X$ in \thref{prop:diaff-minor} is
  finite-dimensional, then the reverse implication is also true.  To
  validate this assertion, we restrict a given proper, convex
  representative $\breve F$ to the convex hull of
  $\{(u,u \otimes u) : u \in X\}$, which here means that we set
  $\breve F(v,w) = +\infty$ for $(v,w)$ outside the convex hull.  Due
  to the fact that the relative interior of a convex set is non-empty
  in finite dimensions, see \cite[Theorem~6.2]{Roc70}, there exists a
  point $(v,w)$ where the classical subdifferential
  $\uppartial \breve F(v,w)$ of the proper, convex representative
  $\breve F$ is non-empty, see for instance
  \cite[Theorem~23.4]{Roc70}.  In other words, we find a dual element
  $(\xi, \Xi) \in X^* \times (X \otimes_{\uppi,\mathrm{sym}} X)^*$
  such that
  \begin{equation*}
    \breve F (v',w') \ge 
    F(v,w) + \iProd{\xi}{v'-v} + \iProd{\Xi}{w'-w}
  \end{equation*}
  for all $v' \in X$ and $w' \in X \otimes_{\uppi,\mathrm{sym}} X$.
  Obviously, the functional $\breve A$ given by
  \begin{equation*}
    \breve A(v',w') \coloneqq
    F(v,w) + \iProd{\xi}{v'-v} + \iProd{\Xi}{w'-w}
  \end{equation*}
  defines a continuous, affine minorant of $\breve F$, and the
  restriction $A(u) \coloneqq \breve A(u, u \otimes u)$ thus a
  diaffine minorant $A$ of $F$.  \qed
\end{Remark}

\begin{Remark}
  \label{rem:diaff-minor:2}
  Looking back at the proof of \thref{prop:diaff-minor}, we can
  directly answer the question: why does the representative of a
  diconvex mapping has to be proper?  If we would allow inproper
  representatives as well, every mapping
  $F \colon X \rightarrow \overline \R$ would be diconvex with the
  convex but inproper representative
  \begin{equation*}
    \breve F(u,w) \coloneqq
    \begin{cases}
      F(u) & w = u \otimes u, \\
      -\infty & w \ne u \otimes u,
      ~\text{but}~ (u,w) \in \conv \{ (u,u \otimes u) : u \in X\},\\
      + \infty & \text{else;}
    \end{cases}
  \end{equation*}
  so the diconvex mappings would simply embrace all possible mappings
  between the \PBanach\ space $X$ and $\overline \R$.  The condition
  that the representative is proper will be needed at several points
  for the generalized subdifferential calculus and the developed
  regularization theory.  \qed
\end{Remark}

After this preliminary considerations, we now generalize the classical
subgradient and subdifferential to the class of (proper) diconvex
mappings.

\begin{Definition}[Dilinear subgradient]
  \label{def:dil-subgrad}
  Let $F \colon X \rightarrow \overline \R$ be a diconvex mapping on
  the real \PBanach\ space $X$.  The dual element
  $(\xi,\Xi) \in X^* \times (X \otimes_{\uppi, \mathrm{sym}} X)^*$ is
  a \emph{dilinear subgradient} of $F$ at $u$ if $F(u)$ is finite and
  if
  \begin{equation*}
    F(v) \ge F(u) + \iProd{\xi}{v-u} + \iProd{\Xi}{v \otimes v - u
      \otimes u}
  \end{equation*}
  for all $v$ in $X$.  The union of all dilinear subgradients of $F$
  at $u$ is the \emph{dilinear subdifferential}
  $\uppartial_\upbeta F(u)$.  If no dilinear subgradient exists, the
  dilinear subdifferential is empty.
\end{Definition}

If the mapping $F$ is convex, then the dilinear subdifferential
obviously contains the usual subdifferential of $F$.  More precisely,
we have $\uppartial_\upbeta F(u) \supset \uppartial F(u) \times \{0\}$.
Where the usual subgradient consists of all linear functionals
entirely lying below the mapping $F$, the dilinear subgradient
consists of all dilinear mappings below $F$. In this context, the
dilinear subdifferential can be interpreted as the $W$-subdifferential
introduced in \cite{Gra10} with respect to the family of dilinear
functionals whereas the bilinear part $\Xi$ here does not have to be
negative semi-definite.  If we think at the one-dimensional case
$F \colon \R \rightarrow \overline \R$, the dilinear subdifferential
embraces all parabolae beneath $F$ at a certain point $u$.  Since each
diconvex mapping $F$ has at least one convex representative
$\breve F$, we next investigate how the representative $\breve F$ may
be used to compute the dilinear subdifferential
$\uppartial_\upbeta F(u)$.

\begin{Definition}[Representative subgradient]
  \label{def:rep-subgrad}
  Let $F \colon X \rightarrow \overline \R$ be a diconvex mapping on
  the real \PBanach\ space $X$ with representative
  $\breve F \colon X \times (X \otimes_{\uppi, \mathrm{sym}} X)
  \rightarrow \overline \R$.
  The dual element
  $(\xi, \Xi) \in X^* \times (X \otimes_{\uppi, \mathrm{sym}} X)^*$ is
  a \emph{representative subgradient} of $F$ at $u$ with respect to
  $\breve F$ if $(\xi, \Xi)$ is a subgradient of the representative
  $\breve F$ at $(u, u \otimes u)$.  The union of all representative
  subgradients of $F$ at $u$ is the \emph{representative
    subdifferential} $\breve \uppartial F(u)$ with respect to
  $\breve F$.
\end{Definition}

Since the representative $\breve F$ of a diconvex mapping $F$ may not
be unique, the representative subgradient usually depends on the
choice of the mapping $\breve F$.  Nevertheless, a representative
subgradient is as well a dilinear subgradient.

\begin{Lemma}[Inclusion of subdifferentials]
  \label{lem:inc-subdiff}
  Let $F \colon X \rightarrow \overline \R$ be a diconvex mapping on
  the real \PBanach\ space $X$ with representative $\breve F$.  Then
  the representative and dilinear subdifferential are related by
  \begin{equation*}
    \breve \uppartial F(u) \subset \uppartial_\upbeta F(u).
  \end{equation*}
\end{Lemma}

\begin{Proof}
  Since each representative subgradient $(\xi, \Xi)$ fulfils
  \begin{equation*}
    \breve F(v, w) \ge \breve F(u, u \otimes u) + \iProd{\xi}{v-u} +
    \iProd{\Xi}{w-u \otimes u}
  \end{equation*}
  for all $(v,w)$ in $X \times (X \otimes_{\uppi, \mathrm{sym}} X)$
  and thus especially for $(v,w) = (v, v \otimes v)$, the asserted
  inclusion follows.  \qed
\end{Proof}

In view of the fact that the representative $\breve F$ of a diconvex
mapping $F$ is not unique, the question arises whether there exists a
certain representative $\breve F$ such that the representative
subdifferential coincides with the dilinear subdifferential.  Indeed,
we can always construct an appropriate representative by considering
the convexification of $F$ on
$X \times (X \otimes_{\uppi, \mathrm{sym}} X)$.  In this context, the
convexification $\conv G$ of an arbitrary functional
\raisebox{0pt}[0pt][0pt]{$G \colon X \rightarrow \overline \R$} on the
\PBanach\ space $X$ is the greatest convex function majorized by $G$
and can be determined by
\begin{equation}
  \label{eq:def-conv-hull-map}
  \conv G (u) = \inf \biggl\{ \sum_{n=1}^N \alpha_n \, G(u_n) \colon u
  = \sum_{n=1}^N \alpha_n \, u_n
\biggr\},
\end{equation}
where the infimum is taken over all convex representations
$u = \sum_{n=1}^N \alpha_n \, u_n$ with $n\in \N$, $u_n \in X$, and
$\alpha_n \in [0,1]$ so that $\sum_{n=1}^N \alpha_n = 1$, see for
instance \cite{Roc70}.  For a diconvex functional
$F \colon X \rightarrow \overline \R$, we now consider the
convexification of
\begin{equation}
  \label{eq:F-otimes}
  F_\otimes (u,w) \coloneqq
  \begin{cases}
    F(u) & w = u \otimes u, \\
    +\infty & \text{else}.
  \end{cases}
\end{equation}
Obviously, the mapping $\conv F_\otimes$ as supremum of all convex
functionals majorized by $F_\otimes$ is a valid representative of $F$
since there exists at least one convex representative $\breve F$ with
$F(u) = \breve F (u, u \otimes u)$ for all $u$ in $X$.  Further, the
convexification $\conv F_\otimes$ is also proper since the
representative $\breve F$ has to be proper.

\begin{Theorem}[Equality of subdifferentials]
  \label{the:equa-subdiff}
  Let $F \colon X \rightarrow \overline \R$ be a diconvex mapping on
  the real \PBanach\ space $X$.  Then the representative
  subdifferential with respect to $\conv F_\otimes$ and the dilinear
  subdifferential coincide, \ie
  \begin{equation*}
    \breve \uppartial F(u) = \uppartial_\upbeta F(u).
  \end{equation*}
\end{Theorem}

\begin{Proof}
  Let $(\xi, \Xi)$ be a dilinear subgradient of $F$ at $u$, which
  means
  \begin{equation*}
    F(v) \ge F(u) + \iProd{\xi}{v-u} + \iProd{\Xi}{v \otimes v - u
      \otimes u}
  \end{equation*}
  for all $v$ in $X$.  Hence, for each convex combination
  $\sum_{n=1}^N \alpha_n (u_n, u_n \otimes u_n)$, we have
  \begin{align*}
    &\sum_{n=1}^N \alpha_n \, F_\otimes (u_n, u_n \otimes u_n)
    \\[\fskip]
    &\qquad\ge \sum_{n=1}^N \alpha_n \, \bigl[ F_\otimes (u, u \otimes u) +
      \iProd{\xi}{u_n - u} + \iProd{\Xi}{u_n \otimes u_n - u \otimes
      u} \bigr]
    \\[\fskip]
    &\qquad= F_\otimes(u, u\otimes u) 
      + \iProd{\xi}{\sum_{n=1}^N \alpha_n \, u_n - u} 
      + \iProd{\Xi}{\sum_{n=1}^N \alpha_n \, u_n \otimes u_n - u \otimes
      u}.
  \end{align*}
  Taking the infimum over all convex combinations
  $(v,w) = \sum_{n=1}^N \alpha_n (u_n, u_n \otimes u_n)$ for an
  arbitrary point $(v,w)$ in
  $X \times (X \otimes_{\uppi, \mathrm{sym}} X)$, we obtain
  \begin{equation*}
    \conv F_\otimes (v,w) \ge \conv F_\otimes (u, u \otimes u) 
    + \iProd{\xi}{v-u} + \iProd{\Xi}{w-u \otimes u}.
  \end{equation*}
  Consequently, the dilinear subgradient $(\xi, \Xi)$ is contained in
  the representative subdifferential $\breve \uppartial F (u)$ with
  respect to $\conv F_\otimes$.  \qed
\end{Proof}

Similarly to the classical (linear) subdifferential, the dilinear
subdifferential of a sum $F+G$ contains the sum of the single dilinear
subdifferentials of $F$ and $G$.

\begin{Proposition}[Sum rule]
  \label{prop:sum-rule}
  Let $F \colon X \rightarrow \overline \R$ and
  $G \colon X \rightarrow \overline \R$ be diconvex mappings on the
  real \PBanach\ space $X$.  Then the dilinear subdifferential of
  $F+G$ and the dilinear subdifferentials of $F$ and $G$ are related
  by
  \begin{equation*}
    \uppartial_\upbeta F(u) + \uppartial_\upbeta G(u)
    \subset \uppartial_\upbeta [F+G](u)
  \end{equation*}
\end{Proposition}

\begin{Proof}
  This is an immediate consequence of \thref{def:dil-subgrad}.  \qed
\end{Proof}

Differently from the classical (linear) subdifferential, we cannot
transfer the chain rule to the dilinear/diconvex setting.  The main
reason is that the composition of a diconvex mapping $F$ and a
dilinear mapping $K$ has not to be diconvex, since a representative of
$F \circ K$ cannot simply be constructed by composing the
representatives of $F$ and $K$.  Therefore, the chain rule can only be
transferred partly.  For an arbitrary bounded linear operator
$K \colon X \rightarrow Y$, we recall that there exists a unique
bounded linear operator
$K \otimes_\uppi K \colon X \otimes_\uppi X \rightarrow Y
\otimes_\uppi Y$
such that $(K \otimes_\uppi K) (u \otimes v) = K(u) \otimes K(v)$, see
for instance \cite{Rya02}.  In the following, the restriction of the
lifted operator $K \otimes_\uppi K$ to the symmetric subspace is
denoted by
$K \otimes_{\uppi,\mathrm{sym}} K \colon X \otimes_{\uppi,
  \mathrm{sym}} X \rightarrow Y \otimes_{\uppi, \mathrm{sym}} Y$.  
Moreover, we recall that the \PCartes{}ian product
\begin{equation*}
  K \times (K \otimes_{\uppi,\mathrm{sym}} K) 
  \colon X \times (X \otimes_{\uppi,\mathrm{sym}} X) 
  \rightarrow Y \times (Y \otimes_{\uppi,\mathrm{sym}} Y)
\end{equation*}
of the mappings $K$ and $K \otimes_{\uppi,\mathrm{sym}} K$ is defined
by
\begin{equation*}
  [K \times (K \otimes_{\uppi,\mathrm{sym}} K)](u,w) 
  \coloneqq (K(u), K \otimes_{\uppi,\mathrm{sym}} K (w)).
\end{equation*}

\begin{Proposition}[Chain rule for linear operators]
  \label{prop:chain-rule:lin-op}
  Let $K \colon X \rightarrow Y$ be a bounded linear mapping and
  $F \colon Y \rightarrow \overline \R$ be a diconvex mapping on the
  real \PBanach\ spaces $X$ and $Y$.  Then the dilinear
  subdifferential of $F \circ K$ is related to the dilinear
  subdifferential of $F$ by
  \begin{equation*}
    (K \times (K \otimes_{\uppi, \mathrm{sym}} K))^*  \uppartial_\upbeta F(K(u)) 
    \subset \uppartial_\upbeta (F \circ K) (u).
  \end{equation*}
\end{Proposition}

\begin{Proof}
  Firstly, we notice that the functional $F \circ K$ is diconvex with
  the representative
  $\breve F \circ (K \times (K \otimes_{\uppi, \mathrm{sym}} K))$.
  Next, let us assume $(\xi, \Xi)\in \uppartial_\upbeta F(K(u))$, which
  is equivalent to
  \begin{equation*}
    F(v) \ge F(K(u)) + \iProdn{\xi}{v - K(u)} + \iProdn{\Xi}{v \otimes v
    - K(u) \otimes K(u)}
  \end{equation*}
  for all $v$ in $Y$.  Replacing $v \in Y$ by $K(v)$ with $v \in X$,
  we obtain
  \begin{equation*}
    (F \circ K) (v) \ge (F \circ K) (u) + \iProdn{\xi}{K(v-u)} +
    \iProdn{\Xi}{(K \otimes_{\uppi, \mathrm{sym}} K)
      (v \otimes v - u \otimes u)}
  \end{equation*}
  for all $v$ in $X$.  Thus,
  $(K \times (K \otimes_{\uppi, \mathrm{sym}} K))^* (\xi, \Xi)$ is
  contained in the subdifferential $\uppartial_\upbeta F \circ K (u)$.
  \qed
\end{Proof}

\begin{Proposition}[Chain rule for convex functionals]
  \label{prop:chain-rule:con-fun}
  Let $K \colon X \rightarrow Y$ be a bounded dilinear mapping and
  $F \colon Y \rightarrow \overline \R$ be a convex mapping on the
  real \PBanach\ spaces $X$ and $Y$.  Then the dilinear
  subdifferential of $F \circ K$ is related to the linear
  subdifferential of $F$ by
  \begin{equation*}
    \breve K^* \uppartial F(K(u)) \subset \uppartial_\upbeta (F
    \circ K) (u).
  \end{equation*}
\end{Proposition}

\begin{Proof}
  The functional $F \circ K$ is diconvex with convex representative
  $F \circ \breve K$, since the lifted operator $\breve K$ is linear.
  Next, we consider a linear subgradient $\xi$ of $F$ at $K(u)$, which
  means
  \begin{equation*}
    F(v) \ge F(K(u)) + \iProd{\xi}{v - K(u)}
    \addmathskip
  \end{equation*}
  for all $v$ in $Y$.  Replacing $v \in Y$ by $K(v)$ with $v \in X$,
  we obtain
  \begin{equation*}
    (F \circ K)(v) \ge (F \circ K) (u) + \iProd{\xi}{\breve K (v-u, v
      \otimes v - u \otimes u)}
  \end{equation*}
  for all $v$ in $X$.  Consequently, $\breve K^* (\xi)$ is contained
  in the subdifferential $\uppartial_\upbeta (F \circ K) (u)$.  \qed
\end{Proof}

Since the representative subdifferential is based on the classical
subdifferential on the lifted space, the classical sum and chain rules
obviously remain valid whenever the representatives fulfil the
necessary conditions.  For instance, one functional of the sum or the
outer functional of the composition is continuous at some point.  At
least, for finite-dimensional \PBanach\ spaces, the representative of
a diconvex functional, which is finite on some open set, is continuous
on the interior of its effective domain.  The central idea to prove
this conjecture is that the manifold of rank-one tensors is curved in
a way such that the convex hull of each open set of the manifold
contains an inner point with respect to the surrounding tensor
product.  In the following, we denote by $\overline B_\epsilon$ the
closed $\epsilon$-ball around zero.

\begin{Lemma}[Local convex hull of rank-one tensors]
  \label{lem:inn-pt-conv-hull}
  Let $u$ be a point of the real finite-dimensional \PBanach\ space
  $X$.  The interior of the convex hull of the set
  \begin{equation}
    \label{eq:inn-pt-lift-ball}
    \bigl\{ (u+v, (u+v)\otimes (u+v)) : v \in \overline B_\epsilon \bigr\}
    \subset X \times (X \otimes_{\uppi, \mathrm{sym}} X)
    \submathskip
  \end{equation}
  is not empty for every $\epsilon > 0$.
\end{Lemma}

\begin{Proof}
  To determine a point in the interior, we construct a suitable
  simplex by using a normalized basis $(e_n)_{n=1}^N$ of the \PBanach\
  space $X$.  Obviously, the convex hull contains the points
  $(u + \epsilon e_n, (u + \epsilon e_n) \otimes (u + \epsilon e_n))$.
  Viewing the point $(u, u \otimes u)$ as new origin of
  $X \times (X \otimes_{\uppi, \mathrm{sym}} X)$, we obtain the
  vectors
  \begin{equation}
    \label{eq:inn-pt-bas-vec-1}
    \begin{aligned}
      &(u + \epsilon e_n, (u + \epsilon e_n) \otimes (u + \epsilon e_n))
      - (u, u \otimes u) 
      \\[\fskip]
      &\qquad= (\epsilon e_n, \epsilon (u \otimes e_n) + \epsilon(e_n \otimes
      u) + \epsilon^2 (e_n \otimes e_n)),
    \end{aligned}
  \end{equation}
  where the first components again form a basis of the first component
  space $X$.  Next, we consider the convex combination of
  $(u + \epsilon e_n, (u+ \epsilon e_n) \otimes (u + \epsilon e_n))$
  and
  $(u - \epsilon e_n, (u - \epsilon e_n) \otimes (u - \epsilon e_n))$
  with weights $\nicefrac{1}{2}$.  In this manner, we obtain
  \begin{equation}
    \label{eq:inn-pt-bas-vec-2}
    \begin{aligned}
      &\tfrac{1}{2} \, (u + \epsilon e_n, (u+ \epsilon e_n) \otimes (u +
      \epsilon e_n)) 
      \\[\fsmallskip]
      &\qquad+ \tfrac{1}{2} \, (u - \epsilon e_n, (u - \epsilon e_n) \otimes (u
      - \epsilon e_n)) 
      - (u, u \otimes u) 
      \\[\fskip]
      &\quad= (0, \epsilon^2 (e_n \otimes e_n)).
    \end{aligned}
  \end{equation}
  Similarly, by considering the corresponding convex combination of
  $(u + \nicefrac{\epsilon}{2} (e_n + e_m), (u +
  \nicefrac{\epsilon}{2} (e_n + e_m) \otimes (u +
  \nicefrac{\epsilon}{2} (e_n + e_m))$
  and
  $(u - \nicefrac{\epsilon}{2} (e_n + e_m), (u -
  \nicefrac{\epsilon}{2} (e_n + e_m) \otimes (u -
  \nicefrac{\epsilon}{2} (e_n + e_m))$
  with $n \ne m$, we get the vectors
  \begin{equation}
    \label{eq:inn-pt-bas-vec-3}
    \begin{aligned}
      &(0, \nicefrac{\epsilon^2}{4} ((e_n + e_m) \otimes (e_n + e_m)))
      \\[\fskip]
      &\quad= (0, \nicefrac{\epsilon^2}{4} (e_n \otimes e_n)) + (0,
      \nicefrac{\epsilon^2}{4} (e_n \otimes e_m + e_m \otimes e_n)) +
      (0, \nicefrac{\epsilon^2}{4} (e_m \otimes e_m)).
    \end{aligned}
  \end{equation}
  Since the vectors $(e_n \otimes e_m)_{n,m =1}^N$ form a basis of the
  tensor product $X \otimes_\uppi X$, see \cite{Rya02}, the vectors in
  \eqref{eq:inn-pt-bas-vec-2} and \eqref{eq:inn-pt-bas-vec-3} span the
  second component space $X \otimes_{\uppi, \mathrm{sym}} X$.  Thus,
  the vectors (\ref{eq:inn-pt-bas-vec-1}--\ref{eq:inn-pt-bas-vec-3})
  form a maximal set of independent vectors, and the convex hull of
  them cannot be contained in a true subspace of
  $X \times (X \otimes_{\uppi, \mathrm{sym}} X)$.  Consequently, the
  simplex spanned by the vectors
  (\ref{eq:inn-pt-bas-vec-1}--\ref{eq:inn-pt-bas-vec-3}) and zero
  contains an inner point.  Since the constructed simplex  shifted
  by $(u, u \otimes u)$ is contained in the convex hull of
  \eqref{eq:inn-pt-lift-ball}, the assertion is established.  \qed
\end{Proof}

Unfortunately, \thref{lem:inn-pt-conv-hull} does not remain valid for
infinite-dimensional \PBanach\ spaces.  For example, if the \PBanach\
space $X$ has a normalized \PSchauder\ basis $(e_n)_{n \in \N}$, we
can explicitly construct a vector not contained in the convex hull of
\eqref{eq:inn-pt-lift-ball} but arbitrarily near to a given element of
the convex hull.  For this, we notice that the tensor product
$X \otimes_{\uppi} X$ possesses the normalized \PSchauder\ basis
$(e_n \otimes e_m)_{n,m \in \N}$ with respect to the square ordering,
see \cite{Rya02}, and that the coordinates of an arbitrary rank-one
tensor
$u \otimes u = \sum_{(n,m) \in \N^2} a_{nm} \, (e_n \otimes e_m) $
with $u = \sum_{n \in \N} u_n \, e_n$ are given by $a_{nm} = u_n u_m$.
Consequently, the coordinates $a_{nn}$ on the diagonal have to be
non-negative.  This implies that the convex hull of
\eqref{eq:inn-pt-lift-ball} only contains tensors with non-negative
diagonal.

Now, let $(v,w)$ be an arbitrary element of the convex hull of
\eqref{eq:inn-pt-lift-ball}.  Since the representation
$\sum_{(n,m) \in \N^2} b_{nm} \, (e_n \otimes e_m)$ of the tensor $w$
converges with respect to the square ordering, the coordinates
$b_{nm}$ form a zero sequence.  Therefore, for each given
$\delta > 0$, we find an $N \in \N$ such that $b_{nn} < \delta$
whenever $n \ge N$.  Subtracting the vector
$(0,\delta \, (e_N \otimes e_N))$ from $(v,w)$, we obtain an
arbitrarily near vector to $(v,w)$ that is not contained in the convex
hull, since one coordinate on the diagonal is strictly negative.
Thus, the convex hull of \eqref{eq:inn-pt-lift-ball} has an empty
interior.

\begin{Proposition}[Continuity of the representative]
  \label{prop:cont-repre}
  Let
  $\breve F \colon X \times (X \otimes_{\uppi, \mathrm{sym}} X)
  \rightarrow \overline \R$
  be a representative of the diconvex mapping
  $F \colon X \rightarrow \overline \R$ on the real finite-dimensional
  \PBanach\ space $X$.  If $F$ is finite on a non-empty, open set,
  then the representative $\breve F$ is continuous in the interior of
  the effective domain $\dom(\breve F)$.
\end{Proposition}

\begin{Proof}
  By assumption, there is a point $u \in X$ such that $F$ is finite on
  an $\epsilon$-ball $\overline B_\epsilon(u)$ around $u$ for some
  $\epsilon > 0$.  Consequently, the representative $\breve F$ is
  finite on the set
  \begin{equation}
    \label{eq:cont-repre:ball}
    \bigl\{(u+v, (u+v) \otimes (u+v)) : v \in \overline B_\epsilon\bigr\}.
  \end{equation}
  Using the construction in the proof of \thref{lem:inn-pt-conv-hull},
  we find a simplex with vertices in \eqref{eq:cont-repre:ball} that
  contains an inner point $(v,w)$ of the convex hull of
  \eqref{eq:cont-repre:ball} and hence of the effective domain
  $\dom(\breve F)$.  Since $\breve F$ is convex and finite on the
  vertices of the constructed simplex, the representative $\breve F$
  is bounded from above on a non-empty, open neighbourhood around
  $(v,w)$, which is equivalent to the continuity of $\breve F$ on the
  interior of the effective domain $\dom(\breve F)$, see for instance
  \cite{ET79, Sho97}.  \qed
\end{Proof}

On the basis of this observation, we obtain the following computation
rules for the representative subdifferential on finite-dimensional
\PBanach\ spaces, which follow immediately form
\thref{prop:cont-repre} and the classical sum and chain rules, see for
instance \cite{BC11, ET79, Roc70, Sho97}.

\begin{Proposition}[Representative sum rule]
  \label{prop:repre-sum-rule}
  Let $F \colon X \rightarrow \overline \R$ and
  $G \colon X \rightarrow \overline \R$ be diconvex functionals on the
  real finite-dimensional \PBanach\ space $X$ with
  representatives $\breve F$ and $\breve G$.  If there
    exists a non-empty, open set where $F$ and $G$ are finite, then
  \begin{equation*}
    \breve \uppartial [ F + G ](u) = \breve \uppartial F(u) +
    \breve \uppartial G(u)
    \addmathskip
  \end{equation*}
  for all $u$ in $X$ with respect to the representative
  $\breve F + \breve G$ of $F+G$.
\end{Proposition}

\begin{Proof}
  In the proof of \thref{prop:cont-repre}, we have constructed an
  appropriate simplex to prove the existence of a point $(v,w)$ in
  $X \times (X \otimes_{\uppi,\mathrm{sym}} X)$ where the
  representative $\breve F$ is finite and continuous.  Using the same
  simplex again, we see that the representative $\breve G$ is finite
  and continuous in the same point $(v,w)$.  Applying the classical
  sum rule -- see for instance \cite[Proposition~II.7.7]{Sho97} -- to
  the representative $\breve F + \breve G$, we obtain
  \begin{equation*}
    \uppartial (\breve F + \breve G)(u,u\otimes u) = \uppartial \breve
    F(u,u\otimes u) + \uppartial \breve G(u,u \otimes u)
  \end{equation*}
  for all $(u, u \otimes u)$ in
  $X \times (X \otimes_{\uppi,\mathrm{sym}} X)$ and thus the
  assertion.  \qed
\end{Proof}

\begin{Remark}
  \label{rem:repre-sum-rule}
  In order to apply the classical sum rule in the proof of
  \thref{prop:repre-sum-rule}, it would be sufficient if only one of
  the functionals $\breve F$ and $\breve G$ is continuous in $(v,w)$.
  The assumption that $F$ and $G$ are finite at some non-empty, open
  set is thus stronger than absolutely necessary.  On the other side,
  this assumption is needed to ensure that the effective domain of
  $\breve G$ and the interior of the effective domain of $\breve F$
  have some point in common.  \qed
\end{Remark}

\begin{Proposition}[Representative chain rule for linear operators]
  \label{prop:repre-chain-rule:lin-opt}
  Let $K \colon X \rightarrow Y$ be a bounded, surjective linear
  mapping and $F \colon Y \rightarrow \overline \R$ be a diconvex
  functional with representative $\breve F$ on the real \PBanach\
  spaces $X$ and $Y$.  If\enspace $Y$ is finite-dimensional, and if
  there exists a non-empty, open set where $F$ is finite, then
  \begin{equation*}
    \breve \uppartial (F \circ K) (u)
    = (K \times (K \otimes_{\uppi, \mathrm{sym}} K))^* \breve \uppartial
    F(K(u))
  \end{equation*}
  for all $u$ in $X$ with respect to the representative
  $\breve F \circ (K \times (K \otimes_{\uppi, \mathrm{sym}} K))$ of
  $F \circ K$.
\end{Proposition}

\begin{Proof}
  Like in the proof of \thref{prop:chain-rule:lin-op}, the mapping
  $\breve F \circ (K \times (K \otimes_{\uppi, \mathrm{sym}} K))$ is a
  proper, convex representative of $F \circ K$.  Since the linear
  operator $K$ is surjective, the mapping
  $K \otimes_{\uppi, \mathrm{sym}} K$ is surjective too.  In more
  detail, there exists finitely many vectors $e_1, \dots, e_N$ such
  that the images $f_n = K(e_n)$ form a basis of the
  finite-dimensional \PBanach\ space $Y$.  Since the symmetric tensors
  \begin{equation*}
    f_n \otimes f_m + f_m \otimes f_n = (K \otimes_{\uppi,
      \mathrm{sym}} K)(e_n \otimes e_m + e_m \otimes e_n)
  \end{equation*}
  form a basis of $Y \otimes_{\uppi,\mathrm{sym}} Y$, \cf\
  \cite[Proposition~1.1]{Rya02}, the bounded linear mapping
  $K \otimes_{\uppi, \mathrm{sym}} K$ is also surjective.  Using
  \thref{prop:cont-repre}, we thus always find a point
  $(K \times ( K \otimes_{\uppi, \mathrm{sym}} K))(v,w)$ where
  $\breve F$ is continuous.  Now, the classical chain rule -- see for
  instance \cite[Proposition~II.7.8]{Sho97} -- implies
  \begin{align*}
    &\uppartial ( \breve F \circ (K \times (K \otimes_{\uppi,
    \mathrm{sym}} K)))(u,u \otimes u)
    \\[\fskip]
    &\qquad= (K \times (K \otimes_{\uppi, \mathrm{sym}} K)^* \uppartial \breve
    F  ((K \times (K \otimes_{\uppi,
      \mathrm{sym}} K))(u,u \otimes u))
    \\[\fskip]
    &\qquad= (K \times (K \otimes_{\uppi, \mathrm{sym}} K)^* \uppartial \breve
    F  (K (u),K (u) \otimes K (u))
  \end{align*}
  for all $u$ in $X$, which establishes the assertion.  \qed
\end{Proof}

\begin{Remark}
  \label{rem:repre-chain-rule:lin-opt}
  In the proof of \thref{prop:repre-chain-rule:lin-opt}, the
  surjectivity of $K$ implies the non-emptiness of the intersection
  between the interior of $\dom(\breve F)$ and the range of
  $K \times (K \otimes_{\uppi, \mathrm{sym}} K)$.  So long as this
  intersection is not empty, \thref{prop:repre-chain-rule:lin-opt}
  remains valid even for non-surjective operators $K$.  The
  non-emptiness of the intersection then depends on the representative
  $\breve F$ and the mapping
  $K \times (K \otimes_{\uppi, \mathrm{sym}} K)$.  The intention
  behind \thref{prop:repre-chain-rule:lin-opt} has been to give a
  chain rule that only depends on properties of the given $F$ and $K$.
  \qed
\end{Remark}

\begin{Proposition}[Representative chain rule for convex functionals]
  \label{prop:repre-chain-rule:con-fun}
  Let $K \colon X \rightarrow Y$ be a bounded dilinear mapping and
  $F \colon Y \rightarrow \overline \R$ be a proper, convex functional
  on the real \PBanach\ spaces $X$ and $Y$.  If there exists a
  non-empty, open set where  $F$ is bounded from
  above, and if the interior of the effective domain $\dom(F)$ and the
  range $\ran(K)$ are not disjoint, then
  \begin{equation*}
    \breve \uppartial (F \circ K) (u)
    = \breve K^* \uppartial F(K(u))
    \addmathskip
  \end{equation*}
  for all $u$ in $X$ with respect to the representative $F \circ
  \breve K$ of $F \circ K$.
\end{Proposition}

\begin{Proof}
  Since the proper and convex mapping $F$ is bounded from above on
  some non-empty, open set, the function $F$ is continuous on the
  interior its effective domain $\dom F$, see for instance
  \cite[Proposition~2.5]{ET79}.  Consequently, we always find a point
  $\breve K (v,v \otimes v)$ where $F$ is continuous, which allows us
  to apply the classical chain rule -- see for instance
  \cite[Proposition~II.7.8]{Sho97} -- to the representative
  $F \circ \breve K$.  In this manner, we obtain
  \begin{equation*}
    \uppartial (F \circ \breve K)(u, u \otimes u)
    = \breve K^* \uppartial F(\breve K(u,u \otimes u))
    \submathskip
  \end{equation*}
  for all $u$ in $X$.  \qed
\end{Proof}

Although the dilinear and the representative subdifferential of a
diconvex functional $F$ coincide with respect to the representative
$\conv F_\otimes$, the established computation rules for the
representative subdifferential cannot be transferred to the dilinear
subdifferential in general.  A counterexample for the sum rule is
given below.  The main reason for this shortcoming is that the
convexification of $(F+G)_\otimes$ does not have to be the sum of the
convexifications of $F_\otimes$ and $G_\otimes$.  An analogous problem
occurs for the composition.

\begin{Counterexample}
  One of the simplest counterexamples, where the sum rule is failing
  for the dilinear subdifferential, is the sum of the absolute value
  function $\absn{\cdot} \colon \R \rightarrow \overline \R$ and
  the indicator function $\chi \colon \R \rightarrow \overline \R$ of
  the interval $[-1,1]$ given by
  \begin{equation*}
    \chi(t) \coloneqq 
    \begin{cases}
      0 & t \in [-1,1] \\
      + \infty & \text{else.}
    \end{cases}
  \end{equation*}
  As mentioned above, for a function from $\R$ into $\overline \R$, the dilinear
  subdifferential consists of all parabolae beneath that function at a
  certain point.  Looking at the point zero, we have the dilinear
  subdifferentials
  \begin{equation*}
    \uppartial_\upbeta \absn{\cdot} (0) 
    = \{t \mapsto c_2 t^2 + c_1 t : c_2 \le 0, c_1 \in [-1,1] \}
    \quad\text{and}\quad
    \uppartial_\upbeta \chi(0)
    = \{t \mapsto c_2 t^2 : c_2 \le 0 \}.
  \end{equation*}
  The sum of the dilinear subdifferentials thus consists of all
  parabolae with leading coefficient $c_2 \le 0$ and linear
  coefficient $c_1 \in [-1,1]$.  However, the dilinear subdifferential of the
  sum also contains parabolae with positive leading coefficient, see
  the schematic illustrations in \autoref{fig:dil-subdiff-coun-ex}.  \qed
  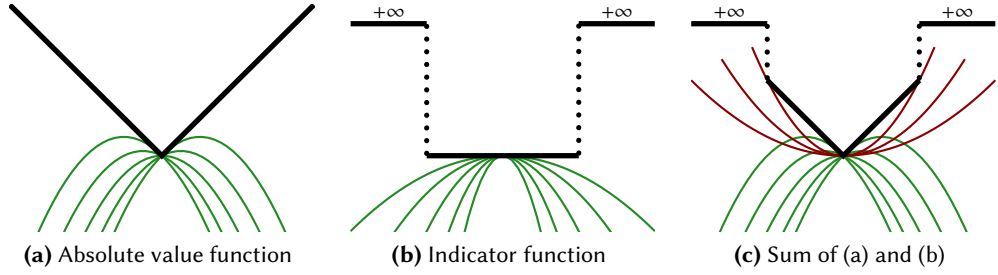
\begin{figure}
    \centering
    \subfloat[Absolute value function]{
      \begin{pspicture*}(-2,-1)(2,2)
        \psplot[algebraic=true,linecolor=ForestGreen]{-2}{2}{-x^2-x}
        \psplot[algebraic=true,linecolor=ForestGreen]{-2}{2}{-x^2+x}
        \psplot[algebraic=true,linecolor=ForestGreen]{-2}{2}{-x^2+x/2}
        \psplot[algebraic=true,linecolor=ForestGreen]{-2}{2}{-x^2-x/2}
        \psplot[algebraic=true,linecolor=ForestGreen]{-2}{2}{-x^2}
        \psline[linewidth=2pt](-2,2)(0,0)(2,2)
      \end{pspicture*}}\quad
    \subfloat[Indicator function]{
      \begin{pspicture*}(-2,-1)(2,2)
        \psplot[algebraic=true,linecolor=ForestGreen]{-2}{2}{-1/4*x^2}
        \psplot[algebraic=true,linecolor=ForestGreen]{-2}{2}{-1/2*x^2}
        \psplot[algebraic=true,linecolor=ForestGreen]{-2}{2}{-x^2}
        \psplot[algebraic=true,linecolor=ForestGreen]{-2}{2}{-2*x^2}
        \psplot[algebraic=true,linecolor=ForestGreen]{-2}{2}{-4*x^2}
        \psline[linewidth=2pt](-2,1.75)(-1,1.75)
        \psline[linewidth=2pt,linestyle=dotted](-1,1.75)(-1,0)
        \psline[linewidth=2pt](-1,0)(1,0)
        \psline[linewidth=2pt,linestyle=dotted](1,0)(1,1.75)
        \psline[linewidth=2pt](1,1.75)(2,1.75)
        \uput{2pt}[u](-1.5,1.75){\scriptsize$+\infty$}
        \uput{2pt}[u](1.5,1.75){\scriptsize$+\infty$}
      \end{pspicture*}}\quad
    \subfloat[Sum of (a) and (b)]{
      \begin{pspicture*}(-2,-1)(2,2)
        \psplot[algebraic=true,linecolor=ForestGreen]{-2}{2}{-x^2-x}
        \psplot[algebraic=true,linecolor=ForestGreen]{-2}{2}{-x^2+x}
        \psplot[algebraic=true,linecolor=ForestGreen]{-2}{2}{-x^2+x/2}
        \psplot[algebraic=true,linecolor=ForestGreen]{-2}{2}{-x^2-x/2}
        \psplot[algebraic=true,linecolor=ForestGreen]{-2}{2}{-x^2}
        \psplot[algebraic=true,linecolor=DarkRed]{-2}{2}{1/4*x^2}
        \psplot[algebraic=true,linecolor=DarkRed]{-1.6}{1.6}{1/2*x^2}
        \psplot[algebraic=true,linecolor=DarkRed]{-1.2}{1.2}{x^2}
        \psline[linewidth=2pt](-2,1.75)(-1,1.75)
        \psline[linewidth=2pt,linestyle=dotted](-1,1.75)(-1,1)
        \psline[linewidth=2pt](-1,1)(0,0)(1,1)
        \psline[linewidth=2pt,linestyle=dotted](1,1)(1,1.75)
        \psline[linewidth=2pt](1,1.75)(2,1.75)
        \uput{2pt}[u](-1.5,1.75){\scriptsize$+\infty$}
        \uput{2pt}[u](1.5,1.75){\scriptsize$+\infty$}
      \end{pspicture*}}
    \caption{Dilinear subdifferential of the absolute value function,
      the indicator function, and the sum.}
    \label{fig:dil-subdiff-coun-ex}
  \end{figure}
\end{Counterexample}

Generalizing the classical \PFermat\ rule, see for instance
\cite{BC11}, we obtain a necessary and sufficient optimality criterion
for the minimizer of a diconvex functional based on the dilinear
subdifferential calculus.

\begin{Theorem}[\PFermat's rule] 
  \label{the:fermat-rule}
  Let $F \colon X \rightarrow \overline R$ be a proper functional on
  the real \PBanach\ space $X$.  Then
  $u^* \in \argmin \{F(u) : u \in X\}$ is a minimizer of\, $F$ if and
  only if\, $0 \in \uppartial_\upbeta F(u^*)$.
\end{Theorem}

\begin{Proof}
  By definition, zero is contained in the dilinear subgradient of $F$
  at $u^*$ if and only if $F(v) \ge F(u^*)$ for all $v$ in $X$, which
  is equivalent to $u^*$ being a minimizer of $F$.  \qed
\end{Proof}

Using the dilinear subgradient, we now generalize the classical
\PBregman\ distance, see for instance \cite{BC11,IJ15}, to the
dilinear/diconvex setting.

\begin{Definition}[Dilinear \PBregman\ distance]
  \label{def:dil-Breg-dist}
  Let $F \colon X \rightarrow \overline \R$ be a diconvex mapping on
  the real \PBanach\ space $X$.  The \emph{dilinear \PBregman\
    domain} $\Delta_{\upbeta, \dom} (F)$ of the mapping $F$ is the
  union of all points in $X$ with non-empty dilinear subdifferential, \ie
  \begin{equation*}
    \Delta_{\upbeta, \dom} (F) 
    \coloneqq \{ u \in X : \uppartial_\upbeta F(u) \ne \emptyset\}.
  \end{equation*}
  For every $u \in \Delta_{\upbeta, \dom} (F)$ and $(\xi, \Xi)
  \in \uppartial_\upbeta F(u)$, the \emph{dilinear \PBregman\ distance}
  of $v$ and $u$ with respect to $F$ and $(\xi, \Xi)$ is given by
  \begin{equation*}
    \Delta_{\upbeta,(\xi,\Xi)}(v,u) \coloneqq F(v) - F(u) -
    \iProd{\xi}{v-u} - \iProd{\Xi}{v \otimes v - u \otimes u}.
  \end{equation*}
\end{Definition}

\section{Dilinear inverse problems}
\label{sec:dil-inverse-probl}

During the last decades, the theory of inverse problems has become one
of the central mathematical tools for data recovery problems in
medicine, engineering, and physics.  Many inverse problems are
ill posed such that finding numerical solutions is challenging.
Although the regularization theory of linear inverse problems is well
established, especially with respect to convergence and corresponding
rates, solving non-linear inverse problems remains problematic, and
many approaches depend on assumptions and source conditions that are
difficult to verify or to validate.  Based on the tensorial lifting,
we will show that the linear regularization theory on \PBanach\ spaces
can be extended to the non-linear class of dilinear inverse problems.

To be more precise, we consider the \PTikhonov\ regularization for the
dilinear inverse problem
\begin{equation*}
  K(u) = g^\dagger,
  \addmathskip
\end{equation*}
where $K \colon X \rightarrow Y$ is a bounded dilinear operator
between the real \PBanach\ spaces $X$ and $Y$, and where $g^\dagger$
denotes the given data without noise.  In order to solve this type of
inverse problems, we study the \PTikhonov\ functional
\begin{equation}
  \label{eq:tikho-func}
  J_\alpha(u) \coloneqq \underbracket{\pNormn{K(u) -
      g^\delta}^p}_{\eqqcolon S(u)} + \alpha \, R(u),
\end{equation}
where $g^\delta$ represent a noisy version of the exact data
$g^\dagger$, where $p \ge 1$, and where $R$ is some appropriate
diconvex regularization term.  To verify the well-posedness and the
regularization properties of the \PTikhonov\ functional $J_\alpha$, we
rely on the well-established non-linear theory, see for instance
\cite{HKPS07}.  For this, we henceforth make the following assumptions,
which are based on the usual requirements for the linear case, \cf\
\cite[Assumption~3.1]{IJ15}.

\begin{Assumption}
  \label{ass:lift-prop}
  Let $X$ and $Y$ be real \PBanach\ spaces with predual $X_*$ and
  $Y_*$, where $X_*$ is separable or reflexive.  Assume that the data
  fidelity functional
  $S(\cdot) \coloneqq \pNormn{K(\cdot) - g^\delta}^p$ with the
  dilinear mapping $K \colon X \rightarrow Y$ and the non-negative,
  proper, diconvex regularization functional
  $R \colon X \rightarrow \overline \R$ satisfy:
  \begin{enumerate}[(i)]
  \item The \PTikhonov\ functional $J_\alpha$ is coercive in the sense
    that $J_\alpha (u) \rightarrow +\infty$ whenever
    $\pNormn{u} \rightarrow +\infty$.
  \item The functional $R$ is sequentially weakly$*$ lower
    semi-continuous.
  \item The dilinear operator $K$ is sequentially weakly$*$
    continuous.
  \end{enumerate}
\end{Assumption}

\begin{Remark}
  \label{rem:lift-prop}
  Since $X_*$ is a separable or reflexive \PBanach\ space, we can
  henceforth conclude that every bounded sequence $(u_n)_{n \in \N}$
  in $X$ contains a weakly$*$ convergent subsequence, see for instance
  \cite[Theorem~2.6.23]{Meg98} and \cite[Theorem~10.7]{Lax02}
  respectively.  \qed
\end{Remark}

For the non-linear regularization theory in \cite{HKPS07}, which
covers a much more general setting, the needed requirements are much
more sophisticated and comprehensive.  Therefore, we briefly verify
that \thref{ass:lift-prop} is compatible with these requirements.

\begin{Lemma}[Verification of required assumptions]
  \label{lem:check-ass}
  The requirements on the dilinear operator $K$ and the regularization
  functional $R$ made in \thref{ass:lift-prop} satisfy the
  requirements in \cite[Assumption~2.1]{HKPS07}.
\end{Lemma}

\begin{Proof}
  We verify the six required assumptions in \cite{HKPS07} step by
  step.
  \begin{enumerate}[(i)]
  \item We have to equip the \PBanach\ spaces $X$ and $Y$ with
    topologies $\tau_X$ and $\tau_Y$ weaker than the norm topology.
    Since $X$ and $Y$ have predual spaces, we simply associate the
    related weak$*$ topologies.
  \item The norm of $Y$ has to be sequentially lower semi-continuous
    with respect to $\tau_Y$, which immediately follows from the
    weak$*$ lower semi-continuity of the norm, see for instance
    \cite[Theorem~2.6.14]{Meg98}.
  \item The forward operator has to be sequentially continuous with
    respect to the topologies $\tau_X$ and $\tau_Y$, which coincides
    with \thref{ass:lift-prop}.iii.
  \item The regularization functional $R$ has to be proper and
    sequentially lower semi-continuous with respect to $\tau_X$, which
    coincides with \thref{ass:lift-prop}.ii.
  \item The domain of the forward operator is sequentially closed with
    respect to $\tau_X$, and the intersection $\dom(K) \cap \dom(R)$
    is non-empty, where $\dom(K)$ denotes the domain of definition of
    $K$.  Both assumptions are satisfied since $K$ is defined on the
    entire \PBanach\ space $X$, which is sequentially weakly$*$
    complete, see for instance \cite[Corollary~2.6.21]{Meg98}, and
    since $R$ is a proper functional.
  \item For every $\alpha > 0$ and $M > 0$, the sublevel sets
    \begin{equation*}
      \mathcal M_\alpha(M)
      \coloneqq \{ u \in X : J_\alpha(u) \le M \}
    \end{equation*}
    have to be sequentially compact with respect to $\tau_X$.  By
    \thref{ass:lift-prop}.i, the \PTikhonov\ functional $J_\alpha$ is
    coercive, which implies that the sublevel sets
    $\mathcal M_\alpha (M)$ are bounded for every $\alpha > 0$ and
    $M > 0$.  Since the functional $J_\alpha$ is weakly$*$ lower
    semi-continuous, the sublevel sets $\mathcal M_\alpha (M)$ are
    also sequentially weakly$*$ closed.  Due to the assumption that
    $X_*$ is separable or reflexive, the sequential weak$*$
    compactness follows, see \thref{rem:lift-prop}.  \qed
  \end{enumerate}
\end{Proof}

\begin{Remark}
  \label{rem:check-ass}
  Since the original proofs of the well-posedness, stability, and
  consistency only employs the sequential versions of the weak$*$
  (semi-)continuity, closure, and compactness, we have weakened these
  assumptions accordingly.  Likewise, we have skipped the unused
  convexity of the regularization functional.  \qed
\end{Remark}

Since the dilinear operator $K$ is the composition of
$u \mapsto (u, u \otimes u)$ and the lifted operator $\breve K$, the
weak$*$ continuity in \thref{ass:lift-prop}.iii may be transferred to
the representative $\breve K$.  The main problem with this approach is
that the predual of the tensor product
$X \otimes_{\uppi, \mathrm{sym}} X$ does not have to exist, even if
the predual $X_*$ of the real \PBanach\ space $X$ is known.
Therefore, we equip the symmetric tensor product
$X \otimes_{\uppi, \mathrm{sym}} X$ with an appropriate topology.

As discussed above, each finite-rank tensor in
$\omega \in X_* \otimes X_*$ defines a bilinear form
$B_\omega \colon X \times X \rightarrow \R$.  If we consider the
closure of $X_* \otimes X_*$ with respect to the norm
\begin{equation*}
  \pNorm{\omega}_\epsilonup \coloneqq \pNorm{B_\omega}
  = \sup \bigl\{ \tfrac{B_\omega(u,v)}{\pNorm{u} \pNorm{v}} : u,v \in
  X \setminus \{0\}  \bigr\} 
  \addmathskip
\end{equation*}
of $\mathcal B(X \times X)$, we obtain the \emph{injective tensor
  product} $X_* \otimes_\epsilonup X_*$, see \cite{DF93,Rya02}.  Since
the space of bilinear forms $\mathcal B(X \times X)$ is the
topological dual of $X \otimes_\uppi X$, the injective tensor product
$X_* \otimes_\epsilonup X_*$ of the predual $X_*$ is a family of
linear functionals on $X \otimes_\uppi X$.

If the \PBanach\ space $X$ has the approximation property, \ie\ the
identity operator can be approximated by finite-rank operators on
compact sets, then the canonical mapping from $X \otimes_\uppi X$ into
$(X_* \otimes_\upepsilon X_*)^*$ becomes an isometric embedding, see
\cite[Theorem~4.14]{Rya02}.  In this case, the injective tensor
product $X_* \otimes_\upepsilon X_*$ even defines a separating family,
and the projective tensor product $X \otimes_\pi X$ together with the
topology induced by the injective tensor product
$X_* \otimes_\epsilon X_*$ thus becomes a \PHausdorff\ space.  In our
setting, the separation property and the approximation property are
even equivalent, \cf\ \cite[Proposition~4.6]{Rya02}.  From a practical
point of view, most \PBanach\ spaces have the approximation property.
For example, every \PBanach\ space with \PSchauder\ basis has this
property, see \cite[Example~4.4]{Rya02}.

In the same manner, the subspace
$X_* \otimes_{\epsilonup, \mathrm{sym}} X_*$ of the symmetric tensors
$\omega$ -- of the symmetric bilinear forms $B_\omega$ -- forms a
family of linear functionals on $X \otimes_{\uppi, \mathrm{sym}} X$.
If the \PBanach\ space $X$ has the approximation property, the
symmetric injective tensor product again defines a separating
family.

On the basis of this observation, we equip the symmetric projective
tensor product $X \otimes_{\uppi, \mathrm{sym}} X$ with the weakest
topology such that each tensor $\omega$ of the symmetric injective
tensor product $X_* \otimes_{\upepsilon,\mathrm{sym}} X_*$ becomes a
continuous, linear functional.  More precisely, the weak$*$ topology
induced by $X_* \otimes_{\upepsilon, \mathrm {sym}} X_*$ is generated
by the family of preimages
\begin{equation*}
  \bigl\{ \omega^{-1}(U) : \omega \in X_* \otimes_{\upepsilon,\mathrm{sym}}
  X_*, U \subset \R, U\text{ is open} \bigr\},
\end{equation*}
see for instance \cite[Proposition~2.4.1]{Meg98}.  Since the symmetric
injective tensor product $X_* \otimes_{\upepsilon, \mathrm{sym}} X_*$
is a subspace of all linear functionals on
$X \otimes_{\uppi,\mathrm{sym}} X$, the induced topology is locally
convex, see \cite[Theorem~2.4.11]{Meg98}.  Further, a sequence
$(w_n)_{n \in \N}$ of tensors in $X \otimes_{\uppi, \mathrm{sym}} X$
converges to an element $w$ in $X \otimes_{\uppi, \mathrm{sym}} X$
with respect to the topology induced by
$X_* \otimes_{\upepsilon,\mathrm{sym}} X_*$ weakly$*$ if and only if
$\omega(w_n)$ converges to $\omega(w)$ for each $\omega$ in
$X_* \otimes_{\upepsilon, \mathrm{sym}} X_*$, see for instance
\cite[Proposition~2.4.4]{Meg98}.

The central reason to choose the injective tensor product
$X_* \otimes_{\upepsilon, \mathrm{sym}} X_*$ as topologizing family
for $X \otimes_{\uppi,\mathrm{sym}} X$ is that, under further
assumptions, like the \PHilbert\ space setting $X=H$, the injective
tensor product actually becomes a true predual of the projective
tensor product, see \cite{Rya02}.

\begin{Lemma}[Weak$*$ continuity of the tensor mapping]
  \label{lem:weak-cont-tensor-map}
  Let $X$ be a real \PBanach\ space with predual $X_*$.  The mapping
  $\otimes \colon X \rightarrow X \otimes_{\uppi, \mathrm{sym}} X$
  with $u \mapsto u \otimes u$ is sequentially weakly$*$ continuous
  with respect to the topology induced by the injective tensor product
  $X_* \otimes_{\epsilonup, \mathrm{sym}} X_*$.
\end{Lemma}

\begin{Proof}
  Let $\phi \otimes \phi$ be a rank-one tensor in
  $X_* \otimes_{\upepsilon, \mathrm{sym}} X_*$, and let
  $(u_n)_{n \in \N}$ be a weakly$*$ convergent sequence in the
  \PBanach\ space $X$.  Without loss of generality, we postulate that
  the sequence $(u_n)_{n \in \N}$ is bounded by $\pNorm{u_n} \le 1$.
  Under the assumption that $u$ is the weak$*$ limit of
  $(u_n)_{n \in \N}$, we observe
  \begin{equation*}
    \lim_{n \rightarrow \infty} (\phi \otimes \phi)(u_n \otimes u_n)
    = \lim_{n \rightarrow \infty} \phi(u_n) \, \phi(u_n)
    = \phi(u) \, \phi(u)
    = (\phi \otimes \phi)(u \otimes u).
  \end{equation*}
  Obviously, this observation remains valid for all finite-rank
  tensors in $X_* \otimes_{\upepsilon, \mathrm{sym}} X_*$.  Now, let
  $\omega$ be an arbitrary tensor in
  $X_* \otimes_{\upepsilon, \mathrm{sym}} X_*$.  For every
  $\epsilon > 0$, we find a finite-rank approximation
  $\widetilde \omega$ of the tensor $\omega$ such that
  $\pNormn{\widetilde \omega - \omega}_\upepsilon \le
  \nicefrac{\epsilon}{3}$.  Hence, for suitable large $n$, we have
  \begin{align*}
    \\[\subalignskip]
    \absn{\omega(u_n \otimes u_n) - \omega(u \otimes u)}
    &\le \absn{\omega(u_n \otimes u_n) - \widetilde \omega (u_n
      \otimes u_n)}
    \\[\fsmallskip]
    &\qquad + \absn{\widetilde \omega(u_n \otimes u_n) - \widetilde
      \omega (u \otimes u)}
    \\[\fsmallskip]
    &\qquad + \absn{\widetilde \omega (u \otimes u) - \omega(u \otimes
      u)} 
      \le \epsilon
    \\[\subalignbelowskip]
  \end{align*}
  and thus
  \begin{equation*}
    \lim_{n \rightarrow \infty} \omega (u_n \otimes u_n) = \omega(u \otimes
    u)
    \addmathskip
  \end{equation*}
  for all linear functionals $\omega$ in
  $X_* \otimes_{\upepsilon, \mathrm{sym}} X_*$.  \qed
\end{Proof}

\begin{Proposition}[Sequential weak$*$ continuity of a dilinear
  operator]
  \label{prop:seq-weak-cont-dil-oper}
  Let $X$ and $Y$ be real \PBanach\ spaces with preduals $X_*$ and
  $Y_*$.  If the representative
  $\breve K \colon X \times (X \otimes_{\uppi, \mathrm{sym}} X)
  \rightarrow Y$
  is sequentially weakly$*$ continuous with respect to the weak$*$
  topology on $X$ and the weak$*$ topology on
  $X \otimes_{\uppi, \mathrm{sym}} X$ induced by
  $X_* \otimes_{\upepsilon, \mathrm{sym}} X_*$, then the related
  dilinear operator $K \colon X \rightarrow Y$ is sequentially
  weakly$*$ continuous.
\end{Proposition}

\begin{Proof}
  Since the related dilinear operator $K$ is given by
  $K(u) \coloneqq \breve K(u, u \otimes u)$, the assertion immediately
  follows from the sequential weak$*$ continuity of $\breve K$ and the
  sequential weak$*$ continuity of the tensor mapping
  $u \mapsto (u, u \otimes u)$, \cf\ \thref{lem:weak-cont-tensor-map}.  \qed
\end{Proof}

With an entirely analogous argumentation, the required sequential
weak$*$ lower semi-continuity of the regularization functional $R$ may
be inherit from the sequential weak$*$ lower semi-continuity of the
representative $\breve R$.

\begin{Proposition}[Sequential weak$*$ lower semi-continuity of a
  diconvex mapping]
  \label{prop:seq-weak-cont-dic-map}
  Let $X$ be a real \PBanach\ space with predual $X_*$.  If the
  representative $\breve F \colon X \times (X \otimes_{\uppi,
    \mathrm{sym}} X) \rightarrow \overline \R$ is sequentially
  weakly$*$ lower semi-continuous with respect to the weak$*$ topology
  on $X$ and the weak$*$ topology on $X \otimes_{\uppi, \mathrm{sym}}
  X$ induced by $X_* \otimes_{\upepsilon, \mathrm{sym}} X_*$, then the
  related diconvex mapping $F : X \rightarrow \overline \R$ is
  sequentially weakly$*$ lower semi-continuous.
\end{Proposition}

\begin{Proof}
  Being a composition of the sequentially weakly$*$ continuous tensor
  mapping $u \mapsto (u, u \otimes u)$, \cf\
  \thref{lem:weak-cont-tensor-map}, and the sequentially weakly$*$
  lower semi-continuous representative $\breve F$, the related
  diconvex mapping $F$ given by
  $F(u) \coloneqq \breve F(u, u \otimes u)$ is obviously sequentially
  weakly$*$ continuous.  \qed
\end{Proof}

At this point, one may ask oneself whether each sequentially weakly$*$
lower semi-continuous, diconvex mapping possesses a sequentially
weakly$*$ lower semi-continuous, convex representative.  At least for
finite-dimensional spaces, where the weak$*$ convergence coincides
with the strong convergence, this is always the case.  Remembering that
all real $d$-dimensional \PBanach\ spaces are isometrically isomorphic
to $\R^d$, we can restrict our argumentation to $X = \R^d$ equipped
with the \PEuklid{ian} inner product and norm.  Further, the
sequential weak$*$ lower semi-continuity here coincides with the lower
semi-continuity.  The projective tensor product
$\R^d \otimes_{\uppi,\mathrm{sym}} \R^d$ becomes the space of
symmetric matrices $\R^{d \times d}_{\mathrm{sym}}$ equipped with the
\PHilbert-\PSchmidt\ inner product and the \PFrobenius\ norm.
Moreover, the space
\raisebox{0pt}[0pt][0pt]{$\R^{d \times d}_{\mathrm{sym}}$} is spanned
by the rank-one tensors $u \otimes u = u \, u^\T$ with $u \in \R^d$.
The dual space $(\R^d \otimes_{\uppi, \mathrm{sym}} \R^d)^*$ may also
be identified with the space of symmetric matrices
$\R^{d \times d}_{\mathrm{sym}}$.

\begin{Theorem}[Lower semi-continuity in finite dimensions]
  \label{the:lift-lsc-map}
  A diconvex mapping $F \colon \R^d \rightarrow \overline \R$ is lower
  semi-continuous if and only if there exists a lower semi-continuous
  representative
  $\breve F \colon \R^d \times \R^{d \times d}_{\mathrm{sym}}
  \rightarrow \overline \R$.
\end{Theorem}

\begin{Proof} 
  In \thref{prop:seq-weak-cont-dic-map}, we have in particular shown
  that the lower semi-continuity of $\breve F$ implies the lower
  semi-continuity of $F$.  Thus, it only remains to prove that each
  lower semi-continuous diconvex functional possesses a lower
  semi-continuous representative $\breve F$.

  The assertion is obviously true for the constant functional
  $F \equiv + \infty$ with representative $\breve F(u,w) \coloneqq 0$
  for one point $(u,w)$ with $w \ne u \otimes u$ and
  $\breve F(\cdot,\cdot) = +\infty$ otherwise.  For the remaining
  functionals, the central idea of the proof is to show that the lower
  semi-con\-tin\-u\-ous convexification
  $\overline{\conv} \, F_\otimes$ -- the closure of the convex hull
  $\conv F_\otimes$ -- of the mapping $F_\otimes$ in
  \eqref{eq:F-otimes} is a valid representative of a lower
  semi-continuous mapping $F$, which means that
  $\overline \conv \, F_\otimes (u, u \otimes u) = F(u)$ for all
  $u \in \R^d$.  For the sake of simplicity, we assume that
  $F(u) \ge 0$ and thus $F_\otimes(u,u \otimes u) \ge 0$ for all $u$
  in $\R^d$, which can always be achieved by subtracting a
  (continuous) dilinear minorant, see \thref{rem:diaff-minor:1}.

  For a fixed point $(u,u \otimes u)$ on the diagonal, we 
  distinguish the following three cases:
  \begin{enumerate}[(i)]
  \item The point $(u,u \otimes u)$
    is not contained in the relative closure of
    $\dom(\overline \conv \, F_\otimes)$, which implies
    $F(u) = \overline \conv \, F_\otimes(u, u \otimes u) = +\infty$.
  \item The point $(u,u \otimes u)$ lies in the relative interior of
    $\dom(\overline \conv \, F_\otimes)$.  Since the effective domain
    $\dom(\overline \conv \, F_\otimes)$ is a subset of
    $\conv\{(u, u \otimes u) : u \in \R^d\}$, and since
    $(u,u \otimes u)$ is an extreme point of the latter set, see
    \thref{lem:extr-points-diag}, the point $(u,u \otimes u)$ thus
    has to be extreme with respect to the effective domain.  Since the
    extreme points of a convex set are, however, contained in the
    relative boundary except for the zero-dimensional case, see
    \cite[Corollary~18.1.3]{Roc70}, the effective domains of
    $\overline \conv \, F_\otimes$ and $F_\otimes$ have to consist
    exactly of the considered point $(u,u \otimes u)$.  In this
    instance, the closed convex hull $\overline \conv \, F_\otimes$
    equals $F_\otimes$ and is trivially a sequentially weakly$*$
    continuous representative.
  \item The point $(u, u \otimes u)$ is contained in the relative
    boundary of $\dom(\overline \conv \, F_\otimes)$.
  \end{enumerate}
  To finish the proof, we have to show
  $F(u) = \overline \conv \, F_\otimes(u, u \otimes u)$ for the third
  case.

  In order to compute
  $\overline \conv \, F_\otimes (u , u \otimes u)$, we apply
  \cite[Theorem~7.5]{Roc70}, which implies
  \begin{equation*}
    \overline \conv \, F_\otimes (u, u \otimes u) 
    = \lim_{\lambda \nearrow 1} \conv F_\otimes \bigl(\lambda (u, u \otimes
    u) + (1-\lambda) (v,w)\bigr),
  \end{equation*}
  where $(v,w)$ is some point in the non-empty relative interior of
  $\dom(\conv F_\otimes)$.  Next, we take a sequence
  $(v_k, w_k) \coloneqq \lambda_k (u, u \otimes u) + (1-\lambda_k)
  (v,w)$
  with $\lambda_k \in (0,1)$ so that
  $\lim_{k\rightarrow \infty} \lambda_k = 1$ and consider the limit of
  the function values $\conv \, F_\otimes (v_k, w_k)$.  Since the
  complete sequence $(v_k,w_k)_{k=1}^\infty$ is contained in the
  relative interior of $\dom(\conv F_\otimes)$, \cf\
  \cite[Theorem~6.1]{Roc70}, all functions values
  $\conv \, F_\otimes (v_k, w_k)$ are finite and can be approximated
  by \PCaratheodory's theorem.  More precisely, for every $\rho > 0$
  and for every $k\in \N$, we can always find convex combinations
  \raisebox{0pt}[0pt][0pt]{$(v_k, w_k) = \sum_{n=1}^{N+1}
    \alpha_n^{(k)} \, (u_n^{(k)}, u_n^{(k)} \otimes u_n^{(k)})$},
  where \raisebox{0pt}[0pt][0pt]{$\alpha_n^{(k)}$} is in $[0,1]$ so
  that
  \raisebox{0pt}[0pt][0pt]{$\sum_{n=1}^{N+1} \alpha_n^{(k)} = 1$}, and
  where $N$ is an integer not greater than the dimension of
  \raisebox{0pt}[0pt][0pt]{$\R^d \times \R^{d \times
      d}_{\mathrm{sym}}$}, such that
  \begin{equation*}
    \absbb{\conv \, F_\otimes (v_k, w_k)
      - \sum_{n=1}^{N+1} \alpha_n^{(k)} \, F_\otimes (u_n^{(k)}, u_n^{(k)}
    \otimes u_n^{(k)})} \le \rho,
  \end{equation*}
  \cf\ \cite[Corollary~17.1.5]{Roc70}.

  In the next step, we examine the occurring sequence of convex
  combinations in more detail.  For this purpose, we define the half
  spaces
  \begin{equation*}
    \mathbb H^+_{u,\epsilon} \coloneqq \bigl\{ (v,w) : G_u(v,w) \ge
    \epsilon \bigr\} 
    \qquad\text{and}\qquad
    \mathbb H^-_{u,\epsilon} \coloneqq \bigl\{ (v,w) : G_u(v,w) \le
    \epsilon  \bigr\},
  \end{equation*}
  with $u \in \R^d$, $\epsilon > 0$, and the functional
  $G_u \colon \R^d \times \R^{d \times d}_{\mathrm{sym}} \rightarrow \R$
  given by
  \begin{equation*}
    G_u(v,w) \coloneqq - \iProdn{2u}{v}_{\R^d} + \iProdn{I}{w}_{\R^{d
        \times d}_{\mathrm{sym}}},
  \end{equation*}
  where $I$ denotes the identity matrix.  Obviously, a vector
  $(v, v \otimes v)$ is contained in the shifted half space
  $\mathbb H^+_{u,\epsilon} + (u, u \otimes u)$ if and only if
  \begin{equation}
    \label{eq:char-H-epsilon}
    \begin{aligned}
      G_u(v-u,v \otimes v - u \otimes u) 
      &=- \iProd{2u}{v-u}_{\R^d} +
      \iProd{I}{v \otimes v}_{\R^{d \times d}_{\mathrm{sym}}} -
      \iProd{I}{u \otimes u}_{\R^{d \times d}_{\mathrm{sym}}}
      \\[\fskip]
      &=- \iProd{2u}{v-u}_{\R^d} + \iProd{v}{v}_{\R^d} -
      \iProd{u}{u}_{\R^d}
      \\[\fskip]
      &= \pNorm{v-u}^2 \ge \epsilon.
    \end{aligned}
  \end{equation}
  Similarly, the vector $(v, v \otimes v)$ is contained in
  $\mathbb H^-_{u,\epsilon} + (u,u \otimes u)$ if and only if
  $\pNorm{v-u}^2 \le \epsilon$.

  As mentioned above, we now consider a sequence of convex
  combinations
  \begin{equation*}
    \sum_{n=1}^{N+1} \alpha_n^{(k)} \, (u_n^{(k)}, u_n^{(k)} \otimes
    u_n^{(k)})
    \xrightarrow[k\rightarrow \infty]{} (u, u \otimes u),
  \end{equation*}
  where $\alpha_n^{(k)} \in [0,1]$ so that
  $\sum_{n=1}^{N+1} \alpha_n^{(k)} = 1$, and where $N$ is some fixed
  integer independent of $k$.  For $n=1$, either the sequence
  \raisebox{0pt}[0pt][0pt]{$(u_1^{(k)})_{k \in \N}$} has a subsequence
  converging to $u$ or there exists an $\epsilon_1 > 0$ such that
  \raisebox{0pt}[0pt][0pt]{$\pNormn{u_1^{(k)} - u}^2 \ge \epsilon_1$}
  for all $k$ in $\N$.  In the second case, the complete sequence
  \raisebox{0pt}[0pt][0pt]{$(u_1^{(k)})_{k \in \N}$} is contained in
  the shifted half space
  $\mathbb H_{u,\epsilon_1}^+ + (u,u \otimes u)$.  Thinning out the
  sequence of convex combinations by repeating this construction for
  the remaining indices $n=2,\dots,N$ iteratively, we obtain a
  subsequence of the form
  \begin{equation*}
    \sum_{n=1}^L \alpha_n^{(\ell)} \, \underbracket{(u_n^{(\ell)}, u_n^{(\ell)}
      \otimes u_n^{(\ell)})}_{\rightarrow (u, u \,\otimes\, u)}
    + \sum_{n=L+1}^{N+1} \alpha_n^{(\ell)} \, \underbracket{(u_n^{(\ell)}, u_n^{(\ell)}
      \otimes u_n^{(\ell)})}_{\in \mathbb H_{u,\epsilon_n}^+ (u, u
      \,\otimes\, u)}
    \xrightarrow[\ell \rightarrow \infty]{} (u, u \otimes u),
  \end{equation*}
  where we assume that the first $L$ sequences
  \raisebox{0pt}[0pt][0pt]{$(u_n^{(\ell)})_{\ell \in \N}$} possess the
  accumulation point $u$ without loss of generality -- if necessary,
  we rearrange the indices $n=1, \dots, N+1$ accordingly.  Thinning
  out the subsequence even further, we can also ensure that the
  coefficients
  \raisebox{0pt}[0pt][0pt]{$(\alpha_n^{(\ell)})_{\ell \in \N}$}
  converge for every index $n$.

  In the following, we have to take special attention of the
  subsequence \raisebox{0pt}[0pt][0pt]{$(u_n^{(\ell)})_{\ell \in \N}$}
  not converging to $u$.  Therefore, we consider the case where
  \raisebox{0pt}[0pt][0pt]{$\beta^{(\ell)} \coloneqq
    \sum_{n=L+1}^{N+1} \alpha_n^{(\ell)}$}
  does not become constantly zero after some index $\ell_0$ more
  precisely.  Taking a subsequence with $\beta^{(\ell)} \ne 0$, and
  re-weighting the sequence of convex combinations, we now obtain
  \begin{equation*}
    (1-\beta^{(\ell)})
     \! \underbracket{\sum_{n=1}^L \tfrac{\alpha_n^{(\ell)}}{1-\beta^{(\ell)}} \,
    (u_n^{(\ell)}, u_n^{(\ell)} \otimes u_n^{(\ell)})}_{\rightarrow (u, u \,\otimes\, u)}
    + \beta^{(\ell)} \!\! \underbracket{
    \sum_{n=L+1}^{N+1} \tfrac{\alpha_n^{(\ell)}}{\beta^{(\ell)}} \,
    (u_n^{(\ell)}, u_n^{(\ell)} 
      \otimes u_n^{(\ell)})}_{\in \mathbb H_{u,\epsilon}^+ + (u, u
      \,\otimes\, u)}
    \xrightarrow[\ell \rightarrow \infty]{} (u, u \otimes u),
  \end{equation*}
  where $\epsilon > 0$ is chosen smaller than
  $\epsilon_{L+1}, \dots, \epsilon_{N+1}$.  Since the second sum is a
  convex combination of points in
  $\mathbb H_{u,\epsilon}^+ + (u, u \,\otimes\, u)$, the value of the
  sum is also contained in
  $\mathbb H_{u,\epsilon}^+ + (u, u \,\otimes\, u)$. Since
  $(u, u \otimes u)$ is not contained in the closed set
  $\mathbb H_{u,\epsilon}^+ + (u,u \otimes u)$ for $\epsilon > 0$ by
  construction, the coefficients $\beta^{(\ell)}$ could neither become
  constantly one.  For an appropriate subsequence, the re-weighted
  convex combinations are thus well defined.  Obviously, the first sum
  converges to $(u,u \otimes u)$.  If we now assume that the sequence
  \raisebox{0pt}[0pt][0pt]{$(\beta^{(\ell)})_{\ell \in \N}$} does not
  converge to zero, the second sum has to converge to some point in
  $\mathbb H_{u,\epsilon}^+ +(u,u \otimes u)$.  Consequently,
  $(u, u \otimes u)$ is a non-trivial convex combination of itself
  with a further point in
  $\mathbb H_{u,\epsilon}^+ + (u,u \otimes u)$.  Since
  $(u, u \otimes u)$ is not contained in
  $\mathbb H_{u,\epsilon}^+ + (u, u \otimes u)$, see
  \eqref{eq:char-H-epsilon}, this is not possible.  Hence, the
  coefficient $\beta^{(\ell)}$ converges to zero, which is our main
  observation in this step.

  Applying the subsequence construction above to the sequence of
  function values
  \begin{equation*}
    \biggl(\sum_{n=1}^{N+1} \alpha_n^{(k)} \, F_\otimes (u_n^{(k)}, u_n^{(k)}
    \otimes u_n^{(k)}) \biggr)_{k=1}^\infty,
  \end{equation*}
  and exploiting the lower semi-continuity and non-negativity of $F$, we
  can finally estimate the limit of the function values
  $\conv F_\otimes (v_k, w_k)$ by
  \begin{align*}
    \lim_{\ell \rightarrow \infty} \conv F_\otimes(v_\ell, w_\ell) + \rho
    &\ge \liminf_{\ell\rightarrow \infty}  \sum_{n=1}^{N+1} \alpha_n^{(\ell)} \,
    F_\otimes (u_n^{(\ell)}, u_n^{(\ell)} \otimes u_n^{(\ell)})
    \\[\fskip]
    &\ge \liminf_{\ell \rightarrow \infty}\ (1-\beta^{(\ell)}) \, \sum_{n=1}^L
    \tfrac{\alpha_n^{(\ell)}}{1 - \beta^{(\ell)}} \,
    F_\otimes(u_n^{(\ell)}, u_n^{(\ell)} \otimes u_n^{(\ell)})
    \\[\fskip]
    &\ge \sum_{n=1}^L \liminf_{\ell \rightarrow \infty} \bigl[
      \tfrac{\alpha_n^{(\ell)}}{1 - \beta^{(\ell)}} \, F(u_n^{(\ell)})
      \bigr]
      \ge F(u)
  \end{align*}
  because $1 - \beta^{(\ell)} = \sum_{n=1}^L \alpha_n^{(\ell)}$
  converges to one as discussed above.  Since the accuracy $\rho$ of
  the approximation can be chosen arbitrarily small, we thus have
  \begin{equation*}
    \overline \conv \, F_\otimes (u, u \otimes u) 
    = \lim_{k \rightarrow \infty} \conv F_\otimes (v_k, w_k) 
    \ge F(u).
  \end{equation*}
  Hence, $\overline \conv F_\otimes (u, u \otimes u)$ equals $F(u)$
  for all $u \in \R^d$, and the lower semi-continuous convex hull
  $\overline{\conv} \, F_\otimes$ is a valid representative of $F$.
  \qed
\end{Proof}

\begin{Remark}
  \label{rem:lift-lsc-map}
  Using an analogous argumentation, one may extend
  \thref{the:equa-subdiff} to the lower semi-continuous convex hull
  $\overline \conv \, F_\otimes$.  More precisely, if
  $F \colon \R^d \rightarrow \R$ is a lower semi-continuous, diconvex
  mapping, then the representative subdifferential with respect to
  $\overline \conv \, F_\otimes$ and the dilinear subdifferential
  coincide, \ie
  \begin{equation*}
    \breve \uppartial F(u) = \uppartial_\upbeta F(u).
    \tag*{\qed}
  \end{equation*}
\end{Remark}

\subsection{Well-posedness and regularization properties}
\label{sec:well-posedn-regul}

We now return to the well-posedness, stability, and consistency of the
\PTikhonov\ regularization of the dilinear inverse problem
$K(u) = g^\dagger$.  In other words, we study the regularization
properties of the variational regularization
\begin{equation*}
  \operatorname*{minimize}_{u \in X} J_\alpha (u)
  \qquad\text{with}\qquad
  J_\alpha (u) = \pNormn{K(u) - g^\delta}^p + \alpha \, R(u).
\end{equation*}
Since the introduced dilinear regularization is particularly an
instance of the non-linear \PTikhonov\ regularization in \PBanach\
spaces, see for instance \cite{HKPS07}, the well-posed\-ness,
stability, and consistency immediately follow from the
well-established non-linear regularization theory.  For the sake of
completeness and convenience, we briefly summarize the central results
with respect to our setting.

Firstly, the \PTikhonov\ functional $J_\alpha$ in
\eqref{eq:tikho-func} is well posed in the sense that the minimum of
the regularized problem $\min \{ J_\alpha (u) : u \in X \}$ is
attained, and that the related minimizer is thus well defined.

\begin{Theorem}[Well-posedness]
  \label{the:well-posed-tikho}
  Under \thref{ass:lift-prop}, for every $\alpha > 0$, there exists at
  least one minimizer $u_\alpha^\delta$ to the functional $J_\alpha$
  in \eqref{eq:tikho-func}.
\end{Theorem}

\begin{Proof}
  The well-posedness immediately follows from
  \cite[Theorem~3.1]{HKPS07} since the required assumptions are
  fulfilled by \thref{lem:check-ass}.  \qed
\end{Proof}

Stability of a variational regularization method means that the
minimizers $u_\alpha^\delta$ of the \PTikhonov\ functional $J_\alpha$
weakly$*$ depend on the noisy data $g^\delta$.  If the regularization
functional $R$ satisfies the so-called H-property, see for instance
\cite{IJ15,Wer02}, the dependence between the solution of the
regularized problem and the corrupted data is actually strong.

\begin{Definition}[H-property]
  \label{def:H-property}
  A functional $R \colon X \rightarrow \overline \R$ possesses the
  \emph{H-property on the space $X$} if any weakly$*$ convergent
  sequence $(u_n)_{n \in \N}$ in $X$ with limit $u$ and with
  $R(u_n) \rightarrow R(u)$ converges to $u$ strongly.
\end{Definition}

\begin{Theorem}[Stability]
  \label{the:stab-tikho}
  Let the sequence $(g_n)_{n \in \N}$ in $Y$ be convergent with limit
  $g^\delta \in Y$.  Under \thref{ass:lift-prop}, the sequence
  $(u_n)_{n \in \N}$ of minimizers of the \PTikhonov\ functional
  $J_\alpha$ in \eqref{eq:tikho-func} with $g_n$ in place of
  $g^\delta$ contains a weakly$*$ convergent subsequence to a
  minimizer $u_\alpha^\delta$ of $J_\alpha$.  If the minimizer of
  $J_\alpha$ is unique, then the complete sequence $(u_n)_{n \in \N}$
  converges weakly$*$.  If the functional $R$ possesses the
  H-property, then the sequence $(u_n)_{n \in \N}$ converges in norm
  topology.
\end{Theorem}

\begin{Proof}
  The existence of a weakly$*$ convergent subsequence, whose limit is
  a minimizer of $J_\alpha$, follows directly from the stability of
  the non-linear \PTikhonov\ regularization
  \cite[Theorem~3.2]{HKPS07}.  Again, the required assumptions are
  fulfilled by \thref{lem:check-ass}.  Moreover, if the minimizer
  $u_\alpha^\delta$ is unique, each subsequence of $(u_n)_{n \in \N}$
  contains a subsequence weakly$*$ converging to $u_\alpha^\delta$,
  and hence, the entire sequence $(u_n)_{n \in \N}$ converges
  weakly$*$, \cf\ \cite[Theorem~3.2]{IJ15}.

  Completely analogously to the proof of \cite[Theorem~3.2]{IJ15}, one
  can now show that the sequence $R(u_n)$ must converge to $R(u^*)$.
  In this case, the H-property additionally yields the convergence of
  $(u_n)_{n \in \N}$ in norm topology.  \qed
\end{Proof}

Finally, the \PTikhonov\ regularization of a dilinear inverse problem
is consistent; so the minimizer $u_\alpha^\delta$ weakly$*$ converges
to a solution $u^\dagger$ of the unperturbed problem $K(u)=g^\dagger$
if the noise level $\delta$ goes to zero.  More precisely, the
solution $u_\alpha^\delta$ converges to an $R$-minimizing solution,
see for instance \cite[Definition~3.3]{HKPS07} or
\cite[Definition~3.2]{IJ15} for the definition.

\begin{Definition}[$R$-minimizing solution]
  \label{eq:R-minim-sol}
  A point $u^\dagger \in X$ is an \emph{$R$-minimizing solution} to the
  problem $K(u) = g^\dagger$ if it satisfies
  \begin{equation*}
    K(u^\dagger) = g^\dagger 
    \qquad \text{and} \qquad
    R(u^\dagger) \le R(u)
    \addmathskip
  \end{equation*}
  for all further solutions $u$ of the dilinear equation
  $K(u) = g^\dagger$.
\end{Definition}

\begin{Theorem}[Existence of an $R$-minimizing solution]
  \label{the:R-minim-sol}
  Under \thref{ass:lift-prop}, there exists at least one
  $R$-minimizing solution to the dilinear problem $K(u)=g^\dagger$.
\end{Theorem}

\begin{Proof}
  The existence of an $R$-minimizing solution immediately follows from
  \cite[Theorem~3.4]{HKPS07}.  The needed assumptions are satisfied by
  \thref{lem:check-ass}.  \qed
\end{Proof}

\begin{Theorem}[Consistency]
  \label{the:consistency}
  Let $(g^{\delta_n})_{n \in \N}$ be a sequence of noisy data with
  $\delta_n \coloneqq \pNormn{g^\dagger - g^{\delta_n}} \rightarrow
  0$.
  Under \thref{ass:lift-prop}, the sequence of minimizers
  $(u_{\alpha_n}^{\delta_n})_{n \in \N}$ contains a weakly$*$
  convergent subsequence whose limit is an $R$-minimizing solution
  $u^\dagger$ if the sequence of regularization parameters
  $(\alpha_n)_{n \in \N} = (\alpha(\delta_n))_{n \in \N}$ satisfies
  \begin{equation*}
    \lim_{n \rightarrow \infty} \tfrac{\delta_n^p}{\alpha_n} = 0
    \qquad \text{and} \qquad
    \lim_{n \rightarrow \infty} \alpha_n = 0.
  \end{equation*}
  If the $R$-minimizing solution $u^\dagger$ is unique, then the
  entire sequence $(u_{\alpha_n}^{\delta_n})_{n \in \N}$ converges
  weakly$*$.  If the functional $R$ possesses the H-property, then the
  sequence $(u_{\alpha_n}^{\delta_n})_{n \in \N}$ converges in norm
  topology.
\end{Theorem}

\begin{Proof}
  Due to the consistency property of the non-linear \PTikhonov\
  regularization, see \cite[Theorem~3.5]{HKPS07}, the sequence of
  minimizers $(u_{\alpha_n}^{\delta_n})_{n \in \N}$ possesses a
  subsequence converging to an $R$-minimizing solution weakly$*$.
  Once more, the required assumption in \cite{HKPS07} are implied by
  \thref{lem:check-ass}.

  If the $R$-minimizing solution $u^\dagger$ is unique, each
  subsequence of $(u_{\alpha_n}^{\delta_n})_{n \in \N}$ has a
  weakly$*$ converging subsequence.  Consequently, the entire
  sequence
  \raisebox{0pt}[0pt][0pt]{$(u_{\alpha_n}^{\delta_n})_{n \in \N}$}
  converges weakly$*$.  As verified in the proof of
  \cite[Theorem~3.5]{HKPS07}, the sequence $R(u_n)$ converges to
  $R(u^*)$.  Under the assumption that the regularization functional
  $R$ has the H-property, the convergence of
  $(u_{\alpha_n}^{\delta_n})_{n \in \N}$ is thus strong.  \qed
\end{Proof}

\subsection{Convergence analysis under source-wise representations}
\label{sec:conv-source-repr}

Based on the dilinear subdifferential calculus, we now analyse the
convergence behaviour of the variational regularization method for
dilinear inverse problems $K(u) = g^\dagger$ if the noise level
$\delta$ goes to zero.  In the following, we assume that the dilinear
operator $K \colon X \rightarrow Y$ maps from a real \PBanach\ space
into a real \PHilbert\ space $Y$.  Further, we restrict ourselves to the
squared \PHilbert\ norm $\pNormn{\cdot}^2$ as data fidelity functional
$S$ in \eqref{eq:tikho-func}, which means that we consider the
\PTikhonov\ functional
\begin{equation}
  \label{eq:tikho-func-hilbert}
  J_\alpha (u) \coloneqq \pNormn{K(u)-
      g^\delta}^2 + \alpha \, R(u).
\end{equation}

Looking back at \thref{the:R-minim-sol}, we recall that the dilinear
inverse problem $K(u) = g^\dagger$ always possesses an $R$-minimizing
solution $u^\dagger$ with respect to the regularization functional
$R$.  Consequently, \PFermat's rule (\thref{the:fermat-rule}) implies
that zero must be contained in the dilinear subdifferential
\raisebox{0pt}[0pt][0pt]{$\uppartial_\upbeta (R +
  \chi_{\{K(u)=g^\dagger\}})(u^\dagger)$},
where the indicator function $\chi_{\{K(u)=g^\dagger\}}(u)$ is $0$ if
$u$ is a solution of the inverse problem $K(u)=g^\dagger$ and
$+\infty$ else.  Applying the sum and chain rule in
\thref{prop:sum-rule} and \ref{prop:chain-rule:con-fun}, we have
\begin{equation*}
  \ran \breve K^* + \uppartial_\upbeta R(u^\dagger)
  \subset \uppartial_\upbeta (R + \chi_{\{K(u)=g^\dagger\}})(u^\dagger).
\end{equation*}
Without further assumptions, the sum and chain rule here only yield an
inclusion.  Although we always have
$0 \in \uppartial_\upbeta (R + \chi_{\{K(u)=g^\dagger\}})(u^\dagger)$,
we hence cannot guarantee
$0 \in \ran \breve K^* + \uppartial_\upbeta R(u^\dagger)$ or,
equivalently,
$\ran \breve K^* \cap \uppartial_\upbeta R(u^\dagger) \ne \emptyset$.
Against this background, we postulate the regularity assumption
that the range of the adjoint operator and the dilinear
subdifferential are not disjoint.  In other words, we assume the
existence of a source-wise representation
\begin{equation}
  \label{eq:source-cond}
  \breve K^* \omega 
  = (\xi^\dagger, \Xi^\dagger) 
  \in \uppartial_\upbeta R (u^\dagger)
  \submathskip
\end{equation}
for some $\omega$ in $Y$.

\begin{Theorem}[Convergence rate]
  \label{the:conv-rate}
  Let $u^\dagger$ be an $R$-minimizing solution of the dilinear
  inverse problem $K(u)=g^\dagger$.  Under the source condition
  $\breve K^* \omega = (\xi^\dagger, \Xi^\dagger) \in
  \uppartial_\upbeta R(u^\dagger)$
  for some $\omega$ in $Y$ and under \thref{ass:lift-prop}, the
  minimizers $u_\alpha^\delta$ of the \PTikhonov\ regularization
  $J_\alpha$ in \eqref{eq:tikho-func-hilbert} converges to $u^\dagger$
  in the sense that the dilinear \PBregman\ distance between
  $u_\alpha^\delta$ and $u^\dagger$ with respect to the regularization
  functional $R$ is bounded by
  \begin{equation*}
    \Delta_{\upbeta,(\xi^\dagger, \Xi^\dagger)}(u_\alpha^\delta,
    u^\dagger) \le \left( \tfrac{\delta}{\sqrt \alpha} +
      \tfrac{\sqrt \alpha}{2} \pNorm{\omega} \right)^2
  \end{equation*}
  and the data fidelity term by
  \begin{equation*}
    \pNormn{K (u_\alpha^\delta) - g^\delta} 
    \le \delta + \alpha \pNorm{\omega}.
  \end{equation*}
\end{Theorem}

\begin{Proof}
  Inspired by the proof for the usual subdifferential in \cite{IJ15},
  the desired convergence rate for the dilinear subdifferential can be
  established in the following manner.  Since $u_\alpha^\delta$ is a
  minimizer of the \PTikhonov\ functional $J_\alpha$ in
  \eqref{eq:tikho-func-hilbert}, $u_\alpha^\delta$ and $u^\dagger$
  satisfy
  \begin{equation*}
    \pNormn{K(u_\alpha^\delta) - g^\delta}^2 + \alpha \, R
    (u_\alpha^\delta)
    \le \pNormn{K(u^\dagger) - g^\delta}^2 + \alpha \, R (u^\dagger).
    \addmathskip
  \end{equation*}
  Remembering $K(u^\dagger) = g^\dagger$, we can bound the norm on the
  right-hand side by $\pNormn{g^\dagger - g^\delta}^2 \le \delta^2$.
  Rearranging the last inequality and exploiting the source condition,
  we get
  \begin{align*}
    &\pNormn{K(u_\alpha^\delta) - g^\delta}^2 + \alpha \,
    \Delta_{\upbeta,(\xi^\dagger,\Xi^\dagger)}(u_\alpha^\delta, u^\dagger)
    \\[\fskip]
    &\qquad\le \delta^2 - \alpha \iProd{\xi^\dagger}{u_\alpha^\delta -
      u^\dagger} - \alpha \iProd{\Xi^\dagger}{u_\alpha^\delta \otimes
      u_\alpha^\delta - u^\dagger \otimes u^\dagger}
    \\[\fskip]
    &\qquad= \delta^2 - \alpha \iProd{\omega}{\breve
      K(u_\alpha^\delta - u^\dagger, u_\alpha^\delta \otimes
      u_\alpha^ \delta - u^\dagger \otimes u^\dagger)}
    \\[\fskip]
    &\qquad= \delta^2 - \alpha \iProd{\omega}{K(u_\alpha^\delta) -
      g^\dagger}
    \\[\fskip]
    &\qquad= \delta^2 - \alpha \iProd{\omega}{K(u_\alpha^\delta) -
      g^\delta} - \alpha \iProd{\omega}{g^\dagger - g^\delta}.
  \end{align*}
  Rearranging the terms, completing the square, and applying
  \PCauchy-\PSchwarz's inequality, we obtain
  \begin{equation*}
    \pNorm{K(u_\alpha^\delta) - g^\delta + \tfrac{\alpha}{2} \omega}^2
    + \alpha \, \Delta_{\upbeta,(\xi^\dagger,\Xi^\dagger)}(u_\alpha^\delta, u^\dagger)
    \le \delta^2 + \alpha \delta \pNorm{\omega} +
    \tfrac{\alpha^2}{4} \pNorm{\omega}^2
    = \left( \delta + \tfrac{\alpha}{2} \pNorm{\omega} \right)^2,
  \end{equation*}
  which proves the convergence rate for the \PBregman\ distance.  The
  second convergence rate follows immediately by applying the reverse
  triangle inequality.  \qed
\end{Proof}

\begin{Remark}
  \label{rem:conv-rate}
  If we apply the a priori parameter choice rule $\alpha \sim \delta$,
  then the \PBregman\ distance
  $\Delta_{\upbeta,(\xi^\dagger, \Xi^\dagger)}(u_\alpha^\delta,
  u^\dagger)$
  as well as the data fidelity term
  $\pNormn{K(u_\alpha^\delta) - g^\delta}$ converges to zero with a
  rate of $\Landau(\delta)$.  \qed
\end{Remark}

\section{Deautoconvolution problem}
\label{sec:deaut-probl}

In order to give a non-trivial example for the practical relevance of
the developed dilinear regularization theory, we consider the
deautoconvolution problem, where one wishes to recover an unknown signal
$u \colon \R \rightarrow \C$ with compact support from its
kernel-based autoconvolution
\submathskip
\begin{equation}
  \label{eq:autoconv}
  \mathcal A_k[u](t) \coloneqq \int_{-\infty}^\infty k(s,t) \, u(s) \,
  u(t-s) \diff s 
  \qquad(t \in \R),
\end{equation}
where $k \colon \R^2 \rightarrow \C$ is an appropriate kernel
function.  Problems of this kind occur in spectroscopy, optics, and
stochastics for instance, see \cite{BH15a}.

Following the model in \cite{BH15a}, we assume that $u$ is a
square-integrable complex-valued signal on the interval $[0,1]$, or,
in other words, $u \in L^2_\C([0,1])$.  To ensure that the integral in
\eqref{eq:autoconv} is well defined, we extend the signal $u$ outside
the interval $[0,1]$ with zero and restrict ourselves to bounded
kernels $k \in L_\C^\infty([0,1] \times [0,2])$.  Considering the
support of $u$, we may moreover assume
\begin{equation}
  \label{eq:ass-ker-supp}
  \supp k \subset \bigl\{ (s,t) : 0 \le s \le 1 ~\text{and}~ s \le t \le
  s+1 \bigr\} 
  \submathskip
\end{equation}
and
\begin{equation}
  \label{eq:ass-ker-sym}
  k(s,t) = k(t-s,t)
  \qquad(0 \le s \le 1, s \le t \le s+1).
  \addmathskip
\end{equation}
The symmetry property \eqref{eq:ass-ker-sym} can be demanded in
general because of the identity
\begin{equation*}
  \mathcal A_k [u] (t)
  = \int_{-\infty}^\infty \tfrac{k(s,t)+k(t-s,t)}{2} \, u(s) \, u(t-s)
  \diff s.
\end{equation*}

After these preliminary considerations, we next verify that the
kernel-based autoconvolution $\mathcal A_k$ is a bounded dilinear
operator.  For this, we exploit that the autoconvolution
$\mathcal A_k$ is the quadratic mapping related to the symmetric
bilinear mapping
$\mathcal B_k \colon L_\C^2([0,1]) \times L_\C^2([0,1]) \rightarrow
L_\C^2([0,2])$ given by
\begin{equation*}
  \mathcal B_k[u,v](t)
  \coloneqq
  \int_{-\infty}^\infty k(s,t) \, u(s) \, v(t-s) \diff s 
  \qquad(t \in \R).
\end{equation*}
The well-definition and boundedness of $\mathcal B_k$ immediately
follows from \PJensen's inequality.  More precisely, we obtain
\begin{align*}
  \pNorm{\mathcal B_k[u,v]}^2_{L^2_\C([0,2])}
  &\le \int_{-\infty}^\infty \int_{-\infty}^\infty \bigl[ k(s,t) \,
    u(s) \, v(t-s) \bigr]^2 \diff s \diff t 
  \\[\fskip]
  &\le \pNorm{k}^2_{L^\infty_\C([0,1] \times [0,2])}
    \pNorm{u}^2_{L^2_\C([0,1])} \pNorm{v}^2_{L^2_\C([0,1])},
\end{align*}
In connection with \thref{ex:quad-map}, this shows that the
kernel-based autoconvolution $\mathcal A_k$ is a bounded dilinear
mapping from $L^2_\C([0,1])$ into $L^2_\C([0,2])$.  In order to apply
the developed dilinear regularization theory, we interpret the complex
\PHilbert\ spaces $L^2_\C([0,1])$ and $L^2_\C([0,2])$ as real
\PHilbert\ spaces with the inner product
$\iProd{\cdot}{\cdot}_\R \coloneqq \Re \iProd{\cdot}{\cdot}_\C$.

Although the deautoconvolution problem is ill posed in general, the
unperturbed inverse problem can have at most two different solutions,
see \cite{GHB+14}.  The original proof is based on \PTitchmarsh's
convolution theorem.

\begin{Proposition}[Ambiguities deautoconvolution]
  Let $g^\dagger$ be a function in $L^2_\C([0,2])$, and let the kernel
  of the autoconvolution be $k \equiv 1$.  If $u \in L^2_\C([0,1])$ is
  a solution of the deautoconvolution problem
  $\mathcal A_k[u] = g^\dagger$, then $u$ is uniquely determined up to
  the global sign.
\end{Proposition}

\begin{Proof}
  In difference to \cite{GHB+14}, our proof is mainly based on
  \PFourier\ analysis.  Using the (\PFourier) convolution theorem, we
  notice that $u$ is a solution of the deautoconvolution problem
  $\mathcal A_k[u] = g^\dagger$ if and only if
  $\Fourier [u]^2 = \Fourier [g^\dagger]$, where $\Fourier$ denotes
  the Fourier transform given by
  \raisebox{0pt}[0pt][0pt]{$\Fourier[u](\omega) \coloneqq
    \int_{-\infty}^\infty u(t) \, \e^{- \I \omega t} \diff t$}
  for $\omega \in \R$. If $\Fourier [g^\dagger]$ is constantly zero,
  the assertion is trivially true.  Otherwise, the \PFourier\
  transforms of the compactly supported signals $u$ and $g^\dagger$
  are restrictions of entire functions by the theorem of
  \PPaley-\PWiener, and the related entire functions are completely
  determined by $\Fourier[u]$ and $\Fourier[g^\dagger]$.  Therefore,
  we may always find a point $\omega_0$ such that $\Fourier [u]$ and
  $\Fourier [g^\dagger]$ are non-zero in an appropriately small
  neighbourhood $U$ around $\omega_0$.  On this neighbourhood,
  there exist exactly two roots of $\Fourier [g^\dagger]$, which are
  restrictions of holomorphic functions, see for instance
  \cite[Section~5]{FL12}.  Consequently, the Fourier transform of $u$
  is either
  \begin{equation*}
    \Fourier [u] (\omega) = \sqrt{\Fourier [g^\dagger](\omega)}
    \qquad\text{or}\qquad
    \Fourier [u] (\omega) = - \sqrt{\Fourier [g^\dagger](\omega)}
  \end{equation*}
  for all $\omega \in U$.  Extending $\Fourier [u]$ from the
  neighbourhood to the whole real line by using the unique,
  corresponding entire function, one can conclude that $u$ and $-u$
  are the only possible solutions.  \qed
\end{Proof}

We proceed by considering the perturbed deautoconvolution problem
\begin{equation*}
  \mathcal A_k [u] = g^\delta
  \qquad\text{with}\qquad
  \pNormn{g^\delta - g^\dagger} \le \delta.
\end{equation*}
Since the autoconvolution is a bounded dilinear mapping, we can apply
the developed dilinear regularization theory to the \PTikhonov\
functional
\begin{equation}
  \label{eq:tik-autoconv}
  J_\alpha(u) \coloneqq \pNormn{\mathcal A_k [u] - g^\delta}^2 +
  \alpha \pNormn{u}^2,
  \addmathskip
\end{equation}
where we choose the squared norm of $L^2_\C([0,1])$ as regularization
term.  In order to verify that the autoconvolution fulfils the
required assumptions and to analyse the source-wise representation
\eqref{eq:source-cond}, firstly, we need a suitable representation of
the kernel-based autoconvolution and of the dual and predual spaces of
$L^2_\C([0,1]) \otimes_{\uppi, \mathrm{sym}} L^2_\C([0,1])$.

For this, we notice that each tensor $w$ in the tensor product
$H \otimes_\uppi H$, where $H$ is an arbitrary \PHilbert\ space, may
be interpreted as a nuclear operator in sense of the mapping $L_w$ or
$R_w$ in (\ref{eq:tensor-lin-op}--\ref{eq:left-right-map}).
\bq{Nuclear} here means that the singular values of $L_w$ and $R_w$
are absolutely summable.  More precisely, the projective tensor
product of a \PHilbert\ space with itself is isometrically isomorphic
to the space of nuclear operators from $H$ to $H$, \ie\
\begin{equation*}
  H \otimes_{\uppi} H \simeq \mathcal N(H),
  \addmathskip
\end{equation*}
see for instance \cite[Section~VI.5]{Wer02}.  Further, the injective
tensor product is here isometrically isomorphic to the space of
compact operators, \ie\
\begin{equation*}
  H \otimes_\upepsilon H \simeq \mathcal K(H),
\end{equation*}
see \cite[Corollary~4.13]{Rya02}, and thus becomes a true predual of
the projective tensor product, which means
\begin{equation*}
  (H \otimes_\epsilon H)^* \simeq H \otimes_\uppi H,
  \addmathskip
\end{equation*}
see \cite[Satz~VI.6.4]{Wer02}.  Finally, the dual space of the
projective tensor product coincides with the space of linear
operators, \ie\
\begin{equation*}
  (H \otimes_\uppi H)^* \simeq \mathcal L(H),
\end{equation*}
see \cite[Satz~VI.6.4]{Wer02}.  The action of an operator $T$ in
$\mathcal K(H)$ or $\mathcal L(H)$ on an arbitrary tensor $w$ in
$H \otimes_\uppi H$ is given by the tensorial lifting of the bilinear
mapping
\begin{equation*}
  (u,v) \mapsto \iProd{Tu}{v}
\end{equation*}
or by the \PHilbert-\PSchmidt\ inner product if $w$ is interpreted as
nuclear operator.  The corresponding spaces for the symmetric
projective tensor product $H \otimes_{\uppi, \mathrm{sym}} H$
coincide with the related self-adjoint operators.  In particular,
these identifications are valid for the real \PHilbert\ space
$L^2_\C([0,1])$.

To handle the autoconvolution $\mathcal A_k$, we split this mapping
into a linear integral part
\raisebox{0pt}[0pt][0pt]{$\mathcal I_k \colon L^2_{\C, \mathrm{sym}}
  ([0,1]^2) \rightarrow L^2_\C([0,2])$}
and into a quadratic part
\raisebox{0pt}[0pt][0pt]{$\odot \colon L^2_\C([0,1]) \rightarrow
  L^2_{\C, \mathrm{sym}} ([0,1]^2)$},
where $L^2_{\C, \mathrm{sym}} ([0,1]^2)$ denotes the subspace of the
symmetric square-in\-te\-grable functions on the square
$[0,1] \times [0,1]$ defined by
\begin{equation*}
  L^2_{\C,\mathrm{sym}}([0,1]^2) \coloneqq \bigl\{ u \in L^2_\C([0,1]^2) :
  u(s,t)=u(t,s) ~\text{for almost every}~ 0\le s,t \le 1 \bigr\}. 
\end{equation*}
More precisely, we employ the factorization
$\mathcal A_k = \mathcal I_k \circ {\odot}$ where the operators
$\mathcal I_k$ and $\odot$ are given by
\begin{equation}
  \label{eq:int-part}
  \mathcal I_k[w](t) \coloneqq \int_{-\infty}^\infty k(s,t) \,
  w(s,t-s) \diff s 
  \qquad (t \in \R)
\end{equation}
and
\begin{equation}
  \label{eq:quad-part}
  \odot [u](s,t) \coloneqq u(s) \, u(t)
  \qquad (s,t \in \R).
  \addmathskip
\end{equation}
On the basis of these mappings, we can determine the lifting of the
autoconvolution $\mathcal A_k$.

\begin{Lemma}[Lifting of the autoconvolution]
  \label{lem:lift-autoconv}
  The unique (quadratic) lifting of the kernel-based autoconvolution
  $\mathcal A_k = \mathcal I_k \circ {\odot}$ is given by
  \begin{equation*}
    \breve{\mathcal A}_k = \mathcal I_k \circ \breve \odot,
  \end{equation*}
  the unique dilinear lifting by
  \begin{equation*}
    (0, \breve{\mathcal A_k}) = (0, \mathcal I_k \circ \breve \odot).  
  \end{equation*}
\end{Lemma}

For the dilinear lifting, we here use a matrix-vector-like
representation.  More precisely, the mapping
$(0, \breve{\mathcal A}_k)$ is defined by
\begin{equation*}
  (0, \breve{\mathcal A}_k )
  \begin{pmatrix}
    u \\ w
  \end{pmatrix}
  \coloneqq 0[u] + \breve{\mathcal A}_k[w]
\end{equation*}
with $u \in L^2_\C([0,1])$ and
$w \in L^2_\C([0,1]) \otimes_{\uppi, \mathrm{sym}} L^2_\C([0,1])$.

\begin{Proof}[\thref{lem:lift-autoconv}]
  Firstly, we notice that the related bilinear mapping of the
  quadratic mapping $\odot$ is given by
  \begin{equation*}
    (u,v) \mapsto \bigl( (s,t) \mapsto u(s) \, v(t) \bigr)
    \addmathskip
  \end{equation*}
  for arbitrary $u$ and $v$ in $L^2_\C([0,1])$.  Obviously, this
  bilinear map is bounded with norm one, which means that there exists
  a unique bilinear lifting, see \thref{prop:lift-boun-bil-map}.
  Restricting this lifting to the symmetric projective tensor product,
  we obtain the unique quadratic lifting
  \raisebox{0pt}[0pt][0pt]{$\breve \odot \colon L^2_\C ([0,1])
    \otimes_{\uppi, \mathrm{sym}} L^2_\C([0,1]) \rightarrow
    L^2_\C([0,1]^2)$}
  with $\breve \odot[u \otimes u] = {\odot}[u]$.  Since the integral
  operator $\mathcal I_k$ is bounded with norm $\pNormn{k}_\infty$,
  the complete mapping
  $\breve{\mathcal A}_k = \mathcal I_k \circ {\breve \odot}$ is bounded
  too.  Moreover, we have
  \begin{equation*}
    \breve{\mathcal A}_k[u \otimes u] 
    = (\mathcal I_k \circ \breve \odot) [u \otimes u] 
    = (\mathcal I_k \circ {\odot})[u]
    = \mathcal A_k [u];
  \end{equation*}
  so the defined mapping $\breve{\mathcal A}_k$ is the unique
  quadratic lifting of the autoconvolution $\mathcal A_k$.

  In \thref{ex:quad-map}, we have already seen that the unique
  dilinear lifting with domain of definition
  $L^2_\C([0,1]) \times (L^2_\C([0,1]) \otimes_{\uppi, \mathrm{sym}}
  L^2_\C([0,1]))$
  of an arbitrary quadratic mapping is completely independent of the
  first component space.  In other words, the first component is always
  mapped to zero, which completes the proof.  \qed
\end{Proof}

\begin{Remark}
  \label{rem:lift-autoconv}
  Although the mapping $\breve \odot$ looks like a continuous
  embedding of the symmetric projective tensor product
  $L^2_\C([0,1]) \otimes_{\uppi, \mathrm{sym}} L^2_\C([0,1])$ into
  $L^2([0,1]^2)$, this is not the case.  Since we consider
  $L^2_\C([0,1])$ as an \R-linear space, the functions $u$ and $\I u$
  are orthogonal and thus linear independent.  The reason for this is
  simply that $\I$ is not a real scalar.  As a consequence, for
  $u \ne 0$, the tensors $u \otimes u$ and $\I u \otimes \I u$ are
  also linear independent, see \cite[Proposition~1.1]{Rya02}.  Now,
  the image of $u \otimes u + \I u \otimes \I u$ is the zero function,
  which shows that $\breve \odot$ is not injective.  If one restricts
  oneself to the function space of real functions $L^2_\R([0,1])$, or
  if one consider $L^2_\C([0,1])$ as $\C$-linear space, then
  $\breve \odot$ truly becomes a continuous embedding.  \qed
\end{Remark}

With the factorization of the kernel-based autoconvolution in mind, we
next determine the related adjoint operator of the lifting
$\breve{\mathcal A}_k$.

\begin{Lemma}[Adjoint of the autoconvolution]
  \label{lem:adj-autoconv}
  The adjoint operator
  $\breve{\mathcal A}_k^* \colon L^2_\C([0,2]) \rightarrow
  (L^2_\C([0,1]) \otimes_{\uppi,\mathrm{sym}} L^2_\C([0,1]))^*$
  of the quadratic lifting of the kernel-based autoconvolution
  $\mathcal A_k$ is given by
  $\breve{\mathcal A}_k^* = \breve \odot^* \circ \mathcal I_k^*$ with
  \begin{equation*}
    \mathcal I_k^* [\phi](s,t) = \overline{k(s,s+t)} \, \phi(s+t)
    \qquad (s,t \in [0,1]),
  \end{equation*}
  where $\mathcal I_k^*$ is a mapping between $L^2_\C([0,2])$ and
  $L^2_{\C, \mathrm{sym}}([0,1]^2)$, and the adjoint operator of the
  dilinear lifting by
  \begin{equation*}
    (0, \breve{\mathcal A}_k^*)^\T 
    = (0, \breve \odot^* \circ \mathcal I_k^*)^\T.
  \end{equation*}
\end{Lemma}

For the adjoint dilinear lifting, we again use a matrix-vector-like
representation.  In more detail, the mapping $(0, \breve{\mathcal
  A}_k^*)^\T$ is defined by
\begin{equation*}
  \begin{pmatrix}
    0 \\ \breve{\mathcal A}_k^*
  \end{pmatrix}
  \phi
  \coloneqq
  \begin{pmatrix}
    0[\phi] \\ \breve{\mathcal A}_k^*[\phi]
  \end{pmatrix}
  \submathskip
\end{equation*}
with $\phi \in L^2_\C([0,2])$.

\begin{Proof}[\thref{lem:adj-autoconv}]
  The representation of the adjoint integral operator $\mathcal I_k^*$
  can be directly verified by
  \begin{align*}
    \iProd{\phi}{\mathcal I_k[w]}_\R
    &= \Re \biggl[\int_{-\infty}^\infty \int_{-\infty}^\infty
      \overline{w(s,t-s) \, k(s,t)} \, \phi(t) \diff s \diff t \biggr]
    \\[\fskip] 
    &= \Re \biggl[\int_{-\infty}^\infty \int_{-\infty}^\infty
      \overline{w(s,t) \, k(s,t+s)} \, \phi(t+s) \diff s \diff t \biggr] 
      =  \iProd{\mathcal I_k^*[\phi]}{w}_\R
  \end{align*}
  for all $w \in L^2_\C([0,1]^2)$ and $\phi \in L^2_\C([0,2])$.  The
  remaining assertion immediately follows from the properties of the
  adjoint.  \qed
\end{Proof}

\begin{Remark}
  \label{rem:adj-autoconv:1}
  The lifted operator $\breve \odot$ maps every tensor in the
  symmetric projective tensor product to a square-integrable function.
  The adjoint mapping $\breve \odot^*$ thus allows us to interpret
  each function $\omega$ in $L^2_{\C,\mathrm{sym}}([0,1]^2)$ as a
  bounded linear functional on the projective tensor product
  \raisebox{0pt}[0pt][0pt]{$L^2_\C([0,1]) \otimes_{\uppi,
      \mathrm{sym}} L^2_\C([0,1])$}.
  More precisely, the action of the function $\omega$ on an arbitrary
  tensor $w$ is given by
  \begin{equation}
    \label{eq:act-adj-odot}
    w \mapsto \iProd{\omega}{\breve \odot [w]}_\R 
    = \Re \iProd{\omega}{\breve \odot[w]}_\C,
  \end{equation}
  which we will exploit later.  \qed
\end{Remark}

\begin{Remark}
  \label{rem:adj-autoconv:2}
  For the trivial kernel $k \equiv 1$, the image of the adjoint
  operator $\mathcal I_k^*$ consists of all functions $w$ in
  $L^2_{\C,\mathrm{sym}}([0,1]^2)$ with \PHankel\ structure, where the
  diagonal $t \mapsto w(t,t)$ coincides with the given $\phi$ up to
  scaling.  \qed
\end{Remark}

The next ingredients for the dilinear regularization theory in
\autoref{sec:dil-inverse-probl} are the required conditions in
\thref{ass:lift-prop}.  As a preparing step, we prove that the lifted
quadratic operator $\breve \odot$ is sequentially weakly$*$
continuous.

\begin{Lemma}[Sequentially weak$*$ continuity]
  \label{lem:weak-cont-quad-part}
  The quadratic operator $\odot$ defined in \eqref{eq:quad-part} is
  sequentially weakly$*$ continuous.
\end{Lemma}

\begin{Proof}
  The central idea of the proof is to exploit
  \thref{prop:seq-weak-cont-dil-oper}, which means that we have to
  show the sequential weak$*$ continuity of the dilinear lifting
  $\breve \odot$ with respect to the topology induced by
  $L^2_\C([0,1]) \otimes_{\epsilonup, \mathrm{sym}} L^2_\C([0,1])$.
  Since the dilinear lifting of a quadratic operator is independent of
  the first component space, \cf\ \thref{ex:lin-amp}, it is enough to
  show the sequential weak$*$ continuity of the quadratic lifting
  $\breve \odot$ in \thref{lem:lift-autoconv}.  For this, we have to
  show
  \begin{equation}
    \label{eq:omega-test}
    \iProd{\omega}{\breve \odot[w_n]}_\R \rightarrow
    \iProd{\omega}{\breve\odot[w]}_\R
    \addmathskip
  \end{equation}
  for every $\omega$ in $L^2_{\C, \mathrm{sym}}([0,1]^2)$, and for
  every weakly$*$ convergent sequence $w_n \weakrightharpoonup w$.

  As mentioned above, the symmetric injective tensor product
  $L^2_\C([0,1]) \otimes_{\epsilonup} L^2_\C([0,1])$ is here
  isometrically isomorphic to the self-adjoint compact operators,
  where the action of a self-adjoint compact operator $\Phi$ is given
  by the lifting of the quadratic form
  $(u,u) \mapsto \iProd{\Phi[u]}{u}_\R$.  This observation is the key
  component to establish the assertion.  More precisely, if we can
  show that the action of an arbitrary symmetric function $\omega$ see
  \thref{rem:adj-autoconv:1}, which is equivalent to the testing in
  \eqref{eq:omega-test}, corresponds to the lifting of a self-adjoint
  compact operator, then the assertion is trivially true.

  For this purpose, we consider the action of $\omega$ to a symmetric
  rank-one tensor $u \otimes u$, which may be
  written as
  \begin{equation}
    \label{eq:act-adj-quad-part:rank-one}
    \iProd{\omega}{\breve \odot[u \otimes u]}_\R
    = \Re \biggl[ \int_0^1 \!\! \int_0^1 \omega(s,t) \,
    \overline{u(s)} \; \overline{u(t)} \diff 
    s \diff t \biggr]
    = \iProd{\Phi_\omega[u]}{u}_\R,
  \end{equation}
  where the operator $\Phi_\omega \colon L^2_\C([0,1]) \rightarrow
  L^2_\C([0,1])$ is defined by
  \begin{equation}
    \label{eq:act-adj-quad-part:int-op}
    \Phi_\omega[u](t) \coloneqq \int_0^1 \omega (s,t) \, \overline{u(s)} \diff s
    \qquad (t \in [0,1]).
  \end{equation}
  The central observation is that the $\R$-linear operator
  $\Phi_\omega$ resembles a \PFredholm\ integral operator.  Due to the
  occurring conjugation, for a fixed point $(s,t)$, the multiplication with
  $w(s,t)$ here only acts as an $\R$-linear mapping  instead of an
  usually $\C$-linear mapping on $u(s)$.  Similarly to the
  classical theory, see for instance \cite{Wer02}, we can approximate
  the kernel function $\omega$ by a sequence of appropriate step
  functions $\omega_n$ on a rectangular partition of $[0,1]^2$ such
  that $\omega_n \rightarrow \omega$ in
  $L^2_{\C,\mathrm{sym}}([0,1]^2)$.  The ranges of the related
  operators $\Phi_{\omega_n}$ are finite-dimensional, and because of
  the convergence of $\Phi_{\omega_n}$ to $\Phi_\omega$, the
  compactness of $\Phi_\omega$ follows.  Since $\omega$ is symmetric,
  the operator $\Phi_\omega$ is moreover self-adjoint.  Consequently,
  the quadratic lifting
  $w \mapsto \iProdn{\omega}{\breve \odot[w]}_\R$ with
  $w \in L^2_\C([0,1]) \otimes_{\uppi, \mathrm{sym}} L^2_\C ([0,1])$
  of the quadratic form $u \mapsto \iProdn{\Phi_\omega[u]}{u}_\R$ is
  contained in the symmetric injective tensor product
  $L^2_\C([0,1]) \otimes_{\upepsilon, \mathrm{sym}} L^2_\C([0,1])$.

  For a weakly$*$ convergent sequence $w_n \weakrightharpoonup w$, we
  thus have \eqref{eq:omega-test} for all $\omega$ in
  $L^2_\C([0,1]^2)$, which shows the sequential weak$*$ continuity of
  the quadratic lifting $\breve \odot$ and hence of the dilinear
  lifting $(0, \breve \odot)$.  The sequential weak$*$ continuity of
  the dilinear operator $\odot$ now immediately follows from
  \thref{prop:seq-weak-cont-dil-oper}.  \qed
\end{Proof}

Using \thref{lem:weak-cont-quad-part}, we can now verify that the
\PTikhonov\ regularization of the deautoconvolution problem satisfies
the required assumptions for the developed regularization theory in
\autoref{sec:dil-inverse-probl}.

\begin{Lemma}[Verification of required assumptions]
  \label{lem:check-ass-autoconv}
  The \PTikhonov\ functional $J_\alpha$ in \eqref{eq:tik-autoconv}
  related to the kernel-based autoconvolution $\mathcal A_k$ fulfils
  the requirements in \thref{ass:lift-prop}.
\end{Lemma}

\begin{Proof}
  We briefly verify the needed requirements step by step.
  \begin{enumerate}[(i)]
  \item Obviously, the \PTikhonov\ functional $J_\alpha$ in
    \eqref{eq:tik-autoconv} is coercive since the regularization term
    coincides with the squared \PHilbert\ space norm.
  \item For the same reason, the regularization term is sequentially
    weakly$*$ lower semi-continuous, see for instance
    \cite[Theorem~2.6.14]{Meg98}.
  \item In \thref{lem:weak-cont-quad-part}, we have already proven
    that the quadratic operator $\odot$ of the factorization
    $\mathcal A_k = \mathcal I_k \circ {\odot}$ is sequentially
    weakly$*$ continuous.  Further, the obvious norm-to-norm
    continuity of the integral operator $\mathcal I_k$ implies the
    weak-to-weak continuity of $\mathcal I_k$, see for instance
    \cite[Proposition~2.5.10]{Meg98}.  Since $L^2_\C([0,1]^2)$ and
    $L^2_\C([0,2])$ are \PHilbert\ spaces, the mapping $\mathcal I_k$
    is weakly$*$-to-weakly$*$ continuous as well.  Consequently, the
    composition $\mathcal A_k = \mathcal I_k \circ {\odot}$ is
    sequentially weakly$*$ continuous as required.  \qed
  \end{enumerate}
\end{Proof}

\begin{Remark}
  \label{rem:check-ass-autoconv}
  Since the \PTikhonov\ functional $J_\alpha$ for the
  deautoconvolution problem fulfils all constraints in
  \thref{ass:lift-prop} by \thref{lem:check-ass-autoconv}, we can
  employ the developed dilinear/diconvex regularization theory in
  \autoref{sec:dil-inverse-probl}.  Therefore, the related
  minimization problem
  \begin{equation*}
    \minimize_{u \in L^2_\C([0,1])} \; \pNormn{\mathcal A_k [u] -
      g^\delta}^2 + \alpha \pNormn{u}^2
    \addmathskip
  \end{equation*}
  is well posed, stable, and consistent.  \qed
\end{Remark}

Besides the well-known well-posedness, stability, and consistency of
the regularized deautoconvolution problem, the convergence rate
introduced by \thref{the:conv-rate} is far more interesting.  In order
to study the employed source condition, we have to compare the range
of the adjoint autoconvolution operator
$(0,\breve{\mathcal A}_k^*)^\T$, see \thref{lem:adj-autoconv}, and the
dilinear subdifferential of the squared \PHilbert\ norm on
$L^2_\C([0,1])$.

More generally, we initially determine the dilinear subdifferential
for the norm on an arbitrary \PHilbert\ space $H$.  For this purpose,
we exploit that the dual space $(H \otimes_{\uppi, \mathrm{sym}} H)^*$
is isometrically isomorphic to the space of self-adjoint, bounded
linear operators.  As mentioned above, the action of a specific
self-adjoint, bounded linear operator $\Phi$ on an arbitrary symmetric
tensor $w$ in $H \otimes_{\uppi, \mathrm{sym}} H$ is given by the
lifting of the quadratic map $u \mapsto \iProd{\Phi[u]}{u}_H$ with
$u \in H$.  In the following, we use the dual pairing notation
$\iProd{\Phi}{w}_{H \otimes_{\uppi, \mathrm{sym}} H}$ to refer to this
action.  With this preliminary considerations, the dilinear
subdifferential of an arbitrary \PHilbert\ norm is given in the
following manner, where the operator $\Id$ denotes the identity and
$\mathcal S_-(H)$ the set of all self-adjoint and negative
semi-definite operators on $H$.

\begin{Theorem}[Dilinear Subdifferential of \PHilbert\ norms]
  \label{the:dil-sub-hilb-norm}
  Let $H$ be a real \PHilbert\ space with inner product
  $\iProdn{\cdot}{\cdot}_H$ and norm $\pNorm{\cdot}_H$.  The dilinear
  subdifferential of the squared \PHilbert\ norm is given by
  \begin{equation*}
    \uppartial_\upbeta \pNorm{\cdot}^2_H (u) = \bigl\{ (-2 \, T u, \Id
    + T)  : T \in \mathcal S_-(H) \bigr\},
  \end{equation*}
  and of the $p$th power with $p \in (1,2)$ by
  \begin{equation*}
    \uppartial_\upbeta \pNorm{\cdot}^p_H(u) 
    = \bigl\{ \bigl((p \pNorm{u}^{p-2}_H \Id +
    2T)\, u, T \bigr) : T \in \mathcal S_-(H) \bigr\}.
  \end{equation*}
\end{Theorem}

\begin{Proof}
  Firstly, we notice that the $p$th power of the \PHilbert\ norm is
  differentiable with \PGateaux\ derivative
  \begin{equation*}
    D \pNorm{\cdot}^p_H (u) = p \pNorm{u}^{p-2}_H u
    \qquad(u \in H, p>1), 
  \end{equation*}
  see for instance \cite{BL11}.  Further, the squared \PHilbert\ norm
  is itself a dilinear mapping with representative
  $(u,w) \mapsto \iProdn{\Id}{w}_{H \otimes_{\uppi, \mathrm{sym}} H}$
  and
  $\pNormn{u}^2 = \iProdn{\Id}{u \otimes u}_{H \otimes_{\uppi,
      \mathrm{sym}} H}$,
  and the \PTaylor\ expansion at some point $u$ in $H$ is given by
  \begin{equation*}
    \pNorm{v}^2_H = \pNorm{u}^2_H + \iProd{2u}{v-u}_H +
    \iProd{\Id}{(v-u)\otimes(v-u)}_{H \otimes_{\uppi, \mathrm{sym}} H}.  
    \addmathskip
  \end{equation*}
  Since all dilinear mappings are differentiable in the \PHilbert\
  space setting, the linear part of the \PTaylor\ expansion is fixed,
  and the only possibility to construct a dilinear mapping beneath the
  squared \PHilbert\ norm is to add a negative semi-definite quadratic
  part $T$.  In this manner, we obtain
  \begin{align*}
    \pNorm{v}^2_H
    & \ge \pNorm{u}^2_H + \iProd{2u}{v-u}_H + \iProd{\Id +T}{(v-u)
      \otimes (v-u)}_{H \otimes_{\uppi, \mathrm{sym}} H} 
    \\[\fskip]
    & = \pNorm{u}^2_H + \iProd{-2 \, T u}{v-u}_H + \iProd{\Id + T}{v
      \otimes v - u \otimes u}_{H \otimes_{\uppi, \mathrm{sym}} H},
  \end{align*}
  which yields the claimed dilinear subgradients.

  Next, we consider the case $p\in(1,2)$.  Again, at an arbitrary
  point $u$ in $H$, the derivative of all dilinear mappings beneath
  $\pNorm{\cdot}^p_H$ in $u$ have to be $p \pNormn{u}^{p-2}_H \, u$.  In
  other words, in the \PTaylor\ expansion form, the dilinear
  subgradients have to satisfy
  \begin{equation*}
    \pNorm{v}^p_H \ge \pNorm{u}^p_H + \iProdb{p \pNorm{u}^{p-2}_H u}{v-u}_H
    + \iProdb{T}{(v-u) \otimes (v-u)}_{H \otimes_{\uppi, \mathrm{sym}} H}
    \qquad(v \in H)
  \end{equation*}
  for some suitable self-adjoint operator $T \colon H \rightarrow H$.
  Replacing $v$ by $t v$ with $t \in \R$, we obtain the condition
  \begin{equation*}
    \begin{aligned}
      t^p \pNorm{v}^p_H &\ge \bigl[ \pNorm{u}^p_H - \iProdb{p
        \pNorm{u}^{p-2}_H u}{u}_H + \iProdb{T}{u \otimes u}_{H
        \otimes_{\uppi, \mathrm{sym}} H} \bigr] 
      \\[\fsmallskip]
      &\quad + t \bigl[ \iProdb{p \pNorm{u}^{p-2}_H u}{v}_H - 2
      \iProdb{T}{v \otimes u}_{H \otimes_{\uppi, \mathrm{sym}} H}
      \bigr]
      \\[\fsmallskip]
      &\quad+ t^2 \bigl[ \iProdb{T}{ v \otimes v}_{H \otimes_{\uppi,
          \mathrm{sym}} H} \bigr] 
    \end{aligned}
    \qquad(v \in H, t \in \R).
  \end{equation*}
  Since the left-hand side is growing with a power of $p<2$, the
  coefficient
  $\iProdn{T}{v \otimes v}_{H \otimes_{\uppi, \mathrm{sym}} H} =
  \iProd{T v}{v}_H$
  on the right-hand side have to be non-positive for all $v$ in $H$.
  The self-adjoint operator $T$ has thus to be negative semi-definite.
  Rewriting the subgradient condition in the usual form, we have
  \begin{equation*}
    \pNorm{v}^p_H \ge \pNorm{u}^p_H + \iProdb{(p \pNorm{u}^{p-2}_H \Id - 2T)
      \, u}{v-u}_H + \iProdb{T}{v \otimes v - u \otimes u}_{H
      \otimes_{\uppi, \mathrm{sym}} H} 
  \end{equation*}
  and hence the dilinear subdifferential in the assertion.  \qed
\end{Proof}

Similarly to the dilinear subdifferential of the squared \PHilbert\
norm, we use the identification of the dual space
$(L^2_\C([0,1]) \otimes_{\uppi, \mathrm{sym}} L^2_\C([0,1]))^*$ with
the space of self-adjoint, bounded linear operators to describe the
range of the adjoint $(0,\breve{\mathcal A}_k^*)^\T$ of the dilinearly
lifted autoconvolution.  More precisely, using the unique identification
(\ref{eq:act-adj-quad-part:rank-one}--\ref{eq:act-adj-quad-part:int-op}),
and incorporating the factorization
$\breve{\mathcal A}_k^* = \breve \odot^* \circ \mathcal I_k^*$, we may
write the range of the adjoint lifted operator as
\begin{equation}
  \label{eq:ran-adj-autoconv}
  \ran (0,\breve{\mathcal A}_k^*)^\T = \bigl\{ (0, \Phi_{\mathcal
    I_k^*[\phi]}) :  \phi \in L^2_\C([0,2]) \bigr\},
  \addmathskip
\end{equation}
where the self-adjoint, bounded integral operator $\Phi_{\mathcal
  I_k^*[\phi]}$ is given by
\begin{equation}
  \label{eq:ran-adj-autoconv:Phi}
  \Phi_{\mathcal I_k^*[\phi]} [u] (t)
  \coloneqq \int_0^1  \overline{k(s,s+t)} \, \phi(s+t) \,
  \overline{u(s)} \diff s 
  \qquad(t \in [0,1]).
\end{equation}
In our specific setting, the operator $\Phi_{\mathcal I_k^*[\phi]}$ is
additionally compact since the range of $\breve{\odot}^*$ is contained
in the injective tensor product as discussed in the proof of
\thref{lem:weak-cont-quad-part}.

Comparing the range of $(0,\breve{\mathcal A}_k^*)^\T$ and the
subdifferential of the squared \PHilbert\ norm in
\thref{the:dil-sub-hilb-norm}, we notice that an element
$(-2 \, Tu^\dagger, \Id + T)$ of
$\uppartial_\beta \pNormn{\cdot}^2(u^\dagger)$ is contained in the
range of the adjoint $(0, \breve{\mathcal A}_k^*)^\T$ if and only if
\begin{equation*}
  Tu^\dagger = 0
  \qquad \text{and} \qquad
  \Id + T = \Phi_{\mathcal I_k^*[\phi]}
\end{equation*}
for some $\phi$ in $L^2_\C([0,2])$.  Since $T$ is a self-adjoint,
negative semi-definite operator, the spectrum
$\sigma(\Phi_{\mathcal I_k^*[\phi]}) = \sigma(\Id + T)$ of eigenvalues
of $\Phi_{\mathcal I_k^*[\phi]}$ is bounded from above by one.
Considering that 
\begin{equation}
  \label{eq:sym-spec-Phi}
  \Phi_{\mathcal I_k^*[\phi]} [v] = \lambda v
  \qquad \text{implies} \qquad
  \Phi_{\mathcal I_k^*[\phi]} [\I v] = - \lambda \, \I v
\end{equation}
for every eigenfunction $v$ related to the eigenvalue $\lambda$, we
see that the spectrum of $\Phi_{\mathcal I_k^*[\phi]}$ is moreover
symmetric, which means
\begin{equation*}
  \sigma( \Phi_{\mathcal I_k^*[\phi]} ) \subset [-1,1].
  \submathskip
\end{equation*}
Further, the equation
\begin{equation*}
  T u^\dagger = (\Phi_{\mathcal I_k^*[\phi]} - \Id) \, u^\dagger = 0
  \addmathskip
\end{equation*}
yields that $u^\dagger$ has to be an eigenfunction of
$\Phi_{\mathcal I_k^*[\phi]}$ with respect to the eigenvalue one.  In
summary, the source-wise representation \eqref{eq:source-cond} may
thus be rewritten in the following form.

\begin{Theorem}[Source condition -- deautoconvolution]
  \label{the:sour-cond-deauto}
  For a norm-min\-i\-miz\-ing solution $u^\dagger$ of the kernel-based
  deautoconvolution problem $\mathcal A_k[u] = g^\dagger$, the source
  condition \eqref{eq:source-cond} is fulfilled if and only if there
  exists a $\phi$ in $L^2_\C([0,2])$ such that
  \begin{equation}
    \label{eq:source-cond-deauto:norm-eigen}
    \pNormn{\Phi_{\mathcal I_k^*[\phi]}} = 1
    \qquad\text{and}\qquad
    \Phi_{\mathcal I_k^*[\phi]}[u^\dagger] = u^\dagger
    \addmathskip
  \end{equation}
  for the integral operator
  \submathskip
  \begin{equation}
    \label{eq:source-cond-deauto:oper}
    \Phi_{\mathcal I_k^*[\phi]}[u](t)
    \coloneqq \int_{0}^1 \overline{k(s,s+t)} \, \phi(s+t) \,
    \overline{u(s)} \diff s
    \qquad(t \in [0,1]).
  \end{equation}
\end{Theorem}

Starting from an arbitrary integral operator
$\Phi_{\mathcal I_k^*[\phi]}$ or, more precisely, from a function
$\phi \in L^2_\C([0,1])$, we can easily construct solutions to the
deautoconvolution problem that fulfil the source condition
\eqref{eq:source-cond} or, equivalently,
(\ref{eq:source-cond-deauto:norm-eigen}--\ref{eq:source-cond-deauto:oper})
by rescaling the spectrum $\sigma(\Phi_{\mathcal I_k^*[\phi]})$
appropriately and by choosing $u^\dagger$ from the eigenspace of the
eigenvalue one.  For the trivial kernel $k \equiv 1$, the
eigenfunction condition simply means that
\begin{equation}
  \label{eq:for-op-conv-rep}
  \bigl(\phi * \overline{u^\dagger (-\cdot)} \,\bigr)
  \big\vert_{[0,1]} 
  = u^\dagger.
\end{equation}
Hence, the source-wise representation \eqref{eq:source-cond} can only
be fulfilled for continuous functions.  More generally, if the
integral operator $\Phi_{\mathcal I_k^*[\phi]}$ possesses smoothing
properties, the norm-minimizing solution $u^\dagger$ has to be
arbitrary smooth.  Our observations also cover the case of incomplete
measurements, which can be modelled by setting the kernel $k(s,t)$
zero for the corresponding $t \in [0,1]$.

Applying \thref{the:conv-rate} to the deautoconvolution problem, we
obtain the following convergence rate.

\begin{Corollary}[Convergence rate -- deautoconvolution]
  \label{cor:conv-rate-deauto}
  Let $u^\dagger$ be a norm-min\-i\-miz\-ing solution of the
  kernel-based deautoconvolution problem
  $\mathcal A_k[u] = g^\dagger$. If there exist an integral operator
  $\Phi_{\mathcal I_k^*[\phi]}$ fulfilling the conditions
  \eqref{eq:act-adj-quad-part:rank-one}, the minimizer
  $u_\alpha^\delta$ of the \PTikhonov\ regularization $J_\alpha$ in
  \eqref{eq:tik-autoconv} converges to
  \raisebox{0pt}[0pt][0pt]{$u^\dagger$} in the sense that the dilinear
  \PBregman\ distance between $u_\alpha^\delta$ and $u^\dagger$ with
  respect to the norm of $L^2_\C([0,1])$ is bounded by
  \begin{equation*}
    \Delta_{\upbeta, (0,\Phi_{\mathcal I_k^*[\phi]})}(u_\alpha^\delta, u^\dagger)
    \le \bigl( \tfrac{\delta}{\sqrt \alpha} + \tfrac{\sqrt \alpha}{2}
    \, \pNormn{\phi} \bigr)^2
  \end{equation*}
  and the data fidelity term by
  \begin{equation*}
    \pNormn{\mathcal A_k [u_\alpha^\delta] - g^\delta} \le \delta +
    \alpha \, \pNormn{\phi}.
    \addmathskip
  \end{equation*}
\end{Corollary}

Although \thref{cor:conv-rate-deauto} gives a convergence rate for the
deautoconvolution problem, the employed \PBregman\ distance is usually
very weak.  Knowing that, besides the true solution $u^\dagger$, the
function $-u^\dagger$ also solves the deautoconvolution problem, we
notice that the \PBregman\ distance
\raisebox{0pt}[0pt][0pt]{$\Delta_{\upbeta,(0,\Phi_{\mathcal
      I_k^*[\phi]})}$}
cannot measure distances in direction of $u^\dagger$.  Since the
squared \PHilbert\ norm is itself a dilinear functional, in the worst
case, it can happen that the \PBregman\ distance is constantly zero,
and that the derived convergence rate is completely useless.
Depending on the integral operator $\Phi_{\mathcal I_k^*[\phi]}$ or,
more precisely, on the eigenspace $E_1$ of the eigenvalue one, we can
estimate the \PBregman\ distance from below in the following manner.

\begin{Theorem}[\PBregman\ distance -- deautoconvolution]
  \label{the:brag-dist-deauto}
  Let\;$\Phi_{\mathcal I_k^*[\phi]}$ be an integral operator of the
  form
  {\upshape(}\ref{eq:source-cond-deauto:norm-eigen}--\ref{eq:source-cond-deauto:oper}{\upshape)}
  with distinct positive eigenvalues
  $1 = \lambda_1 > \lambda_2 > \dots$ and related finite-dimensional
  eigenspaces $E_1, E_2, \dots$ The corresponding dilinear \PBregman\
  distance $\Delta_{\upbeta, (0, \Phi_{\mathcal I_k^*[\phi]})}$ is
  bounded from below by
  \begin{equation*}
    \Delta_{\upbeta,(0,\Phi_{\mathcal I_k^*[\phi]})} (v,u^\dagger) \ge (1-\lambda_2)
    \pNorm{P_{E_1^\perp}(v - u^\dagger)}^2 \!,
  \end{equation*}
  where $P_{E_1^\perp}$ denotes the orthogonal projection onto the
  orthogonal complement of the eigen\-space $E_1$.
\end{Theorem}

\begin{Proof}
  The self-adjoint, compact integral operator
  $\Phi_{\mathcal I_k^*[\phi]} \colon L^2_\C([0,1]) \rightarrow
  L^2_\C([0,1])$
  possesses a symmetric spectrum, where the eigenspace of the
  eigenvalue $-\lambda_n$ is given by $\I E_n$, see
  \eqref{eq:sym-spec-Phi}; hence, the spectral theorem implies that
  the action of $\Phi_{\mathcal I_k^*[\phi]}$ is given by
  \begin{equation*}
    \Phi_{\mathcal I_k^*[\phi]}[v] = \sum_{n=1}^\infty \lambda_n \,
    P_{E_n} (v) - \lambda_n \, P_{\I E_n}(v)
    \qquad(v \in L^2_\C([0,1])),
  \end{equation*}
  where $P_{E_n}$ and $P_{\I E_n}$ denote the projections onto the
  eigenspace $E_n$ and $\I E_n$ respectively.  Considering that
  $u^\dagger$ is an eigenfunction in $E_1$, we may write the
  \PBregman\ distance for the squared \PHilbert\ norm as
  \begin{align*}
    \Delta_{\upbeta,(0,\Phi_{\mathcal I_k^*[\phi]})}(v, u^\dagger) 
    &= \pNormn{v}^2 - \pNormn{u^\dagger}^2 - \iProd{\Phi_{\mathcal I_k^*[\phi]}}{v
      \otimes v - u^\dagger \otimes u^\dagger}
    \\[\fskip]
    &= \pNormn{v}^2 - \iProd{\Phi_{\mathcal I_k^*[\phi]}[v]}{v}.
  \end{align*}
  Using the spectral representation of $\Phi_{\mathcal I_k^*[\phi]}$,
  and denoting the kernel of $\Phi_{\mathcal I_k^*[\phi]}$ by $E_0$,
  we have the estimation
  \begin{align*}
    \Delta_{\upbeta,(0,\Phi_{\mathcal I_k^*[\phi]})} (v, u^\dagger)
    &= \pNormn{P_{E_0}(v)}^2
      + \sum_{n=1}^\infty (1 - \lambda_n) \, \pNormn{P_{E_n}(v)}^2
      + (1 + \lambda_n) \, \pNormn{P_{\I E_n}(v)}^2
    \\[\fskip]
    &\ge (1-\lambda_2) \, \biggl[ \pNormn{P_{E_0}(v)}^2 
      + \sum_{n=2}^\infty \pNormn{P_{E_n}(v)}^2
      + \sum_{n=1}^\infty \pNormn{P_{\I E_n} (v)}^2
      \biggr]
    \\[\fskip]
    &= (1-\lambda_2) \, \pNormn{P_{E_1^\perp}(v - u^\dagger)}^2,
  \end{align*}
  which yields the assertion.  \qed
\end{Proof}

\begin{Remark}
  \label{rem:brag-dist-deauto}
  The quality of the convergence rate in \thref{cor:conv-rate-deauto}
  thus crucially depends on the integral operator
  $\Phi_{\mathcal I_k^*[\phi]}$ of the source condition.  If the true
  solution $u^\dagger$ is the only eigenfunction to the eigenvalue
  one, the minimizers \raisebox{0pt}[0pt][0pt]{$u_\alpha^\delta$} of
  the regularized perturbed problem converges nearly strongly to the
  true solution $u^\dagger$ or $-u^\dagger$ or, more precisely,
  converges to the \R-linear subspace spanned by $u^\dagger$ in norm
  topology.  If we choose $\alpha \sim \delta$ and assume
  $\dim E_1 = 1$, then $u_\alpha^\delta$ converges in
  $L^2_\C([0,1]) \mathbin/ \Span\{u^\dagger\}$ with a rate of
  $\Landau(\delta^{\nicefrac{1}{2}})$ strongly.  \qed
\end{Remark}

\section{Numerical simulations}
\label{sec:numer-simul}

Besides the theoretical regularization analysis in the previous
sections, we now perform numerical experiments to verify the
established convergence rates and to examine the error that is not
covered by the \PBregman\ distance.  As before, we restrict ourselves
to the deautoconvolution problem and, moreover, to the \bq{kernelless}
setup with the trivial kernel  $k \equiv 1$.

\subsection{Construction of valid source elements}
\label{sec:const-source-elem}

In a first step, we numerically construct signals that satisfy the
source condition in \thref{the:sour-cond-deauto}.  For this, we
approximate the integral operator
$\Phi_{\mathcal I_k^*[\phi]} \colon L_\C^2([0,1]) \rightarrow
L_\C^2([0,1])$
by applying the midpoint rule.  More precisely, we assume
\begin{equation}
  \label{eq:approx-int-oper-phi}
  \Phi_{\mathcal I_k^*[\phi]}[u]\bigl(\tfrac{n}{N}\bigr) \approx \frac{1}{N} \,
  \sum_{k=0}^N \phi\bigl( \tfrac{n+k}{N} \bigr) \, \overline{u\bigl(
    \tfrac{k}{N} \bigr)}
  \qquad (n=0, \dots, N)
\end{equation}
for an appropriate large positive integer $N$.  Starting from the
source element $\phi$, we now determine the eigenfunction $u_1$ to the
major eigenvalue $\lambda_1$.  Exploiting the convolution
representation of $\Phi_{\mathcal I_k^*[\phi]}$ in
\eqref{eq:for-op-conv-rep}, we may compute the action of the integral
operator $\Phi_{\mathcal I_k^*[\phi]}$ efficiently by applying the
\PFourier\ transform.  The eigenfunction $u_1$ itself can be
determined approximately by using the power iteration.  To overcome
the issue that the spectrum of $\Phi_{\mathcal I_k^*[\phi]}$ is
symmetric, which means that $-\lambda_1$ is also an eigenvalue, see
\eqref{eq:sym-spec-Phi}, we apply the power iteration to the operator
$\Phi_{\mathcal I_k^*[\phi]} \circ \Phi_{\mathcal I_k^*[\phi]}$.  In
so doing, we obtain a vector $v_1$ in the span of $E_1 \cup \I E_1$,
where $E_1$ is the eigenspace with respect to the eigenvalue
$\lambda_1$, and $\I E_1$ with respect to $-\lambda_1$, \cf\
\eqref{eq:sym-spec-Phi}.  The projections to $E_1$ and $\I E_1$ can,
however, simply be computed by
\begin{equation*}
  P_{E_1}(v_1) = \tfrac{1}{2} \, 
  \bigl(v_1 + \Phi_{\mathcal I_k^*[\phi]}[v_1]\bigr)
  \qquad\text{and}\qquad
  P_{\I E_1}(v_1) = \tfrac{1}{2} \, 
  \bigl(v_1 - \Phi_{\mathcal I_k^*[\phi]}[v_1]\bigr).
\end{equation*}

\begin{figure}\centering
  \subfloat[Magnitude of the source element $\phi^{(j)}$.]{
    \includegraphics{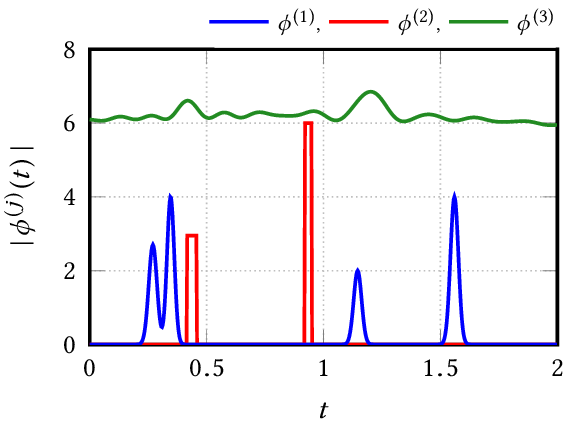}
  }
  \quad
  \subfloat[Phase of the source element $\phi^{(j)}$.]{
    \includegraphics{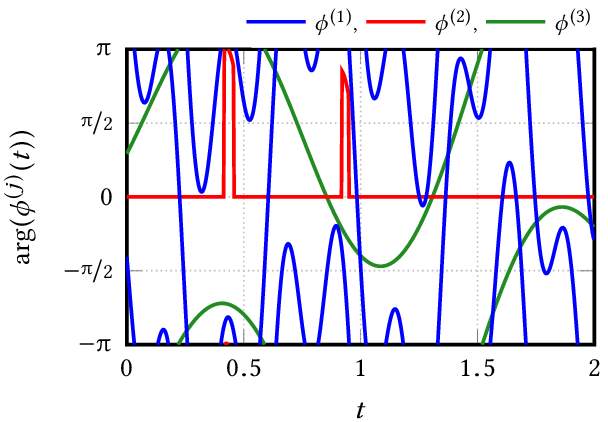}
  }
  \par
  \subfloat[Magnitude of the eigenfunction $u_1^{(j)}$.]{
    \includegraphics{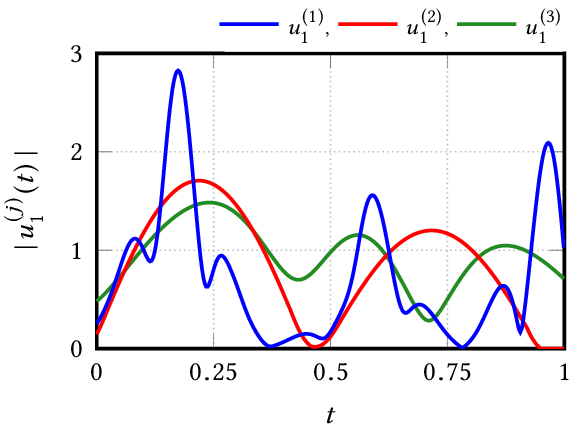}
  }
  \quad
  \subfloat[Phase of the eigenfunction $u_1^{(j)}$.]{
    \includegraphics{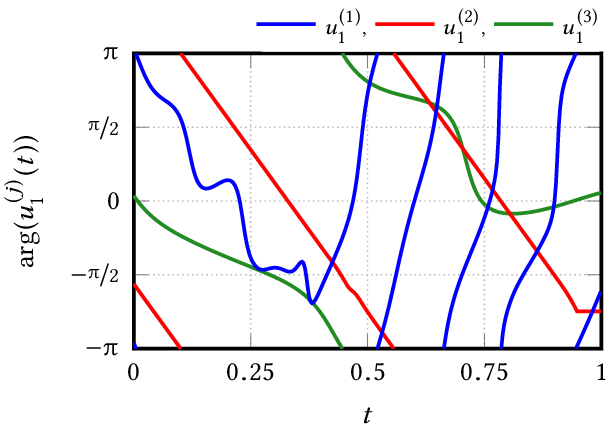}
  }
  \caption{Numerical construction of norm-minimizing solutions
    $u^\dagger \coloneqq u_1^{(j)}$ satisfying the source condition in
    Proposition~\ref{the:sour-cond-deauto} on the basis of an
    explicitly known source element $\phi^{(j)}$.  For the
    approximation of the integral operator
    $\Phi_{\mathcal I_k^*[\phi]}$ in \eqref{eq:approx-int-oper-phi}, a
    discretization with $N=10^{5}$ samples has been used.}
  \label{fig:source-element}
\end{figure}

For the numerical experiments, we choose the signals
$\phi^{(j)} \coloneqq \absn{\phi^{(j)}} \e^{\I\arg(\phi^{(j)})}$ with
\begin{align*}
  \absn{\phi^{(1)}(t)} 
  &\coloneqq 2.7 \, \e^{-\nicefrac{(t-0.27)^2}{0.025^2}} + 4 \,
    \e^{-\nicefrac{(t-1.56)^2}{0.024^2}} 
  \\[\fsmallskip]
  &\qquad+ 4 \,
    \e^{-\nicefrac{(t-0.346)^2}{0.022^2}} + 2 \,
    \e^{-\nicefrac{(t-1.146)^2}{0.024^2}}
  \\[\fskip]
  \arg(\phi^{(1)}(t))
  &\coloneqq 5 \cos(7.867 t) - 2.3 \sin(25.786 t),
  \\[20pt]
  \absn{\phi^{(2)}(t)}
  &\coloneqq 2.95 \, \Ind_{[0.415,0.459]}(t) + 6 \,
    \Ind_{[0.92,0.95]}(t)
  \\[\fskip]
  \arg(\phi^{(2)}(t))
  &\coloneqq 2 \cos(0.867 t) - 1.3 \sin(25.786 t),
  \\
  \shortintertext{and}
  \absn{\phi^{(3)}(t)}
  &\coloneqq 0.645 \sinc(9.855 \, (t-1.2)) + 0.434 \sinc(15.243 \,
    (t-0.42))
  \\[\fsmallskip]
  &\qquad+ 6.234 \sinc(0.143 \, (t-0.85))
  \\[\fskip]
  \arg(\phi^{(3)}(t))
  &\coloneqq 1.2 \cos(2.867t) - 2.3 \sin(4.786 \, (t-0.78)) +
    \e^{0.643 t}.
\end{align*}
Here the indicator function $\Ind_{[t_1,t_2]}(t)$ is $1$ if $t$ is
contained in the interval $[t_1,t_2]$ and $0$ else.  Rescaling the
source element with respect to the major eigenvalue
\raisebox{0pt}[0pt][0pt]{$\lambda_1^{(j)}$} by
\raisebox{0pt}[0pt][0pt]{$\nicefrac{1}{\lambda_1^{(j)}} \,
  \phi^{(j)}$},
we easily obtain norm-minimizing solutions
\raisebox{0pt}[0pt][0pt]{$u^\dagger \coloneqq u_1^{(j)}$} which
satisfy the required source condition
(\ref{eq:source-cond-deauto:norm-eigen}--\ref{eq:source-cond-deauto:oper}).
The source elements \raisebox{0pt}[0pt][0pt]{$\phi^{(j)}$} and the
related eigenfunctions \raisebox{0pt}[0pt][0pt]{$u_1^{(j)}$} are
presented in \autoref{fig:source-element}.  The results of the
simulations here look quite promising in the sense that the class of
functions $u^\dagger$ satisfying
(\ref{eq:source-cond-deauto:norm-eigen}--\ref{eq:source-cond-deauto:oper})
is rich on naturally occurring signals.

\subsection{Validation of the theoretical converence rate}
\label{sec:valid-theor-conv}

\begin{figure}\centering
  \subfloat[Magnitude of the approximation $u^*_j$.]{
    \includegraphics{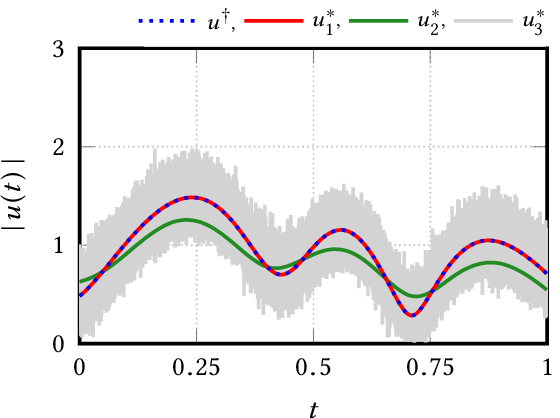}
  }
  \quad
  \subfloat[Phase of the approximation $u^*_j$.]{
    \includegraphics{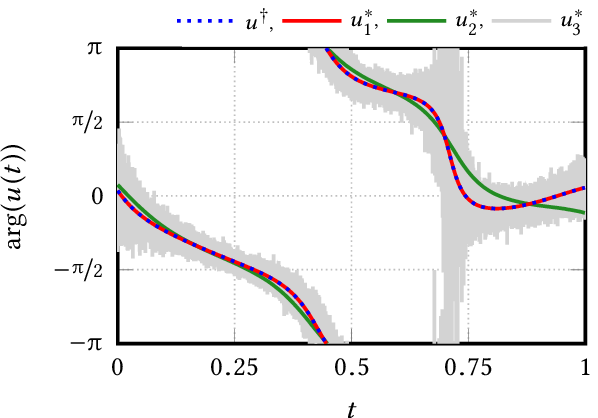}
  }
  \caption{Comparison between the norm-minimizing solution $u^\dagger$
    and the numerically reconstructed signals $u_\alpha^\delta$ with and
    without regularization.  The noise level for $u^*_1$ and $u^*_3$
    amounts to $\delta = 0.05 \pNormn{u^\dagger}$, and for $u^*_2$ to
    $\delta=100\pNormn{u^\dagger}$.  The regularization for $u^*_1$
    and $u^*_2$ corresponds to $\alpha = 100 \delta$.  The
    reconstruction $u^*_3$ is computed without regularization or with
    $\alpha=0$.}
  \label{fig:comp-wiht-without-regl}
\end{figure}

To verify the convergence rate numerically, we have to solve the
deautoconvolution problem for different noise levels. Referring to
\cite{Ger11} and \cite{GHB+14}, we here apply the \PGauss-\PNewton\ method
to the discretized problem with forward operator
\begin{equation*}
  \mathcal A_k[u]\bigl(\tfrac{n}{N}\bigr)
  \approx \frac{1}{N} \sum_{k=0}^N u \bigl( \tfrac{k}{N} \bigr) \, u
  \bigl( \tfrac{n-k}{N} \bigr).
\end{equation*}
In order to solve the occurring equation systems iteratively, we use
the conjugate gradient method with a suitable preconditioner and
exploit the \PToeplitz\ structure of the related system matrix, see
for instance \cite{RZR12}.  For the exact signal $u^\dagger$ arising
from the source element $\phi^{(3)}$, the numerical approximations
$u^*_j$ corresponding to the minimizer
\raisebox{0pt}[0pt][0pt]{$u_\alpha^\delta$} of the regularized and
non-regularized deautoconvolution problem are shown in
\autoref{fig:comp-wiht-without-regl}.  Besides the numerical
ill-posedness of the discretized deautoconvolution problem for minor
levels $\delta \coloneqq \pNormn{g^\delta - g^\dagger}$ of
\PGauss{ian} noise, we can here see the smoothing effect of the
applied \PTikhonov\ regularization.  Even for considerably higher
noise levels like $\delta = 100 \pNormn{u^\dagger}$, the
reconstruction covers the main features of the unknown norm-minimizing
solution $u^\dagger$.  The considered noise level $\delta$ here
depends on the norm $\pNormn{u^\dagger}$ of the true signal and not on
the norm $\pNormn{g^\dagger}$ of the exact data.  Depending on our
underlying numerical implementation, for the considered signal, the
noise level $\delta = 100 \pNormn{u^\dagger}$ corresponds to a
\PGauss{ian} noise, whose norm approximately equals
$\pNormn{g^\dagger}$.

\begin{figure}\centering
  \subfloat[\PBregman\ distance between $u_\alpha^\delta$ and
  $u^\dagger$.\label{fig:conv-rate-ana:a}]{
    \includegraphics{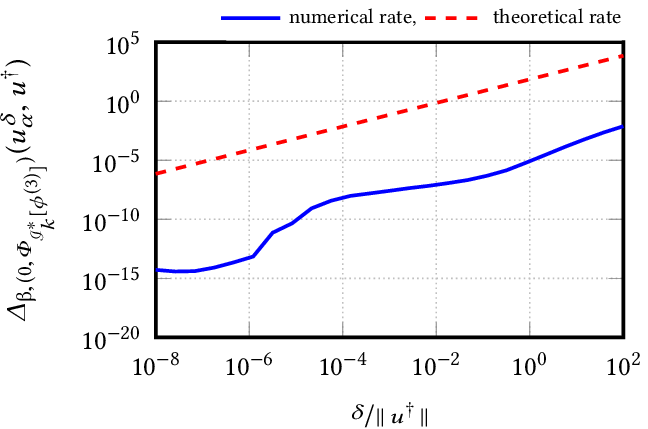}
  }
  \quad
  \subfloat[Distance of $u_\alpha^\delta$ to ray $E_1$ spanned by
  $u^\dagger$.\label{fig:conv-rate-ana:b}]{ 
    \includegraphics{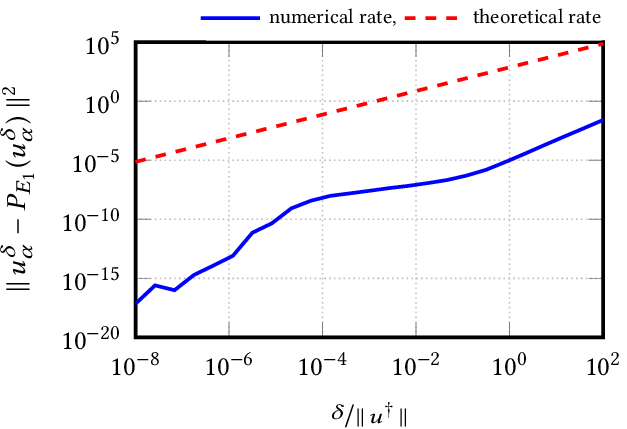}
  }
  \par
  \subfloat[{Discrepancy between $\mathcal A_k[u_\alpha^\delta]$ and
    $g^\delta$.\label{fig:conv-rate-ana:c}}]{ 
    \includegraphics{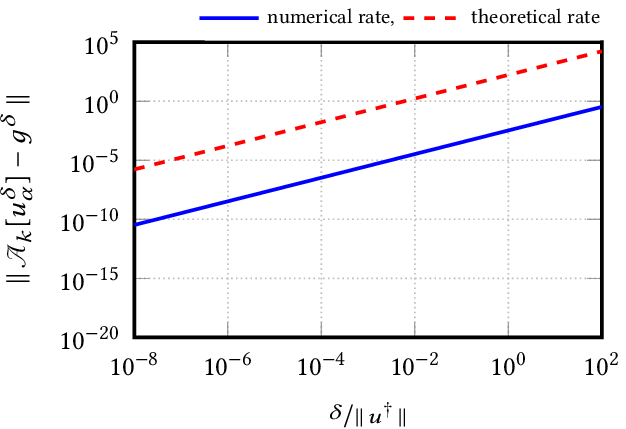}
  } \quad 
  \subfloat[Uncovered error part of $u_\alpha^\delta$ within
  the ray $E_1$ spanned by $u^\dagger$.\label{fig:conv-rate-ana:d}]{
    \includegraphics{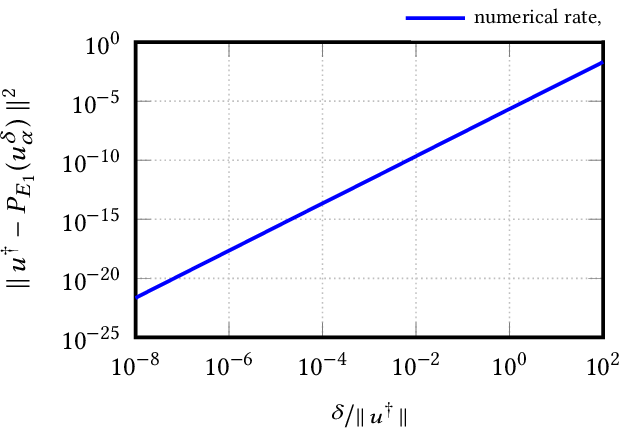}
  }
  \caption{Comparison of the theoretical convergence rates in
    Corollary~\ref{cor:conv-rate-deauto} and
    Theorem~\ref{the:brag-dist-deauto} with the numerical convergence
    rates for the source element $\phi^{(3)}$.}
  \label{fig:conv-rate-ana}
\end{figure}

Unfortunately, the result of the \PGauss-\PNewton\ method usually
strongly depends on the start value, which have to be a very accurate
approximation of the norm-minimizing solution $u^\dagger$ for very low
noise levels $\delta$.  For this reason, we extent the simulations for
the convergence rate analysis by choosing 25 randomly created start
values around the true solution $u^\dagger$ per noise level.  Since we
have constructed the norm-minimizing solution $u^\dagger$ and the
exact data $g^\dagger$ from a specific source element $\phi$, besides
the convergence rates in \thref{cor:conv-rate-deauto}, we have an
explicit upper bound for the \PBregman\ distance between the
regularized solution $u_\alpha^\delta$ and the norm-minimizing
solution $u^\dagger$ as well as for the discrepancy between the
forward operator $\mathcal A_k[u_\alpha^\delta]$ and the perturbed
data $g^\delta$.  The convergence rate analysis for the source element
$\phi^{(3)}$ is presented in \autoref{fig:conv-rate-ana}.
Additionally, \thref{the:brag-dist-deauto} yields an upper bound for
the distance
$\pNormn{P_{E_1^\perp}(u_\alpha^\delta)}^2 = \pNormn{u_\alpha^\delta -
  P_{E_1}(u_\alpha^\delta)}^2$
between \raisebox{0pt}[0pt][0pt]{$u_\alpha^\delta$} and the ray $E_1$
spanned by $u^\dagger$.  Here all three theoretical convergence rates
and upper bounds match with numerical results.

If we consider the discrepancy
$\pNormn{\mathcal A_k[u_\alpha^\delta] - g^\delta}$ in
\autoref{fig:conv-rate-ana}.\subref*{fig:conv-rate-ana:c} more
closely, we notice that the numerical and the theoretical rates
coincides except for a multiplicative constant.  Consequently, we
cannot hope to improve the theoretical rate of $\Landau(\delta)$.  The
\PBregman\ distance and the distance to the ray spanned by $u^\dagger$
have a completely different behaviour.  More precisely, we here have
regions where the convergence rate is faster and regions where the
convergence rate is slower than the theoretical rate of
$\Landau(\delta)$.  Especially for the distance to the ray, it seems
that the overall convergence rate is much faster than the theoretical
rate.  In some circumstances, our theoretical rate could thus be too
pessimistic.

In this instance of the deautoconvolution problem, we can observe that
the error
$\pNormn{P_{E_1}(u_\alpha^\delta - u^\dagger)}^2 =
\pNormn{u_\alpha^\delta - P_{E_1}(u_\alpha^\delta)}^2$
within the ray $E_1$, which is not covered by
\thref{cor:conv-rate-deauto}, numerically converges to zero
superlinearly with a rate of $\Landau (\delta^2)$.  The shown
numerical rate here strongly depends on the starting values of the
\PGauss-\PNewton\ method, which have been chosen in a small
neighbourhood around $u^\dagger$.  Choosing starting values around
$-u^\dagger$, we would observe the same convergence rate to
$-u^\dagger$.  In fact, the sequence $u_\alpha^\delta$ could be
composed of two subsequences -- one converging to $u^\dagger$ and the
other to $-u^\dagger$.

\section{Conclusion}
\label{sec:conclusion}

Starting from the question: how non-linear may a non-linear forward
operator be in order to extend the linear regularization theory, we
have introduced the classes of dilinear and diconvex mappings, which
corresponds to linear, bilinear, and quadratic inverse problems.
Exploiting the tensorial structure behind these mappings, we have
introduced two different concepts of generalized subgradients and
subdifferentials -- the dilinear and the representative
generalization.  We have shown that the classical subdifferential
calculus can be partly transferred to both new settings.  Although the
representative generalization yields stronger computation rules, the
related subdifferential unfortunately strongly depends on the
non-unique convex representative of the considered diconvex mapping.
Besides all differences, there exists several connections between the
dilinear and representative subdifferential.

On the basis of these preliminary considerations, we have examined the
class of dilinear inverse problems.  Analogously to linear inverse
problems, the regularizations in the dilinear setting are well posed,
stable, and consistent.  Using the injective tensor product, which is
nearly a predual space of the projective tensor product, as
topologizing family, we have seen that the required sequential weak$*$
(semi-)continuity of the forward operator and regularization term may
be inherited from the lifted versions.  Moreover, we have derived a
convergence rate analysis very similar to the linear setting under
similar assumptions and requirements.  This enables us to give
explicit upper bounds for the discrepancy and the \PBregman\ distance
between the solution of the regularized problem and the $R$-minimizing
solution.

In a last step, we have applied the developed theory to the
deautoconvolution problem that appears in spectroscopy, optics, and
stochastics.  Although the requirements of the non-linear
regularization theory are not fulfilled, our novel approach yields
convergence rates based on a suitable range source condition and the
dilinear \PBregman\ distance.  Depending on the source element, the
\PBregman\ distance is here surprisingly strong.  In the best case,
the solutions of the regularized problems converge strongly to the ray
spanned by the true signal, which is the best possible rate with
respect to the ambiguities of the problem.  Using numerical
experiments, we have considered different source elements and the
corresponding norm-minimizing solutions, which shows that there exists
signals satisfying the required source-wise representation.  Finally,
we have observed the established error bounds in the numerical
simulations.

\bigskip\noindent
\parbox{\linewidth}{\footnotesize {\sffamily\bfseries
    Acknowledgements.} \quad We gratefully acknowledges the funding of
  this work by the Austrian Science Fund (FWF) within the project
  P28858.  The Institute of Mathematics and Scientific Computing of
  the University of Graz, with which the authors are affiliated, is a
  member of NAWI Graz (http://www.nawigraz.at/).}

\newcommand{\etalchar}[1]{$^{#1}$}

\end{document}